\newcommand\isep{\parskip-2mm  \itemsep3mm}

\documentclass[11pt,twoside,leqno]{aomamlt2e} 
\pageno{367}
\received{January 20, 2004}

\allowdisplaybreaks
\DeclareMathAlphabet{\curly}{U}{rsfs}{m}{n}

\newtheorem{thm}{Theorem}
\newtheorem{cor}{Corollary}
\newtheorem{lem}{Lemma}[section]
\newtheorem{conj}{Conjecture}

\renewcommand{\theequation}{\thesection.\arabic{equation}} 

\newcommand{\nts}{\negthickspace}
\newcommand{\ZZ}{{\mathbb Z}}

\newcommand{\RR}{{\mathbb R}}

\newcommand{\bb}{{\mathbf b}}

\newcommand{\LL}{\curly L}

\newcommand{\KK}{\curly K}
\newcommand{\BB}{\curly B}

\newcommand{\GG}{\curly G}

\newcommand{\TT}{\curly T}
\newcommand{\PP}{\curly P}

\newcommand{\Vol}{\operatorname{Vol}}   

\newcommand{\sg}{\ensuremath{\sigma}}

\newcommand{\g}{\ensuremath{\gamma}}
\newcommand{\del}{\ensuremath{\delta}}
\newcommand{\lam}{\ensuremath{\lambda}}
\newcommand{\eps}{\ensuremath{\varepsilon}}

\newcommand{\bg}{\ensuremath{\mathbf{g}}}

\newcommand{\xx}{\ensuremath{\mathbf{x}}}
\newcommand{\bj}{\ensuremath{\mathbf{j}}}

\newcommand{\bx}{{\ensuremath{\boldsymbol{\xi}}}}

\newcommand{\fl}[1]{{\ensuremath{\left\lfloor {#1} \right\rfloor}}}
\newcommand{\cl}[1]{{\ensuremath{\left\lceil #1 \right\rceil}}}
\newcommand{\pfrac}[2]{{\left(\frac{#1}{#2}\right)}}
\newcommand{\lcm}{\text{lcm}}
\newcommand{\ds}{\displaystyle}

\renewcommand{\b}{\ensuremath{\beta}}
\renewcommand{\a}{\ensuremath{\alpha}}
\renewcommand{\SS}{\curly S}
\renewcommand{\AA}{\curly A}
\renewcommand{\(}{\left(}
\renewcommand{\)}{\right)}
\renewcommand{\th}{\ensuremath{\theta}}

\newif\ifdraft

\begin{document}
\currannalsline{168}{2008} 

\title{The distribution of integers with a divisor \\
in a given interval}
\shorttitle{Integers with a divisor in an interval}

\author{Kevin Ford}

\institution{University
of Illinois at Urbana-Champaign, Urbana, IL\\
\email{ford@math.uiuc.edu}}
 

\centerline{\bf Abstract} 
\vskip12pt We determine the order of magnitude of $H(x,y,z)$, the
number of integers $n\le x$ having a divisor in $(y,z]$, for all
$x,y$ and $z$.  We also study
$H_r(x,y,z)$, the number of integers $n\le x$ having exactly $r$ divisors in
$(y,z]$.  When  \hbox{$r=1$} we establish the order of magnitude of $H_1(x,y,z)$
for all $x,y,z$ satisfying $z\le x^{1/2-\eps}$.  For every $r\ge 2$,
$C>1$ and $\eps>0$, we
determine the order of magnitude of $H_r(x,y,z)$ uniformly for
$y$ large and
$y+y/(\log y)^{\log 4 -1 - \eps} \le z \le \min(y^{C},x^{1/2-\eps})$.
As a consequence of these bounds, we settle a 1960 conjecture of Erd\H os
and some conjectures of Tenenbaum.
One key element of the proofs is a new result on the distribution of
uniform order statistics.
 
%
%
\vskip12pt
\centerline{\bf Contents}
\def\sni#1{\vskip-1pt
\noindent{\phantom{0}#1} \hskip5pt}
\def\sno#1{\vskip-1pt\noindent{#1} \hskip5pt}
\vskip4pt

\sni{1.} Introduction

\sni{2.} Preliminary lemmas

\sni{3.} Upper bounds outline

\sni{4.} Lower bounds outline

\sni{5.} Proof of Theorems \ref{thm1}, \ref{thm:shorts}, \ref{thm:squarefree},
\ref{H1bounds}, and \ref{Hrbounds}

\sni{6.} Initial sums over $L(a;\ensuremath {\sigma })$ and $L_s(a;\ensuremath
  {\sigma })$

\sni{7.} Upper bounds in terms of $S^*(t;\ensuremath {\sigma })$

\sni{8.} Upper bounds: reduction to an integral

\sni{9.} Lower bounds: isolated divisors

\sno{10.} Lower bounds: reduction to a volume

\sno{11.} Uniform order statistics

\sno{12.} The lower bound volume

\sno{13.} The upper bound integral

\sno{14.} Divisors of shifted primes

\vskip-1pt\noindent References

\renewcommand{\theenumi}{\roman{enumi}}

\vglue-25pt\phantom{up}
\section{Introduction}\label{sec:intro}
%
%
%

For $0<y<z$, let
$\tau(n;y,z)$ be the number of divisors $d$ of $n$ which
satisfy $y<d\le z$.  Our focus in this paper is to estimate 
$H(x,y,z)$, the number of positive integers $n\le x$ with
$\tau(n;y,z) > 0$, and $H_r(x,y,z)$, the number of $n\le x$ with
$\tau(n;y,z)=r$.  
By inclusion-exclusion,
$$
H(x,y,z) = \sum_{k\ge 1} (-1)^{k-1} \sum_{y<d_1<\cdots<d_k\le z} \fl{
\frac{x}{\text{lcm}[d_1,\cdots,d_k]}},
$$
but this is not useful for estimating $H(x,y,z)$ unless $z-y$ is small.
With $y$ and $z$ fixed, however, this formula implies that
the set of positive integers having at least one
divisor in $(y,z]$ has an \emph{asymptotic density}, i.e.\ the limit
$$
\eps(y,z) = \lim_{x\to\infty}  \frac{H(x,y,z)}{x}
$$
exists.
Similarly, the exact formula
$$
H_r(x,y,z) = \sum_{k\ge r} (-1)^{k-r} \binom{k}{r}  
\sum_{y<d_1<\cdots<d_k\le z} \fl{\frac{x}{\text{lcm}[d_1,\cdots,d_k]}}
$$
implies the existence of
$$
\eps_r(y,z) = \lim_{x\to\infty}  \frac{H_r(x,y,z)}{x}
$$
for every fixed pair $y,z$.

\Subsec{Bounds for $H(x,y,z)$}
Besicovitch initiated the study of such quantities in 1934, proving in
\cite{Bes} that
\begin{equation}\label{Bes}
\liminf_{y\to \infty} \eps(y,2y) = 0,
\end{equation} 
and using \eqref{Bes} to construct 
an infinite set $\AA$ of positive integers
such that its set of multiples $\BB(\AA)=\{am: a\in \AA, m\ge 1\}$ does
not possess asymptotic density.
Erd\H os in 1935 \cite{Erdos35} showed $\ds
\lim_{y\to\infty} \eps(y,2y)=0$ and 
in 1960 \cite{Erdos60} gave the further refinement
(see also Tenenbaum \cite{Ten77})
$$
\eps(y,2y) = (\log y)^{-\del+o(1)} \quad (y\to \infty),
$$
where
$$
\del = 1 - \frac{1+\log\log 2}{\log 2} = 0.086071\ldots.
$$

Prior to the 1980s, a few other special cases were studied.
In 1936, Erd\H os \cite{Erdos36} established
$$
\lim_{y\to \infty} \eps(y,y^{1+u}) = 0, 
$$
provided that $u=u(y)\to 0$ as $y\to \infty$.
In the late 1970s, Tenenbaum (\cite{TenII}, \cite{TenIII})
showed that
$$
h(u,t) = \lim_{x\to\infty} \frac{H(x,x^{(1-u)/t},x^{1/t})}{x}
$$
exists for $0\le u\le 1$, $t\ge 1$ and gave bounds on $h(u,t)$.

%
%

Motivated by a growing collection of applications for such bounds,
Tenenbaum in the early 1980s turned to the problem of 
bounding $H(x,y,z)$ for \emph{all} $x,y,z$.  In the seminal work \cite{Ten84}
he established reasonably
sharp upper and lower bounds for $H(x,y,z)$ which we list below
(paper \cite{Ten83} announces these results and gives a history of
previous bounds for $H(x,y,z)$;  Hall and Tenenbaum's book 
\emph{Divisors} \cite{Divisors} gives a simpler proof of Tenenbaum's theorem).
We require some additional notation.
For a given pair $(y,z)$ with $4\le y < z$, we define $\eta, u, \b, \xi$ by
\begin{equation}\label{etabetaxi}
z = e^{\eta} y=y^{1+u}, \quad \eta = (\log y)^{-\b}, \quad \b = \log 4 - 1 +
\frac{\xi}{\sqrt{\log\log y}}.
\end{equation} 
Tenenbaum defines $\eta$ by $z=y(1+\eta)$, which is asymptotic to our
$\eta$ when $z-y=o(y)$.  The definition in \eqref{etabetaxi} plays
a natural role in the arguments even when $z-y$ is large.
For smaller $z$, we also need the function
\begin{equation}\label{Gdef}
G(\b) = \begin{cases} 
\frac{1+\b}{\log 2} \log \pfrac{1+\b}{e\log 2} + 1 & 0\le \b \le \log 4 -1\\
\b & \log 4-1 \le \b .\end{cases}
\end{equation} 
When $x$ and $y$ are fixed, Tenenbaum discovered that
$H(x,y,z)$ undergoes a change of behavior in the vicinity of 
$$
z=z_0(y) := y\exp\{(\log y)^{1-\log 4} \} \approx y+y/(\log y)^{\log 4-1},
$$
in the vicinity of $z=2y$ and in the vicinity of $z=y^2$.  

\demo{\sc Theorem T1 {\rm (Tenenbaum \cite{Ten84})}}
(i) {\it Suppose $y\to\infty${\rm ,} $z-y\to \infty${\rm ,} $z\le \sqrt{x}$ and
$\xi\to\infty$.  Then
$$
H(x,y,z) \sim \eta x.
$$

{\rm (ii)} Suppose $2\le y < z \le \min(2y,\sqrt{x})$ and $\xi$ is bounded above.
Then
$$
\frac{x}{(\log y)^{G(\b)} Z(\log y)} \ll H(x,y,z) \ll \frac{x}
{(\log y)^{G(\b)} \max(1,-\xi)}.
$$
Here $Z(v)=\exp\{ c\sqrt{\log(100v) \log\log (100v)} \}$ and
$c$ is some positive constant.

\vskip4pt
{\rm (iii)} Suppose $4\le 2y \le z \le \min(y^{3/2},\sqrt{x})$.  Then
$$
\frac{x u^{\del}}{Z(1/u)}  \ll H(x,y,z) \ll \frac{x u^\del \log\log (3/u)}
{\sqrt{\log(2/u)}}.
$$
Moreover{\rm ,} the term $\log\log(3/u)$ on the right may be omitted if
$z\le By$ for some $B>2${\rm ,} the constant implied by $\ll$ depending on $B$.

\vskip4pt {\rm (iv)} If $2\le y\le z\le x${\rm ,} then
$$
H(x,y,z) = x \( 1 + O\pfrac{\log y}{\log z} \).
$$
}

\demo{Remark}  Since
$$
\sum_{n\le x} \tau(n,y,z) = \sum_{y<d\le z} \fl{\frac{x}{d}} \sim \eta
x \qquad (z-y\to \infty),
$$
in the range of $x,y,z$ given in (i) of Theorem T1,
most $n$ with a divisor in $(y,z]$ have only one such divisor.  By (iv),
when $\frac{\log z}{\log y} \to \infty$, almost all integers have a divisor
in $(y,z]$.  

In 1991, Hall and Tenenbaum \cite{HT91} established the order
of $H(x,y,z)$ in the vicinity of the ``threshhold'' $z=z_0(y)$.  
Specifically, they showed that if $3\le y+1 \le z \le \sqrt{x}$,
 $c>0$ is fixed and $\xi \ge -c (\log\log y)^{1/6}$, then 
$$
H(x,y,z) \asymp \frac{x}{(\log y)^{G(\b)} \max(1,-\xi)},
$$
thus showing that the upper bound given by (ii) of Theorem T1 is the
true order in this range.  
In fact the argument in \cite{HT91} implies that
$$
H_1(x,y,z) \asymp H(x,y,z)
$$
in this range of $x,y,z$.  Specifically, 
Hall and Tenenbaum use a lower estimate
$$
H(x,y,z) \ge \sum_{\substack{n\le x \\ n\in \curly N}} \tau(n,y,z) (2
- \tau(n,y,z))
$$
for a certain set $\curly N$, and clearly the right side is also a
lower bound for $H_1(x,y,z)$.
Later, in a slightly more restricted range, 
Hall  (\cite{Hallbook}, Ch. 7) proved an
asymptotic formula for $H(x,y,z)$ which extends the asymptotic
formula of part (i) of Theorem T1.  Richard Hall has kindly pointed
out an error in the stated range of validity of this asymptotic in
\cite{Hallbook}, which we correct below (in \cite{Hallbook}, the 
range is stated as $\xi \ge -c(\log\log y)^{1/6}$).

\demo{\sc Theorem H {\rm (Hall \cite[Th.~7.9]{Hallbook})}}
{\it Uniformly for $z \le x^{1/\log\log x}$ and for
$\xi \ge - o(\log\log y)^{1/6}${\rm ,}
$$
\frac{H(x,y,z)}{x} = (F(\xi)+O(E(\xi,y))) (\log y)^{-\b},
$$
where
$$
F(\xi) = \frac{1}{\sqrt{\pi}} \int_{-\infty}^{\xi/\log 4} e^{-u^2}\, du
$$
and}
$$
E(\xi,y) = \begin{cases} \dfrac{\xi^2+\log\log\log y}{\sqrt{\log\log y}}
e^{-\xi^2/\log^2 4}, & \xi \le 0 \\ \quad \\
\dfrac{\xi+\log\log\log y}{\sqrt{\log\log y}}, & \xi>0.
\end{cases}
$$
\Enddemo

Note that
$$
F(\xi) (\log y)^{-\b} \asymp \frac{1}{(\log y)^{G(\b)} \max(1,-\xi)}
$$
in Theorem {\rm H}.

We now determine the exact order of $H(x,y,z)$ for all $x,y,z$.
Constants implied by $O$, $\ll$ and $\asymp$ 
are absolute unless otherwise noted, e.g.\ by a subscript.  
The notation $f \asymp g$ means $f\ll g$ and $g\ll f$.
Variables $c_1, c_2, \ldots$ will denote certain specific constants,
$y_0$ is a sufficiently large real number, while $y_0(\cdot)$ will
denote a large constant depending only on the parameters given, e.g.\ $y_0(r,c,c')$, and the meaning may change from statement to statement.
Lastly, $\fl{x}$
denotes the largest integer $\le x$.

%
%

\begin{thm}\label{thm1}  Suppose $1\le y\le z\le x$.  Then{\rm ,}
\begin{enumerate}
\item $H(x,y,z) = 0$ if $z< \fl{y}+1${\rm ;}
\item $H(x,y,z) = \fl{x/(\fl{y}+1)}$ if $\fl{y}+1\le z<y+1${\rm ;}
\item $H(x,y,z) \asymp 1$ if $z\ge y+1$ and $x\le 100000${\rm ;}
\item $H(x,y,z) \asymp x$ if $x\ge 100000${\rm ,} $1\le y\le 100$ and $z\ge
  y+1${\rm ;}
\item If $x>100000${\rm ,} $100 \le y \le z-1$ and $y\le \sqrt{x}${\rm ,} 
$$
\frac{H(x,y,z)}{x} \asymp \begin{cases} \log(z/y)=\eta & y+1 \le z \le
  z_0(y) \\  \\
\dfrac{\b}{\max(1,-\xi) (\log y)^{G(\b)}} & z_0(y) \le z \le 2y \\ \\ 
u^\del (\log \tfrac{2}{u})^{-3/2} & 2y \le z \le y^2 \\ \\
1 & z \ge y^{2}.
\end{cases}
$$
\item If $x>100000${\rm ,} $\sqrt{x} < y < z \le x$ and $z \ge y+1${\rm ,} then
$$
H(x,y,z) \asymp \begin{cases} H\(x, \frac{x}{z}, \frac{x}{y}\) &
  \frac{x}{y} \ge \frac{x}{z} + 1 \\ \eta x & \text{ otherwise}.
\end{cases}
$$
\end{enumerate}
\end{thm}

\begin{cor}\label{cor1}
Suppose $x_1,y_1,z_1,x_2,y_2,z_2$ are real numbers with
$1 \le y_i < z_i \le x_i$ $(i=1,2),$ $z_i \ge y_i+1$ $(i=1,2),$
$\log(z_1/y_1) \asymp \log(z_2/y_2)${\rm ,} $\log y_1 \asymp \log y_2$ and
$\log(x_1/z_1) \asymp \log(x_2/z_2)$.  Then
$$
\frac{H(x_1,y_1,z_1)}{x_1} \asymp \frac{H(x_2,y_2,z_2)}{x_2}.
$$
\end{cor}

\begin{cor}\label{Besi}
If $c>1$ and $\frac{1}{c-1}\le y\le x/c${\rm ,} then
$$
H(x,y,cy) \asymp_c \frac{x}{(\log Y)^{\del} (\log\log Y)^{3/2}} \qquad
(Y=\min(y,x/y)+3)
$$
and
$$
\eps(y,cy) \asymp_c \frac{1}{(\log y)^{\del} (\log\log y)^{3/2}}.
$$
\end{cor}

Items (i)--(iv) of Theorem \ref{thm1} are trivial.
The first and fourth part of item (v) are already known (cf.\ the
papers of Tenenbaum \cite{Ten84} and Hall and Tenenbaum \cite{HT91}
mentioned above).  Item (vi) essentially follows from (v) by observing
that  $d|n$ if and only if  $(n/d)|n$.  However, proving (vi)
requires a version of (v) where $n$ is restricted to a short interval,
which we record below.  The range of $\Delta$ can be considerably improved,
but the given range suffices for the application to
Theorem \ref{thm1} (vi).

%
%

\begin{thm}\label{thm:shorts}
For $y_0 \le y \le \sqrt{x}${\rm ,} $z\ge y+1$  
and $\frac{x}{\log^{10} z} \le \Delta \le x${\rm ,} 
$$
H(x,y,z) - H(x-\Delta,y,z) \asymp \frac{\Delta}{x} H(x,y,z).
$$
\end{thm}

Motivated by an application to gaps in the Farey series, we also
record an analogous result for $H^*(x,y,z)$, the number of
\emph{squarefree} numbers $n\le x$ with $\tau(n,y,z)\ge 1$.

\begin{thm}\label{thm:squarefree}
Suppose $y_0 \le y \le \sqrt{x}$, $y+1\le z\le x$ and 
$\frac{x}{\log y} \le \Delta \le x$.  If $z\ge y + K y^{1/5} \log y${\rm ,}
where $K$ is a large absolute constant{\rm ,} then
$$
H^*(x,y,z)-H^*(x-\Delta,y,z) \asymp \frac{\Delta}{x} H(x,y,z).
$$
If $y+ (\log y)^{2/3} \le z \le y+K y^{1/5} \log y${\rm ,} $g>0$ and 
there are $\ge g(z-y)$ square-free numbers in $(y,z]${\rm ,} then 
$$
H^*(x,y,z)-H^*(x-\Delta,y,z) \asymp_g \frac{\Delta}{x} H(x,y,z).
$$
\end{thm}

To obtain good lower bounds on $H^*(x,y,z)$, it is important that 
$(y,z]$ contain many squarefree integers.  In the
extreme case where $(y,z]$ contains no squarefree integers, clearly
$H^*(x,y,z)=0$.  A theorem of Filaseta and Trifonov \cite{FiTr}
  implies that there are  $\ge\frac12(z-y)$ squarefree numbers in $(y,z]$ 
if $z\ge y + K y^{1/5}\log y$, and this is the best result known of this kind.

\demo{Some applications}
Most of the following applications depend on the distribution
of integers with $\tau(n,y,z)\ge 1$ when $z \asymp y$.  See also
Chapter 2 of
\cite{Divisors} for further discussion of these and other applications.

\vskip2pt
1. Distinct products in a multiplication table, a problem of Erd\H os
from 1955 (\cite{Erdos55}, \cite{Erdos60}).   Let
$A(x)$ be the number of positive integers $n\le x$ which can be written
as $n=m_1 m_2$ with each $m_i \le \sqrt{x}$. 

\begin{cor}
We have
$$
A(x) \asymp \frac{x}{(\log x)^{\del} (\log\log x)^{3/2}}.
$$
\end{cor}

\Proof
Apply Theorem \ref{thm1} and the inequalities
\vskip12pt
\hfill $
\displaystyle{H\( \frac{x}{4}, \frac{\sqrt{x}}{4}, \frac{\sqrt{x}}{2} \) \le A(x) \le
\sum_{k\ge 0} H\(\frac{x}{2^k}, \frac{\sqrt{x}}{2^{k+1}}, 
\frac{\sqrt{x}}{2^k} \). }
$ 
\Endproof\vskip12pt

2.  Distribution of Farey gaps (Cobeli, Ford, Zaharescu \cite{CFZ}).

\begin{cor}
Let $( \frac{0}{1}, \frac{1}{Q}, \dots, \frac{Q-1}{Q}, \frac{1}{1} 
)$ denote the sequence of Farey fractions of order $Q${\rm ,} and let
$N(Q)$ denote the number of distinct gaps between successive terms of
the sequence.  Then
$$
N(Q) \asymp \frac{Q^2}{(\log Q)^{\del} (\log\log Q)^{3/2}}.
$$
\end{cor}

\Proof
The distinct gaps are precisely those products $qq'$ with $1 \le q$,\break $q'\le Q$,
$(q,q')=1$ and $q+q' > Q$.  Thus
$$
H^*(\tfrac{9}{25} Q^2, \tfrac{Q}{2}, \tfrac{3Q}{5}) - 
H^*(\tfrac{3}{10} Q^2, \tfrac{Q}{2}, \tfrac{3Q}{5}) \le N(Q) \le H(Q^2,Q/2,
Q),
$$
and the corollary follows from Theorems \ref{thm1} and \ref{thm:squarefree}.
\Endproof\vskip4pt

3.  Divisor functions.  Erd\H os introduced (\cite{ET81}, \cite{ET83} and 
\S 4.6 of \cite{Divisors}) the function
$$
\tau^+(n) = | \{ k\in \ZZ : \tau(n,2^k,2^{k+1}) \ge 1 \} |.
$$

\begin{cor}\label{tauplus}
For $x\ge 3${\rm ,} 
$$
\frac{1}{x}\sum_{n\le x} \tau^+(n) \asymp \frac{(\log x)^{1-\del}}
{(\log\log x)^{3/2}}.
$$
\end{cor}

\Proof
This follows directly from Theorem \ref{thm1} and
\vskip12pt
\hfill 
$\displaystyle{\sum_{n\le x} \tau^+(n) = \sum_k H(x,2^k,2^{k+1}).}
$ 
\Endproof\vskip12pt

Tenenbaum \cite{Ten76} defines $\rho_1(n)$ to be the largest divisor 
$d$ of $n$ satisfying $d\le \sqrt{n}$.

\begin{cor}\label{rho1}
We have
$$
\sum_{n\le x} \rho_1(n) \asymp \frac{x^{3/2}}{(\log x)^{\del} (\log\log x)
^{3/2}}.
$$
\end{cor}

\Proof
Suppose $x/4^l < n\le x/4^{l-1}$.  Since $\rho_1(n)$ lies in
$(\sqrt{x} 2^{-k}, \sqrt{x} 2^{-k+1}]$ for some integer $k\ge l$,
\begin{align*}
\frac{\sqrt{x}}{4} \( H\( x, \tfrac{\sqrt{x}}{4}, \tfrac{\sqrt{x}}{2} \)
-  H\( \tfrac{x}{4}, \tfrac{\sqrt{x}}{4}, \tfrac{\sqrt{x}}{2} \) \) &\le 
\sum_{n\le x} \rho_1(n) \\
&\le \sum_{l\ge 1} \sum_{k\ge l} \frac{\sqrt{x}}{2^{k-1}}
H\( \tfrac{x}{4^{l-1}}, \tfrac{\sqrt{x}}{2^{k}}, \tfrac{\sqrt{x}}{2^{k-1}} \)
\end{align*}
and the corollary follows from Theorem \ref{thm1}.
\Endproof

4.  Density of unions of residue classes.  Given moduli $m_1,\dots,m_k$,
let $\delta_0(m_1,\dots,m_k)$ be the minimum, over all possible residue
classes $a_1 \bmod m_1$, $\dots, a_k \bmod m_k$, of the density 
of integers which lie in at least one of the classes.  By a theorem of
Rogers (see \cite[p.~242--244]{HRo}), the minimum is achieved by taking
$a_1=\cdots=a_k=0$ and thus $\delta_0(m_1,\dots,m_k)$ is the density of
integers possessing a divisor among the numbers $m_1,\dots,m_k$.
When $m_1,\dots,m_k$ consist of the integers in an interval $(y,z]$, then
$\delta_0(m_1,\dots,m_k)=\eps(y,z)$.

\vskip8pt

5.  Bounds for $H(x,y,z)$ were used in recent work of Heath-Brown 
\cite{HB} on the validity of the Hasse principle for pairs of
quadratic forms.

\vskip8pt

6.  Bounds on $H(x,y,z)$ are central to the study of the function
$$
\max \{ |a-b| : 1\le a,b \le n-1, ab\equiv 1 \pmod{n} \}
$$
in \cite{FKSY}.

%
\Subsec{Bounds for $H_r(x,y,z)$}
In the paper \cite{Erdos60}, Erd\H os made the following 
conjecture:\footnote{Erd\H os 
also mentioned this conjecture in some of his books on unsolved
problems, e.g.\ \cite{EG}, and he wrote it in the Problem Book (page 2)
of the Mathematisches Forschungsinstitut Oberwolfach.}

\begin{conj}[Erd\H os \cite{Erdos60}]\label{erdosconj}
$$
\lim_{y\to \infty} \frac{\eps_1(y,2y)}{\eps(y,2y)} = 0.
$$
\end{conj}

This can be interpreted as the assertion that the conditional probability that 
a random integer has exactly 1 divisor in $(y,2y]$ given that it
has at least one divisor in $(y,2y]$, tends to zero as $y\to\infty$.

In 1987, Tenenbaum \cite{Ten87} gave general bounds on $H_r(x,y,z)$, which
are
of similar strength to his bounds on $H(x,y,z)$ (Theorem T1) when
$z\le 2y$.

\demo{\sc Theorem T2 {\rm (Tenenbaum \cite{Ten87})}}
{\it Fix $r\ge 1${\rm ,} $c>0$.

{\rm (i)} If $y\to \infty${\rm ,} $z-y\to \infty${\rm ,} and $\xi\to\infty${\rm ,}
then 
$$
\frac{H_r(x,y,z)}{H(x,y,z)} \to \begin{cases} 1 & r=1 \\ 0 & r\ge 2
\end{cases}.
$$

\vskip2pt {\rm (ii)} If $y\ge y_0(r)${\rm ,} $z_0(y) \le z \le \min(2y,x^{1/(r+1)-c})${\rm ,} then
$$
\frac{1}{Z(\log y)} \ll_{r,c} \frac{H_r(x,y,z)}{H(x,y,z)} \le 1.
$$
\vskip2pt

{\rm (iii)} If $y_0(r) \le 2y \le z \le \min(y^{3/2},x^{1/(r+1)-c})${\rm ,}
$$
\frac{1}{\log(z/y) Z(\log y)} \ll_{r,c} \frac{H_r(x,y,z)}{H(x,y,z)} \ll_{r}
\frac{Z(\log y)}{(\log (z/y))^{\del}}.
$$

\vskip2pt {\rm (iv)} If $y\ge y_0(r)${\rm ,} $y^{3/2} \le z \le x^{1/2}${\rm ,} then}
$$
\frac{\log\pfrac{\log z}{\log y}}{\log z} \ll_r \frac{H_r(x,y,z)}{H(x,y,z)} 
\ll_r \frac{(\log y)^{1-\del} (\log\log z)^{2r+1}}{\log z}.
$$

\demo{Remarks} In \cite{Ten87}, (ii) and (iii) above are stated with
$c=0$, but the proofs of the lower bounds require $c$ to be positive.
The construction of $n$ with $\tau(n,y,z)=r$ on p. 177 of \cite{Ten87}
requires $z^{\frac{1}{r+3}+r+1} \le x$, but the proof can be modified
to work for $z\le x^{\frac{1}{r+1}-c}$ for any fixed $c>0$.
\Enddemo

Based on the strength of the bounds in (ii) and (iii)  above,
Tenenbaum made two conjectures.  In particular, he asserted that
Conjecture \ref{erdosconj} is false.

\begin{conj}[Tenenbaum \cite{Ten87}]\label{tenconj1}
For every $r\ge 1${\rm ,} $c>0${\rm ,} and $c'>0${\rm ,} if $\xi\to -\infty$ as $y\to\infty${\rm ,}
$y\le x^{1/2-c'}$ and $z\le cy${\rm ,} then
$$
H_r(x,y,z) \gg_{r,c,c'} H(x,y,z).
$$
\end{conj}

\begin{conj}[Tenenbaum \cite{Ten87}]\label{tenconj2}
If $c>0$ is fixed{\rm ,} $y\le x^{1/2-c}${\rm ,} $r\ge 1$ and $z/y\to \infty${\rm ,} then
$$
H_r(x,y,z) = o(H(x,y,z)).
$$
\end{conj}

Using the methods used to prove Theorem \ref{thm1} plus
some additional arguments, we shall prove much stronger bounds
on $H_r(x,y,z)$ which will settle these three conjectures
(except Conjecture \ref{tenconj1} when $z$ is near $z_0(y)$).
When $z\ge 2y $, the order of $H_r(x,y,z)$ depends on $\nu(r)$, the 
exponent of the largest power of 2 dividing $r$ (i.e.\ $2^{\nu(r)} \| r$). 

\begin{thm}\label{H1bounds}
Suppose that $c>0${\rm ,} $y_0(c) \le y${\rm ,}
$y+1 \le z \le x^{5/8}$ and $yz\le x^{1-c}$.  Then
\begin{equation}\label{P1smallz}
\frac{H_1(x,y,z)}{H(x,y,z)} \asymp_{c} \frac{\log\log (z/y+10)}
{\log (z/y+10)}.
\end{equation} 
\end{thm}

\begin{thm}\label{Hrbounds}
Suppose that $r\ge 2${\rm ,} $c>0${\rm ,} $y_0(r,c) \le y${\rm ,}
$z \le x^{5/8}$ and $yz\le x^{1-c}$. 
If $z_0(y) \le z\le 10y${\rm ,} then
\begin{equation}\label{Prsmallz}
\frac{\max(1,-\xi)}{\sqrt{\log\log y}} \ll_{r,c} 
\frac{H_r(x,y,z)}{H(x,y,z)} \le 1.
\end{equation} 
When $C>1$ is fixed and $10y \le z \le y^{C}${\rm ,}  
\begin{equation}\label{Prlargez}
\frac{H_r(x,y,z)}{H(x,y,z)} \asymp_{r,c,C} 
\frac{(\log\log (z/y))^{\nu(r)+1}}{\log (z/y)}.
\end{equation} 
When $y\ge y_0(r)$ and $y^2 \le z \le x^{5/8}${\rm ,} then
\begin{equation}\label{Prhugez}
\frac{H_r(x,y,z)}{H(x,y,z)} \gg_r \frac{(\log\log y)^{\nu(r)+1}}{\log z}.
\end{equation} 
\end{thm}

\begin{cor}\label{epscor}
For every $\lam>1$ and $r\ge 1${\rm ,}
$$
\frac{\eps_r(y,\lam y)}{\eps(y,\lam y)} \gg_{r,\lam} 1.
$$
while for each $r\ge 1${\rm ,} if $z/y\to \infty$ then
$$
\frac{\eps_r(y,z)}{\eps(y,z)}  \to 0.
$$
\end{cor}
\vskip8pt

In particular, Conjecture \ref{erdosconj} is false, Conjecture
\ref{tenconj2} is true, and Conjecture~\ref{tenconj1} is true
provided $z\ge y + y/(\log y)^{\log 4-1-b}$ for a fixed $b>0$.

The upper bounds in Theorems \ref{H1bounds} and \ref{Hrbounds} are 
proved in the wider range $y\le \sqrt{x}, z\le x^{5/8}$.  
The conclusions of the two theorems, however, are not true when
$yz \approx x$.  This is a
consequence of $d|n$ implying $\frac{n}{d}|n$, which shows for
example that $\tau(n,y,n/y)$ is odd only if $n$ is a square
or $y|n$.  For another example, while $H_1(x,x^{1/4},x^{3/5}) \asymp x 
\frac{\log\log x}{\log x}$ by Theorem \ref{H1bounds}, we have
\begin{equation}\label{H1sym}
H_1(x,x^{1/4},x^{3/4}) \asymp \frac{x}{\log x}.
\end{equation} 
The lower bound is obtained by considering $n=ap$ with $a\le x^{1/4}$ and
$p$ a prime in $(\frac12 x^{3/4}, x^{3/4}]$.  Now suppose $n>x^{3/4}$,
$\tau(n,x^{1/4},x^{3/4})=1$, and $d|n$ with $x^{1/4}<d\le x^{3/4}$.
Since $\frac{n}{d}<x^{3/4}$, we have $\frac{n}{d} \le x^{1/4}$, hence
$d>x^{1/2}$.  If $d$ is not prime, then there is a proper divisor of
$d$ that is $ \ge \sqrt{d} > x^{1/4}$, a contradiction.  Thus,
$d$ is prime and $n=da$ with $a\le x^{1/4}$.  The   upper
bound in \eqref{H1sym} follows.

There is an application
to the Erd\H os-Montgomery function $g(n)$, which counts
the number of pairs of
consecutive divisors $d,d'$ of $n$ with $d|d'$ (see \cite{ET81},
\cite{ET83}).  The following sharpens Th\'eor\`eme 2 of \cite{Ten87}.

\begin{cor}\label{EMcor}
We have
$$
\frac{1}{x}\sum_{n\le x} g(n) \asymp \frac{(\log x)^{1-\del}}
{(\log\log x)^{3/2}}.
$$
\end{cor}

\Proof
The upper bound follows from $g(n) \le \tau^+(n)$ and Corollary
\ref{tauplus}.  We also have $g(2n) \ge I(n)$, where 
$I(n)$ is the number of $d|n$ such that if $d'|n, d'\ne d$, then
$d'>2d$ or $d'<d/2$ (see \S \ref{sec:lower}).  We quickly derive
(cf.\ \cite[p.~185]{Ten87})
$$
\sum_{n\le x} g(n) \gg \frac{x}{\log x} \sum_{m\le x^{1/5}} \frac{I(m)}
{m}.
$$
Applying \eqref{Isum} with $g=1$, $y=\sqrt{x}$, $\a=\frac15$ and
$\sg=\log 2$, we obtain
$$
\sum_{m\le x^{1/5}} \frac{I(m)}{m} \gg
\frac{(\log x)^{2-\del}}{(\log\log x)^{3/2}},
$$
and this gives the corollary.
\Endproof\vskip4pt

In a forthcoming paper, the author and G. Tenenbaum \cite{FT} show that
Conjecture \ref{tenconj1} is false when $z$ is close to $z_0(y)$.
Specifically, if $c>0$ is fixed, $g(y)>0$,
$\ds \lim_{y\to\infty} g(y) = 0$, $y\le x^{1/2-c}$ and 
$y+1 \le z \le y + y (\log y)^{1-\log 4 + g(y)}$, then
$$
H_1(x,y,z) \sim H(x,y,z) \qquad (y\to\infty).
$$
Moreover, the lower bound in \eqref{Prsmallz} is the true
order of $\frac{H_r(x,y,z)}{H(x,y,z)}$ for $r\ge 2$.

\Subsec{Divisors of shifted primes}
The methods developed in this paper may also be used to estimate 
a more general quantity
$$
H(x,y,z;\AA) = | \{ n\le x: n\in \AA, \tau(n,y,z)\ge 1 \}
$$
for a set $\AA$ of positive integers which is
well enough distributed in arithmetic progressions so that the
initial reductions (Lemmas \ref{Hupper},
\ref{basiclower}, \ref{Hlower2}) can be made to work.
An example is $\AA$ being an arithmetic progression $\{un+v:n\ge
1\}$, where the modulus $u$ may be fixed or grow at a moderate rate as a
function of~$x$.  Estimates with these $\AA$ are given in \cite{FKSY}.

One example which we shall examine in this paper
is when $\AA$ is a set of \emph{shifted primes} 
(the set $P_\lam = \{q+\lam: q \text{ prime} \}$ 
for a fixed non-zero $\lam$).
Results about the multiplicative structure of
shifted primes play an important role in many number theoretic
applications, especially in the areas of primality testing,
integer factorization and cryptography.
Upper bounds for $H(x,y,z;P_\lam)$ have been given by Pappalardi 
(\cite[Th.~3.1]{Pappalardi}),  Erd\H os and Murty (\cite[Th.~2]{EM}) and 
Indlekofer and Timofeev (\cite[Th.~2 and its corollaries]{IT}). 
Improving on all of these estimates, 
we give upper bounds of the expected order of magnitude, for all
$x,y,z$ satisfying $y\le \sqrt{x}$.

\begin{thm}\label{thm:shifted}
Let $\lam$ be a fixed non-zero integer.  
Let $1 \le y \le \sqrt{x}$ and  $y+1\le z \le x$.  Then
$$
H(x,y,z;P_\lam) \ll_\lam \begin{cases} \dfrac{H(x,y,z)}{\log x} & z \ge y +
 (\log y)^{2/3} \\ \quad \\
\dfrac{x}{\log x} \sum_{y<d\le z} \dfrac{1}{\phi(d)} & y < z \le y+ 
(\log y)^{2/3}. \end{cases}
$$
\end{thm}

Lower bounds are much more difficult, depending heavily on the 
distribution of primes in arithmetic progressions.  The special
case $z=y+1$ already presents difficulties, since then
$H(x,y,y+1;P_\lam)$
counts the primes $\le x$ in the progression $-\lam \pmod{\fl{y}+1}$.
If the interval $(y,z]$ is long, however, we can make use of 
average result for primes in arithmetic progressions.

\begin{thm}\label{thm:shiftedlong}
For fixed $\lambda, a, b$ with $\lambda\ne 0$ and $0\le a<b\le 1${\rm ,}
$$
H(x,x^a,x^b;P_\lam) \gg_{a,b,\lam} \frac{x}{\log x}.
$$ 
\end{thm}

Theorem \ref{thm:shiftedlong} has an application to counting finite
fields for which there is a curve with Jacobian of small exponent \cite{FS}.

%
%
\Subsec{Outline of the paper}
In Section~\ref{sec:prelim} we give a few preliminary lemmas about
 primes and sieve counting functions. 
Sections \ref{sec:upper} and \ref{sec:lower} provide an outline of the
upper and lower bound arguments with
most proofs omitted.  These tools are combined to prove Theorems
\ref{thm1}, \ref{thm:shorts}, \ref{thm:squarefree}, \ref{H1bounds} and
\ref{Hrbounds} in Section~\ref{sec:theorems}.

The first step in all estimations is to relate
the average behavior of\break
$\tau(n,y,z)$, which
contains local information about the divisors of $n$, with average
behavior of functions which measure global distribution of
divisors.  This is accomplished in Section~\ref{sec:initial}.
The upper and lower bound arguments begin to diverge after this point.
In general, the upper bounds are more difficult, since one may restrict
ones attention to
numbers with nice properties for the lower bounds.
The prime divisors of $n$ which are $<z/y$ play an insignificant
role in the estimation of $H(x,y,z)$.
For example, if $y<d\le 2y$, then $md\le z$ for $1\le m\le z/(2y)$.
By the same reason, the prime factors of $n$ which are $\le z/y$
play a very important role in the estimation of $H_r(x,y,z)$.
Quantifying this difference of roles
for the upper bounds in Section~\ref{sec:Sstar} is
much more difficult than for the lower bounds in Section~\ref{sec:isolated}, although the underlying idea is
the same. 

In Section~\ref{sec:Sstar}, both  $H(x,y,z)$ and $H_r(x,y,z)$ are bounded 
above in terms of a quantity $S^*(t;\eta)$, which is an
average over square-free $n$  whose prime factors lie in
$(z/y,z]$ of a global divisor function of $n$.
The contribution to  $S^*(t;\eta)$ from those $n$ with
exactly $k$ prime factors is then estimated
in terms of an integral over $\RR^k$ of an elementary but complicated
function.  Strong estimates for this integral
are proved in Section 13, and depend on new probability bounds for
uniform order statistics given in Lemma 11.1 (see \S 11
for relevant definitions).  

The lower bound argument follows roughly the same outline as the
upper bound, but the details are quite different.  Averages
over the `global' divisor functions are estimated in terms of
averages of a function which counts `isolated' divisors of numbers
(divisors which are not too close to other divisors) in Section~\ref{sec:isolated}.  Averages
over counts of isolated divisors of numbers with $k$ prime factors are bounded
below in terms of the volume of a certain complicated region in
$\RR^k$.  Bounding from below the volume of this region makes use of the
bounds on uniform order statistics from Section~\ref{sec:uos}, and this is
accomplished in Section~\ref{sec:lowvol}.   For $z \ge y+y(\log y)^{1-\log
4+b}$, $b>0$ fixed, we need only take a single value of $k$.

There is an alternative approach to obtaining lower bounds for $H(x,y,z)$
which avoids the use of bounds for order statistics
(see \S 2 of \cite{Hxy2y}), but they
appear to be necessary for our upper
bounds and for our lower bounds for $H_r(x,y,z)$.

Finally in Section~\ref{sec:shifted}, we apply the upper bound tools
developed in the prior sections to give upper   estimates for
$H(x,y,z;P_\lam)$, proving Theorem \ref{thm:shifted}.  Theorem~\ref{thm:shiftedlong} is much simpler
and has a self-contained proof in
Section~\ref{sec:shifted}.

A relatively short, self-contained proof that
$$
H(x,y,2y) \asymp \frac{x}{(\log y)^{\del}(\log\log y)^{3/2}}  \qquad 
(3\le y\le \sqrt{x})
$$
is given in \cite{Hxy2y}.  
Aside from part of the lower bound argument,
the methods are those given here, omitting
complications which arise in the general case.

\bigskip

\Subsec{Heuristic arguments for $H(x,y,z)$}
Since the prime factors of $n$ 
which are $<z/y$ play a very insignificant role,
we essentially must count how many $n\le x$ have $\tau(n',y,z)\ge 1$, where
$n'$ is the product of the primes dividing $n$ lying between $z/y$ and $z$.
For simplicity, assume $n'$ is squarefree, $n'\le z^{100}$
and has $k$ prime factors.  When
$z \ge y+y(\log y)^{1-\log 4+c}$, $c>0$ fixed,
the majority of such $n$ satisfy $k-k_0=O(1)$, where
$$
k_0 = \fl{\frac{\log\log z - \log \eta}{\log 2}}.
$$
For example, most integers $n$ which have a
divisor in $(y,2y]$ have $\frac{\log\log y}{\log 2} + O(1)$ prime
factors $\le 2y$.

To see this, assume for the moment that the set $D(n')=\{ \log d : d|n' \}$ is
uniformly distributed in $[0,\log n_1]$.  Then the probability that
$\tau(n',y,z)\ge 1$ should be about $2^k \frac{\eta}{\log n_1} \asymp
2^k \frac{\eta}{\log z}$.  This is $\gg 1$ precisely when $k\ge k_0+O(1)$.
Using the fact (e.g.\ Theorem 08 of \cite{Divisors}) 
that the number of $n\le x$ with $n'$ having $k$ prime factors is
approximately 
$$
\frac{x \log (z/y+2)}{\log z} \frac{(\log\log z - \log\log (z/y+2))^k}{k!},
$$
we obtain a heuristic estimate for $H(x,y,z)$ which matches the upper bounds
of Theorem T1, sans the $\log\log (3/u)$ factor in (iii).  When $\b=o(1)$ or
$\eta>1$, this is slightly too big.  The reason stems from the uniformity
 assumption about $D(n')$.  In fact, for most $n'$ with about $k_0$ prime
 factors, the set $D(n')$ is far from uniform, possessing many clusters
of close divisors and large gaps between them.  This substantially decreases
the likelihood that $\tau(n',y,z)\ge 1$.  The cause is slight irregularities in
the distribution of prime factors of $n'$ which are guaranteed ``almost
surely''  by large deviation results of probability theory (see e.g.\ Ch. 1 of
\cite{Divisors}). 
The numbers $\log\log p$ over $p|n'$ are well-known to behave like
random numbers in $[\log\log \max(2,z/y), \log\log z]$, and any prime that is
slightly below its expectation leads to ``clumpiness'' in $D(n')$.
What we really should count is
the number of $n$ for which $n'$ has $k$ prime factors \emph{and} $D(n')$ is
roughly uniformly distributed.  This corresponds to asking for the prime
divisors of $n'$ to lie all above their expected values.  An analogy from
probability theory is to ask for the likelihood that a 
random walk on the real numbers,
with each step haveing zero expectation, stays completely to the right of
the origin (or a point just to the the left of the origin) 
after $k$ steps.  In Section \ref{sec:uos} we give estimates for
this probability.  In the case $z=2y$, the desired probability is about
$1/k \asymp 1/\log\log y$, which accounts for the discrepancy between the
upper estimates in Theorem T1 and the bounds in Theorem 1.

\Subsec{Some open problems}
 (i)  Strengthen Theorem \ref{thm1} to an asymptotic formula. \vskip2pt
\begin{itemize} \isep
\ritem{(ii)} Determine the order of $H_r(x,y,z)$ for $r\ge 2$ and $z \ge
  y^{C}$ (see the conjecture at the end of \S \ref{sec:theorems}).

\ritem{(iii)} Establish the order of $H_r(x,y,z)$ when $yz \ge x^{1-c}$. 

\ritem{(iv)} Make the dependence on $r$ explicit in Theorem \ref{Hrbounds}
  and Corollary \ref{epscor}.
Hall and Tenenbaum (\cite{Divisors}, Ch. 2) conjecture that for
  each $r\ge 2$,
$$
\lim_{y\to \infty} \frac{\eps_r(y,2y)}{\eps(y,2y)} = d_r > 0.
$$
In light of \eqref{Prlargez},
the sequence $d_1,d_2,\ldots$ may not be monotone.

\ritem{(v)} Provide lower bounds for $H(x,y,z;P_\lam)$ of the expected
  order for other $y,z$ not covered by Theorem \ref{thm:shiftedlong}.
\end{itemize}

\demo{Acknowledgements}
The author thanks
Bruce Berndt, Valery Nevzorov, Walter Philipp, Steven Portnoy, and Doug West
for helpful conversations regarding probability estimates for uniform order
statistics. The author also thanks G\'erald Tenenbaum
for several preprints of his work and for informing the author
about the theorem of Rogers mentioned above, and thanks\break
Dimitris Koukoulopoulos for discussions which led to a simplification
of the proof of Lemma \ref{sumW}.
The author is grateful to his wife, Denka Kutzarova,
for constant support and many helpful conversations about the paper.
Much of this paper was written while the author enjoyed the hospitality
of the Institute of Mathematics and Informatics, Bulgarian Academy of
Sciences.  Finally, the author acknowledges the referee for a thorough
reading of the paper and for helpful suggestions.

This work was partially
supported by National Science Foundation Grant DMS-0301083.

 \section{Preliminary lemmas}\label{sec:prelim}

{\it Further notation}.
$P^+(n)$ is the largest prime factor of $n$, $P^-(n)$ is the 
smallest prime factor of $n$.  Adopt the   conventions $P^+(1)=0$ and
$P^-(1)=\infty$.  Also, $\omega(n)$ is the number of distinct prime
divisors of $n$, $\Omega(n)$ is the number of prime power divisors of $n$,
$\pi(x)$ is the number of primes $\le x$, $\tau(n)$ is the number of
divisors of $n$.  $\PP(s,t)$ is the set of positive integers composed
of prime factors $p$ satisfying $s <p \le t$.  Note that for all
$s,t$ we have $1\in \PP(s,t)$.  $\PP^*(s,t)$ is the set of square-free
members of $\PP(s,t)$.

We list a  few  estimates from prime number theory and sieve
theory.  The first is the Brun-Titchmarsh inequality
and the second is a consequence of the Prime Number Theorem with
classical de la Val\'ee Poussin error term.

\begin{lem}\label{BrunTitch}
Uniformly in $x>y>1${\rm ,} we have $\pi(x)-\pi(x-y) \ll \frac{y}{\log y}$.
\end{lem}

\begin{lem}\label{sum1p}
For certain constants $c_0, c_1${\rm ,} 
$$ 
\sum_{p\le x} \frac{1}{p} = \log\log x + c_0 + O(e^{-c_1\sqrt{\log x}})
\qquad (x\ge 2).
$$ 
\end{lem}

The next result is a simple
application of Brun's sieve (see \cite{HR}) together with
the Prime Number Theorem with classical error term.
Although not required for this paper, using bounds on the
number of primes in short intervals, the best to date of which is
\cite{BHP}, allows us to take  $\Delta$ as small as $x^{0.525}$ in
the lower bound in the next lemma. 

\begin{lem}\label{Phi}
Let $\Phi(x,z)$ be the number of integers $\le x${\rm ,} all of whose prime
factors are $>z$.
If $1 < z^{1/100} \le \Delta \le x${\rm ,} then
$$ 
\Phi(x,z) - \Phi(x-\Delta,z) \ll \frac{\Delta}{\log z}.
$$ 
If $x\ge 2z${\rm ,} $z$ is sufficiently large and $x e^{-(c_1/2)\sqrt{\log
    x}} \le \Delta\le x${\rm ,} then
$$
\Phi(x,z)-\Phi(x-\Delta,z) \gg \frac{\Delta}{\log z}.
$$
\end{lem}

The second tool is crude but quite useful due to its uniformity. 
A proof may be found in Tenenbaum \cite{Tenbook}.

\begin{lem}\label{Psi}
Let $\Psi(x,y)$ be the number of integers $\le x${\rm ,} all of whose prime
factors are $\le y$.  Then{\rm ,} uniformly in $x\ge y\ge 2${\rm ,} 
$$
\Psi(x,y) \ll x \exp\{ - \tfrac{\log x}{2\log y}\}.
$$
\end{lem}

\begin{lem}\label{Psixyz}
Uniformly in $x>0${\rm ,} $y\ge 2$ and $z\ge 1.5${\rm ,} we have
\begin{align}
\sum_{\substack{n\ge x \\ n\in \PP(z,y)}} \frac{1}{n} &\ll \frac{\log y}{\log z}
  e^{-\frac{\log x}{4\log y}}, \label{Psi1n} \\
\sum_{\substack{n\ge x \\ n\in \PP(z,y)}} \frac{\log n}{n}& \ll \frac{\log y 
  \log(xy)}{\log z}e^{-\frac{\log x}{4\log y}}, \label{Psilogn}
\end{align}
\end{lem}

\Proof Without loss of generality we may assume that $x\ge 1$.  
Put $\a=\frac{1}{4\log y} \le \frac25$.  The
result \eqref{Psi1n} is trivial unless $z<y$, in which case
\begin{align*}
\sum_{\substack{n\ge x \\ n\in \PP(z,y)}} \frac{1}{n} & \le x^{-\a} 
  \sum_{n\in\PP(z,y)} n^{\a-1} = x^{-\a} \prod_{z<p\le y} \( 1 + 
  p^{\a-1} + p^{2(\a-1)} + \cdots \)\\
&= x^{-\a} \exp\biggl\{ \sum_{z<p\le y} \(\frac{1}{p} + O\(p^{-6/5}+
  \frac{\a \log p}{p}\)\) \biggr\} \\
&\ll x^{-\a} \frac{\log y}{\log z},
\end{align*}
where we used Lemma \ref{sum1p} in the final step.   Now set
$$
F(t)=\sum_{\substack{n>t \\ n\in \PP(z,y)}} \frac{1}{n}.
$$
By \eqref{Psi1n} and partial summation,
\begin{align*}
\sum_{\substack{n\ge x \\ n\in \PP(z,y)}} \frac{\log n}{n} &= F(x) \log x +
  \int_x^\infty \frac{F(t)}{t}\, dt \\
&\ll \frac{\log y \, \log x}{\log z}  e^{-\frac{\log x}{4\log y}} + 
  \frac{\log y}{\log z} \int_x^\infty t^{-1-\frac{1}{4\log y}}\, dt \\
&\ll \frac{\log y \log(xy)}{\log z}e^{-\frac{\log x}{4\log y}}.
\end{align*}
\vglue-24pt
\Endproof\vskip8pt

We shall also need Stirling's formula
\begin{equation}\label{Stirling}
k! = \sqrt{2\pi k} (k/e)^k (1 + O(1/k)),
\end{equation} 
although in most estimates weaker bounds will suffice.\smallbreak 

Our last lemma is a consequence of Norton's bounds \cite{Norton} 
for the partial sums of the exponential series.  
It is easily derived from Stirling's formula.

\begin{lem}\label{Nortonlem}
Suppose $0 \le h < m \le x$ and $m-h \ge \sqrt{x}$.  Then
$$
\sum_{h\le k\le m} \frac{x^k}{k!} \asymp \min \( \sqrt{x}, \frac{x}{x-m} \) 
\frac{x^m}{m!}.
$$
\end{lem}

\ifdraft
\vfil\eject
\else
\fi

%
%
%
%
%
\section{Upper bounds outline}\label{sec:upper}
%
%
%
%
%

Initially, we bound $H(x,y,z)$ and $H_r(x,y,z)$ in terms of averages of
the functions
\begin{alignat}{5}
L(a;\sg) &=  \text{meas} (\LL(a;\sg)), &\quad
 \LL(a;\sg) &= \{x: \tau(a,e^x,e^{x+\sg}) \ge 1\} \label{Ldef}, \\
L_r(a;\sg) &=  \text{meas}  (\LL_r(a;\sg)), &\quad
 \LL_r(a;\sg) &=  \{x: \tau(a,e^x,e^{x+\sg}) =r\}\label{Lrdef}.
\end{alignat}
Here $\text{meas}(\cdot)$ denotes Lebesgue measure.
Both functions measure the global distribution of
divisors of $a$.  Before launching into the estimation of $H$ and $H_r$,
we list some basic inequalities for $L(a;\sg)$.

\begin{lem}\label{Lineq}
We have
\begin{enumerate}
\item  $L(a;\sg) \le \min(\sg \tau(a), \sg+\log a)${\rm ;}
\item  If $(a,b)=1${\rm ,} then $L(ab;\sg) \le \tau(b) L(a;\sg)${\rm ;}
\item  If $(a,b)=1${\rm ,} then $L(ab;\sg) \le L(a;\sg+\log b)${\rm ;}
\item  If $\g \le \sg${\rm ,} then $L(a;\sg) \le (\sg/\g) L(a;\g)${\rm ;}
\item  If $p_1 < \cdots < p_k${\rm ,} then
$$
L(p_1\cdots p_k;\sg) \le\min_{0 \le j\le k} 2^{k-j} (\log(p_1\cdots p_j)+\sg).
$$
\end{enumerate}
\end{lem}

\Proof  Part (i) is immediate, since
$$
\LL(a;\sg) = \bigcup_{d|a} [-\sg+\log d,\log d) \subseteq [-\sg,\log a).
$$
Parts (ii) and (iii) follow from
$$ 
\LL(ab;\sg) = \bigcup_{d|b} \{x+\log d : x\in \LL(a;\sg)  \} \subseteq
\LL(a;\sg+\log b).
$$ 
Since $\LL(a;\sg)$ is a union of intervals of length $\sg$, we obtain (iv). 
Combining parts (i) and (ii) with $a=p_1\cdots p_j$ and
$b=p_{j+1}\cdots p_k$ yields (v).
\Endproof\vskip4pt

  Define
\begin{align}
S(t;\sg) &= \sum_{P^+(a) \le te^{\sg}}
  \frac{L(a;\sg)}{a\log^2(t/a+P^+(a))},
  \label{Sdef} \\
S_s(t;\sg) &= \sum_{P^+(a) \le te^{\sg}}
  \frac{L_s(a;\sg)}{a\log^2(t/a+P^+(a))}. \label{Srdef}
\end{align}

\begin{lem}\label{Hupper}
Suppose $100 \le y \le \sqrt{x}${\rm ,}
$z=e^\eta y \le \min(x^{5/8},y^{\log\log y})$ and $\eta\ge \frac{1}{\log y}$.
If $x/\log^{10} z \le \Delta\le x${\rm ,}
$$
H(x,y,z) -H(x-\Delta,y,z) \ll \Delta \max_{y^{1/2} \le t\le x} S(t;\eta).
$$
If in addition $y\ge y_0(r)${\rm ,} then
$$
H_r(x,y,z) \ll_r x \max_{\substack{1\le s\le r \\ \nu(s)\le \nu(r)}}
 \max_{y^{1/2} \le t\le x} S_s(t;\eta).
$$
\end{lem}

Lemma \ref{Hupper} will be proved in Section~\ref{sec:initial}.
If $m<z/y$ then $\tau(n,y,z)\ge 1$ implies
 $\tau(nm,y,z^2/y)\ge 1$ and we expect (and prove) that $H(x,y,z)$ and
$H(x,y,z^2/y)$ have the same order.  Thus, for the problem of bounding
$H(x,y,z)$, the prime
factors of $n$ below $z/y=e^\eta$ can essentially be ignored.  
For the problem of bounding
$H_r(x,y,z)$, the prime factors of $n$ less than $z/y$ cannot be
ignored and they play a different
role in the estimation than the prime factors $>z/y$.  
In the next two lemmas, we estimate 
both $S(t;\sg)$ and $S_s(t;\sg)$ in terms of 
the quantity
\begin{equation}\label{Sstardef}
S^*(t;\sg) =  \sum_{a\in \PP^*(e^\sg,te^{\sg})}
  \frac{L(a;\sg)}{a \log^2(t^{3/4}/a+P^+(a))}.
\end{equation} 
Occasionally we will have need of the trivial lower bound
\begin{equation}\label{Sstarlower}
S^*(t;\sg) \ge \frac{\sg}{\log^2 t},
\end{equation} 
obtained by taking the term $a=1$ in \eqref{Sstardef}.

\begin{lem}\label{Hupper2}
Suppose $t$ is large and $0 < \sg \le \log t$. Then
$$
S(t;\sg) \ll (1+\sg) S^*(t;\sg).
$$
\end{lem}

Particularly important in the estimation of $S_s(t;\sg)$
is the distribution of the
gaps between the first $r+1$ divisors of $a$, which
ultimately depends on the power of 2 dividing $r$. 

\begin{lem}\label{Hrupper}  
Suppose $r\ge 1${\rm ,} $C>1${\rm ,} $y\ge y_0(r,C)${\rm ,} $z=e^\eta y${\rm ,} $z\le x^{5/8}$
and $e^{100rC} y \le z \le y^{C}$.  Then
$$
H_r(x,y,z) \ll_{r,C} x (\log \eta)^{\nu(r)+1} \max_{y^{1/2} \le t\le x} 
S^*(t;\eta).
$$
\end{lem}

Lemmas \ref{Hupper2} and \ref{Hrupper} will be proved in Section~\ref{sec:Sstar}.

To deal with the factor
$\log^2 (t^{3/4}/a+P^+(a))$ appearing in \eqref{Sstardef}, define
\begin{equation}\label{Tdef}
T(\sg,P,Q) = \sum_{\substack{a\in \PP^*(e^\sg,P) \\ a\ge Q}} 
\frac{L(a;\sg)}{a}.
\end{equation} 
If $a\le t^{1/2}$ or $P^+(a)>t^{1/3000}$,
then $\log^2 (t^{3/4}/a+P^+(a)) \gg \log^2 t$. 
Otherwise, $e^{e^{g-1}} < P^+(a) \le e^{e^g}$ for some integer
$g$ satisfying $e^{\sg} \le e^{e^g} \le t^{1/1000}$.
Thus we have
\begin{equation}\label{ST}
S^*(t;\sg) \ll \frac{T(\sg,t e^\sg,1)}{\log^2 t} + \sum_{\substack{g\in\ZZ,
g\ge 1 \\ e^{\sg} \le e^{e^g} \le t^{1/1000}}} e^{-2g} T(\sg,e^{e^g},t^{1/2}).
\end{equation} 
We break up the sum in $T(\sigma,P,Q)$ according to the value of
$\omega(a)$ and define
$$
T_k(\sg,P,Q) = \sum_{\substack{a\in \PP^*(e^\sg,P) \\ a\ge Q \\ \omega(a)=k}} 
\frac{L(a;\sg)}{a}.
$$
Note that the definition of $k$ given here is slightly different from that
mentioned in the heuristic argument of subsection 1.5, but usually differs
only  by $O(1)$.
By Lemma \ref{sum1p} and part (v) of Lemma \ref{Lineq}, 
$T_k(\sg,P,Q)$ will be bounded in terms of
\begin{equation}\label{Udef}
U_k(v;\a) = \int_{R_k} \min_{0\le j\le k} 2^{-j} \( 2^{v\xi_1} +
\cdots + 2^{v\xi_j} + \a \)\, d\bx,
\end{equation} 
where
$$
R_k = \{ \bx \in \RR^k : 0\le \xi_1 \le \cdots \le \xi_k \le 1 \}.
$$
For convenience, let $U_0(v;\a)=\a$.

%
%

\begin{lem}\label{Tnint}
Suppose $P\ge 100${\rm ,} $0 < \sg < \log P${\rm ,} and
$Q\ge 1$.  Let
$$
v = \cl{\frac{\log\log P -\max(0,\log \sg)}{\log 2}} \\
$$
and suppose $0\le k\le 10v$.  Then
$$
T_k(\sg,P,Q) \ll e^{-\frac{\log Q}{\log P}} 
(\sg+1) (2v\log 2)^k U_k(v;\min(1,\sg)).
$$
\end{lem}

Lemma \ref{Tnint} will be proved in Section~\ref{sec:uppertoint}.
As a rough
heuristic, $2^{v\xi_1} + \cdots + 2^{v\xi_j} \ll 2^{v\xi_j}$ most of the
time.  Thus, bounding $U_k(v;\a)$ boils down to determining the distribution in
$R_k$ of the function
$$
F(\bx) = \min_{1\le j\le k} (\xi_j - j/v).
$$
The numbers $\xi_1,\dots,\xi_k$ can be regarded as independent uniformly
distributed random variables on $[0,1]$, relabeled to have the above ordering,
and are
known as \emph{uniform order statistics}.
Making this heuristic precise, and using results about the
distribution of uniform order statistics from Section~\ref{sec:uos},
 leads to the next result, which will be proved in Section~\ref{sec:upperint}.

\begin{lem}\label{Unlem}
Suppose $k,v$ are integers with $0\le k\le 10v$
and $0 < \a \le 1$.  Then
$$
U_k(v;\a) \ll \frac{\a \min \bigl( k+1,(1 +|v-k-\tfrac{\log \a}{\log
    2}|^2) \log(2/\a) \bigr)} {(k+1)! (\a 2^{k-v}+1)}.
$$
\end{lem}

Notice that, as a function of $k$,
the bound in Lemma \ref{Unlem} undergoes a change of behavior
at $k = \fl{v - \frac{\log \a}{\log 2}}$.
It is now straightforward to give 
a relatively simple upper bound for $T(\sg,P,Q)$.

%
%

\begin{lem}\label{Tlem}  Suppose $P$ is sufficiently large{\rm ,} $Q\ge 1$, and
$$
(\log P)^{-1} \le \sg < \log P.
$$
Define $\theta=\theta(\sg,P)$ and $\nu=\nu(\sg,P)$ by $\sg = 
(\log P)^{-\theta}$ and $\theta=\log 4 -1 - \nu (\log\log P)^{-1/2}$
\/{\rm (}\/these quantities are related to those in \eqref{etabetaxi}{\rm ).}
Then
$$
T(\sg,P,Q) \ll \begin{cases} e^{-\frac{\log Q}{\log P}}
\dfrac{\sg^{\del-1} (\log P)^{2-\del}}
{(\log \frac{\log P}{1+\sg} + 1)^{3/2}} & 
   \quad \sg \ge 1 \\ \quad \\  e^{-\frac{\log Q}{\log P}}
\dfrac{(\log P)^{2-G(\th)}\log(2/\sg)}{\max(1,\nu) \log\log P} 
& \quad \sg < 1. \end{cases}
$$
\end{lem}

\Proof Define $v$  as in the statement of Lemma \ref{Tnint} and
 set $\a=\min(1,\sg)$.  Put $\g = e^{-\frac{\log Q}{\log P}}$.
By Lemmas \ref{Tnint}
and \ref{Unlem}, when $0\le k\le 10v$,
\begin{equation}\label{Tn1}
\begin{split}
T_k(\sg,P,Q) &\ll \g \a (\sg+1) Z_k \ll \g \sg Z_k, \\
Z_k &= \frac{\min(k+1,(1 +|v-k-\tfrac{\log \a}{\log 2}|^2) \log(2/\a))}
{(k+1)! (\a 2^{k-v}+1)} (2v\log 2)^k.
\end{split}
\end{equation} 
Put $k_1=\fl{v-\frac{\log\a}{\log 2}}$ and note that $v\le k_1 \le 2v$. 
Now,
\begin{equation}\label{Znsum}
\sum_{k_1\le k\le 10v} Z_k \ll \log(2/\a) \sum_{b\ge 0} \frac{b^2+1}{2^b}\,
\frac{(2v\log 2)^{k_1+b}}{(k_1+b+1)!} \ll \frac{\log(2/\a) (2v\log 2)^{k_1}}
{(k_1+1)!}.
\end{equation} 
By $L(a;\sg) \le 2^{\omega(a)} \sg$, 
\begin{align*}
\sum_{k\ge 10v} T_k(\sg,P,Q) &\le Q^{-1/\log P} \sg \sum_{k\ge 10v} 2^k 
\sum_{\substack{a\in\PP^*(e^\sg,P) \\ \omega(a)=k}} \frac{1}{a^{1-1/\log P}}
\\ 
&\le \g \sg \sum_{k\ge 10v} \frac{2^k}{k!} \biggl( \sum_{e^\sg < p \le P}
\frac{1}{p^{1-1/\log P}} \biggr)^k \\
&= \g \sg \sum_{k\ge 10v} \frac{(2v\log 2 + O(1))^k}{k!} \\
&\ll \g \sg \frac{(2v\log 2)^{10v}}{(10v)!} \\
&\ll \g \sg \frac{(2v\log 2)^{k_1}}{(k_1+1)!}.
\end{align*}
Together with \eqref{Tn1} and \eqref{Znsum}, we conclude that
\begin{equation}\label{Tnsum}
\sum_{k\ge k_1} T_k(\sg,P,Q) \ll \frac{\g \sg \log(2/\a) (2v\log 2)^{k_1}}
{(k_1+1)!}.
\end{equation} 

Suppose that $\theta \le \frac13$, so that $k_1 \le \frac43 v$.  Since
$\frac43 < 2\log 2$, \eqref{Tn1} implies
$$
\sum_{0\le k\le k_1} Z_k \ll \log(2/\a) \sum_{0\le k\le k_1} (k_1-k+1)^2
\frac{(2v\log 2)^k}{(k+1)!} \ll \frac{\log(2/\a) (2v\log 2)^{k_1}}
{(k_1+1)!}.
$$
Combined with \eqref{Tn1}, \eqref{Tnsum} and Stirling's formula, this gives
\begin{equation}\label{smalltheta}
\begin{split}
T(\sg,P,Q) &\ll  \frac{\g \sg \log(2/\a) (2v\log 2)^{k_1}}
{(k_1+1)!} \\
&\ll \frac{\g \sg \log(2/\a)}{v^{3/2}} \pfrac{2ev\log 2}
{k_1}^{k_1} \qquad (\th \le \tfrac13).
\end{split}
\end{equation} 
When $\sg\ge 1$, we have $\th \le 0$, $\a=1$, $k_1=v$, and 
$$
(2e\log 2)^v \asymp \pfrac{\log P}{\sg}^{2-\del},
$$
and so the lemma follows in this case.  We also have
$$
v=\frac{\log\log P}{\log 2}+O(1), \qquad k_1 = (1+\th) \frac{\log\log P}
{\log 2} + O(1) \qquad (0 \le \th \le 1)
$$
and
$$
\pfrac{2ev\log 2}{k_1}^{k_1} \asymp (\log P)^{2+\th-G(\th)} \qquad
(0\le \th \le \log 4 -1).
$$
Thus, if $0\le \th \le \frac13$, then $\nu \asymp (\log\log P)^{1/2}$ and
the lemma follows from \eqref{smalltheta}.
If $\frac13 \le \theta \le 1$, \eqref{Tn1} and Lemma \ref{Nortonlem} give
\begin{align*}
\sum_{0\le k\le k_1} Z_k \ll \sum_{0\le k\le k_1} \frac{(2v\log 2)^k}{k!} &\ll
\begin{cases}
e^{2v\log 2} ,& k_1 \ge 2v\log 2 -\sqrt{v} \\ \quad \\
\dfrac{\sqrt{v} (2v\log 2)^{k_1}}{\nu k_1!} ,& k_1 <  2v\log 2 -\sqrt{v}
\end{cases} \\
&\ll  \frac{(\log P)^{2+\th-G(\th)} \log(2/\sg)}{\max(1,\nu)\log\log P}.
\end{align*}
Together with \eqref{Tn1} and \eqref{Tnsum}, this proves the lemma in
the final case.
\hfill\qed

%
%

\begin{lem}\label{Sstar}
Suppose $t$ is large and $\sg \ge (\log t)^{-1/2}$.  Put
$\th=\th(\sg,t)$ and $\nu=\nu(\sg,t)$.   Then
$$
S^*(t,\sg) \ll \begin{cases} \dfrac{\sg^{\del-1}(\sg+\log t)^{2-\del}}
{(\log t)^2 (\log(\sg+\log t) -\log \sg+1)^{3/2}},
& \sg \ge 1 \\ \quad \\  \dfrac{\log(2/\sg)}{(\log\log t) \max(1,\nu) 
(\log t)^{G(\theta)}}, & \sg\le 1. \end{cases}
$$
\end{lem}

\Proof
First suppose $\sg\ge 1$.  By \eqref{ST} and Lemma \ref{Tlem}, writing
$g=\fl{\log\log t} - \ell$ gives
\begin{eqnarray*}
S^*(t;\sg) &\ll &\frac{\sg^{\del-1}}{(\log t)^{\del}} \biggl[ \pfrac{\sg+\log t}
{\log t}^{2-\del} \frac{1}
  {(\log(\sg+\log t) - \log \sg+1)^{3/2}}   \\&&
+ \sum_{1\le \ell \le \log\log t -\log \sg} \frac{e^{\del \ell}}
  {e^{e^{\ell-1}} (\log\log t - \log \sg + 1 - \ell)^{3/2}}\biggr].
\end{eqnarray*}
The sum on $\ell$ is empty if $\sg>\log t$.  Otherwise, the sum on $\ell$ is
dominated by terms with $\ell\ll 1$, and this proves the lemma in this case.

Suppose that $\sg < 1$. 
By Lemma \ref{Tlem}, the first term in \eqref{ST} is
$$
\ll \frac{\log(2/\sg)}{\max(1,\nu) (\log t)^{G(\theta)}\log\log t}.
$$
We use Lemma \ref{Tlem} when $e^{-g} \le \sg$.
The contribution of these terms (if any) in \eqref{ST} is
\begin{align*}
&\ll \log(2/\sg) \!\!\! 
  \sum_{\log\frac{1}{\sg} < g\le \log\log t} 
  \frac{e^{-gG(-(\log \sg)/g)}}{e^{(e^{-g-1}\log t)}g
  \max(1,\sqrt{g}(1-\log 4 - \tfrac{\log \sg}{g}))} \\
&\ll \frac{\log(2/\sg)}{\max(1,\nu) (\log t)^{G(\theta)} \log\log t}.
\end{align*}
When $g < \log(1/\sg) \le \frac12 \log\log t$, Lemma \ref{Tnint}
gives
$$
T_k(\sg,e^{e^g},t^{1/2}) \ll e^{-\frac12 \sqrt{\log t}}
 (2v\log 2)^k U_k(v;\sg) \le e^{-\frac12 \sqrt{\log t}} (2v\log 2)^k \sg/k!,
$$
where $v=\frac{g}{\log 2}+O(1)$.  Summing on $k$ and $g$ yields
$$
\sum_{g < \log\frac{1}{\sg}} e^{-2g} T(\sg,e^{e^g},t^{1/2}) \ll 
\sum_{g\le \frac12 \log\log t} \sg e^{-\frac12\sqrt{\log t}}
\ll \exp\{ -\frac13 (\log t)^{1/2} \},
$$
which is negligible compared to the contribution of the terms in \eqref{ST}
with $g \ge \log(1/\sg)$.
\Endproof

For fixed $\sigma$, $\th(\sg,t)$ is decreasing,
$\nu(\sg,t)$ is increasing and $G(\th)/\th$ is increasing, as functions
of $t$.  Thus, we have the following.

\begin{lem}\label{Sstarfinal}
Suppose $y$ is large and $\eta \ge (\log y)^{-0.4}$.  Then
$$
\max_{t\ge y^{1/2}} S^*(t;\eta) \ll 
\begin{cases} \dfrac{\eta^{\del-1}(\eta+\log y)^{2-\del}}
{(\log y)^2 (\log(\eta+\log y) -\log \eta+1)^{3/2}},
& \eta \ge 1 \\ \quad \\  \dfrac{\log(2/\eta)}{(\log\log y) \max(1,\nu(\eta,y)) 
(\log y)^{G(\theta(\eta,y))}} ,& \eta\le 1. \end{cases}
$$
\end{lem}

 \section{Lower bounds outline}\label{sec:lower}
 As with the upper bounds, we initially bound $H(x,y,z)$ in terms of
sums over $L(a;\sg)$ and bound $H_r(x,y,z)$ in terms of sums over
$L_s(a;\sg)$ (but only for $s=r$).  
The initial bounds are similar to those in Lemma \ref{Hupper}.

\begin{lem}\label{basiclower} 
Suppose $y_0 \le y <z = e^\eta y${\rm ,} $\frac{1}{\log^{20} y} \le \eta 
\le \frac{\log y}{100}${\rm ,}  $y\le \sqrt{x}$ and $x/\log^{10} z \le 
\Delta\le x$.  Then
\begin{align*}
H(x,y,z)-H(x-\Delta,y,z) &\ge H^*(x,y,z)-H^*(x-\Delta,y,z) \\
&\gg 
\frac{\Delta}{\log^2 y} \sum_{\substack{a\le y^{1/8} \\ \mu^2(a)=1}} 
\frac{L(a;\eta)}{a}.
\end{align*}
Suppose $r\ge 1${\rm ,} $0<c\le \frac18${\rm ,} $C>0${\rm ,} $y_0(r,c,C) \le y < z=e^{\eta} y${\rm ,} 
$\frac{1}{\log^2 y} \le
\eta \le C\log y$ and $z\le x^{1/2-c}$.  Then
$$ 
H_r(x,y,z) \gg_{r,c,C} 
\frac{x}{\log^2 y} \sum_{a\le y^{2c}} \frac{L_r(a;\eta)}{a}.
$$ 
\end{lem}

Lemma \ref{basiclower} will be proved in Section~\ref{sec:initial}.
Both $L(a;\sg)$ and $L_r(a;\sg)$ may be bounded below in terms of the
function
\begin{equation}\label{Idef}
I(n;\sg) = | \{ d|n: \tau(n,de^{-\sg},de^\sg)=1 \}|.
\end{equation} 
Introduced by Tenenbaum \cite{Ten87}, $I(n;\sg)$  counts
\emph{$\sg$-isolated} divisors of $n$.  

In the first part of Lemma \ref{basiclower}, take
square-free $a=h'h$, where
$h' \le z/y \le y^{1/100}$ and $P^-(h) > z/y$.  Clearly 
$$
L(h'h;\eta) \ge L(h;\eta) \ge \eta I(h;\eta),
$$
and summing over $g$ we obtain the following.

\begin{lem}\label{Hlower2}
Suppose $y_0 \le y <z = e^\eta y${\rm ,} 
$\frac{1}{\log^{20} y} \le \eta \le \frac{\log y}{100}${\rm ,}
 $y\le \sqrt{x}$ and $\frac{x}{\log^{10} z} \le \Delta\le x$.  Then
$$
H^*(x,y,z)-H^*(x-\Delta,y,z) \gg \frac{\eta (1 + \eta) \Delta}{\log^2 y}
\sum_{\substack{h\le y^{1/10} \\ P^-(h) > z/y \\ \mu^2(h)=1}}
\frac{I(h;\eta)}{h}.
$$
\end{lem}

We follow two methods for bounding $H_r(x,y,z)$ from below,
the first useful for  $z\ll y$  and the second useful for large $z$.

\begin{lem}\label{Hrlower1}
Suppose $r\ge 1${\rm ,} $0<c'\le \frac18${\rm ,} 
$y_0(r,c') \le y < z=e^{\eta} y \le x^{1/2-c'}$ and
$\frac{1}{\log^2 y} \le \eta \le \frac{c'\log y}{10r}$.  Then
$$ 
H_r(x,y,z) \gg_{r,c'} \frac{\eta^r x}{(\log y)^{r+1}} 
\sum_{a\le y^{c'/100r}} \frac{I(a;\eta)^r}{a}.
$$ 
\end{lem}

Lemma \ref{Hrlower1} and its proof are essentially taken from
Lemme 4 of Tenenbaum
\cite{Ten87}.  The main difference is the upper limit of allowable $z$:
Lemme 4 of \cite{Ten87} requires $z\le x^{\frac1{r+1}-c}$.

In the second method, the prime factors of $a$ which are $< z/y$ play
a special role as in Lemma \ref{Hrupper}.

\begin{lem}\label{Hrlower2}
{\rm (i)} Suppose $r\ge 1${\rm ,} $C>0${\rm ,} $0< c'\le \frac18${\rm ,} $y_0(r,c',C) \le y <
z=e^{\eta} y \le x^{1/2-c'}$ and
$1000r \cdot 3^{2r} \le \eta \le C \log y$.  Then
$$ 
H_r(x,y,z) \gg_{r,c',C} \frac{\eta (\log \eta)^{\nu(r)+1} x}{\log^2 y} 
\sum_{\substack{h\le y^{c'} \\ P^-(h) > e^{2\eta}}}
\frac{I(h;2\eta)}{h}.
$$ 
{\rm (ii)} If $r\ge 1${\rm ,} $0<c\le \frac18${\rm ,} $y\ge y_0(r,c)$ and
$y^2 \le z \le x^{1-c}/y${\rm ,} then
$$
H_r(x,y,z) \gg_{r,c} \frac{x (\log\log y)^{\nu(r)+1}}{\log z}.
$$
\end{lem}

Lemmas \ref{Hrlower1} and \ref{Hrlower2} will be proved in Section
\ref{sec:isolated}.
The number of isolated divisors of a number can be easily bounded
from below in terms of
\begin{equation}\label{Wdef}
W(a;\sg)=| \{ (d_1,d_2) : d_1|a, d_2|a, |\log(d_1/d_2)| \le \sg\}|.
\end{equation} 
This function, introduced by Hall \cite{Hall79}, is essential in the study
of the propinquity of divisors (see also \cite{HT81}, \cite{MT84},
\cite{MT85}, Chapters 4 and 5 of \cite{Divisors}, \cite{RaTen}, and
\cite{Ten03}). 
The following lemma is similar to Lemme 5 of Tenenbaum \cite{Ten87}.

%
%

\begin{lem}\label{Ilow} 
There exists $I(a;\sg)$ such that 
$$
I(a;\sg)^r \ge 2^{-r} \tau(a)^{r-1}(3\tau(a)-2W(a;\sg)).\label{low3}
$$
\end{lem}

\Proof
For each divisor $d$ of $a$ not counted by $I(a;\sg)$
there is at least one other divisor $d'$ satisfying $d/e^{\sg} \le d' \le
d e^{\sg}$, so that the pair $(d,d')$ is counted by $W(a;\sg)$. Thus
$$ 
W(a;\sg) \ge \tau(a) + (\tau(a)-I(a;\sg)) = 2\tau(a) - I(a;\sg).
$$ 
The lemma is trivial when $W(a;\sg) \ge \frac32 \tau(a)$.
Otherwise,
\vskip12pt
\hfill $
\displaystyle{I(a;\sg)^r \ge (2\tau(a)-W(a;\sg))^r \ge \pfrac{\tau(a)}{2}^{r-1} (\tfrac32
\tau(a)-W(a;\sg)).}
$ 
\Endproof \vskip12pt

With Lemma \ref{Ilow}, lower bounds for $H(x,y,z)$ and $H_r(x,y,z)$
are obtained via upper bounds on sums over $W(a;\sg)/a$. 
Such upper bounds are achieved by partitioning the primes into
sets $D_1, D_2, \ldots$ and separately
 considering numbers $a$ with a fixed number of prime factors
in each interval $D_j$.  

Each set $D_j$ will consist of the primes in an interval
$(\lam_{j-1},\lam_j]$, with $\lam_j \approx \lam_{j-1}^2$.  More precisely, 
let $\lam_0=1.9$ and inductively define $\lam_j$ for $j\ge 1$ as the
largest prime so that
\begin{equation}\label{Dj}
\sum_{\lam_{j-1} < p \le \lam_j} \frac{1}{p} \le \log 2.
\end{equation} 
For example, $\lam_1=2$ and $\lam_2=7$. 
Write $\lam_j=\exp \{2^{\mu_j} \}$.

\begin{lem}\label{muj}
There are constants $c_3, c_4$ so that $|\mu_j-j-c_3| \le
c_4 2^{-j}$ for all $j\ge 0$.
\end{lem}

\Proof  Clearly $\lam_j\to \infty$ as $j\to \infty$.
By \eqref{Dj} and Lemma \ref{sum1p} with crude error term,
$$
\log\log \lam_j - \log\log \lam_{j-1} = \log 2 + O(1/\log \lam_{j-1}).
$$
Thus, for large $j$, $\log \lam_j \ge 1.9 \log \lam_{j-1}$ and hence
$\sum_j 1/\log \lam_{j-1}$ converges.  Now,
$\mu_j = j + O(1)$ and for $r>s\ge 1$
\begin{equation}
\mu_r -\mu_s = r-s + O\biggl( \sum_{j\ge s} 2^{-\mu_{j-1}} \biggr) =
r-s + O(2^{-s}).
\end{equation} 
Therefore, the sequence $(\mu_j-j)$ is a Cauchy sequence converging to
some value $c_3$, and $|\mu_j-j-c_3| = O(2^{-j})$.
\Endproof
\vskip4pt

For a vector $\bb=(b_1,\dots,b_h)$ of non-negative integers,
let $\AA(\bb)$ be the set of square-free integers $a$ composed of exactly 
$b_j$ prime factors in $D_j$ for each $j$. 
Denote $k=b_1+\cdots+b_h$.  
For the remainder of this section, $M$ will be a sufficiently large
absolute constant, which we take to be an even integer.

\begin{lem}\label{sumW}
Suppose $\sg>0${\rm ,} $\bb=(b_1,\dots,b_h)$ and define $m=\min \{ j: b_j \ge
1\}$. If $\sg<1${\rm ,} further assume that $m\ge M$ and $b_j \le 2^{j/2}$
for each $j$.  Then
$$
\sum_{a\in \AA(\bb)} \frac{W(a;\sg)}{a} \le \frac{(2\log 2)^k}{b_m! 
  \cdots b_h!} \left[\! 1.01 + 2^{c_5} \sg
  \sum_{j=m}^h 2^{-j+b_m + \cdots + b_j} \right],
$$
where $c_5$ is an absolute constant.
\end{lem}

We next apply Lemma \ref{sumW} for many vectors $\bb$.

\begin{lem}\label{lowvol}
Suppose $0<\a\le 1${\rm ,} $y\ge y_0(\a)$ and $0<\sg \le 2^{-2M-1/\a}\log y$.  Define
\begin{align*}
v&=\fl{\frac{\log\log y-\max(0,\log \sg)}{\log 2} - 2M - 1/\a + 1}, \\
s&=M+\max\(0,\fl{\frac{\log \sg}{\log 2}} \) -c_5-10-\frac{\log \sg}{\log 2}.
\end{align*}
Suppose $k\ge M+1$.  Then{\rm ,} for some subset $\AA$ of the squarefree integers
$a \le y^{\a}$ satisfying $P^-(a) > e^{\sg}$ and $\omega(a)=k${\rm ,} we have
$$
\sum_{a\in \AA}
\frac{3\tau(a)-2W(a;\sg)}{a} \ge \frac13 (2v\log 2)^k \Vol(Y_k(s,v)),
$$
where $Y_k(s,v)$ is the set of $\bx=(\xi_1,\dots,\xi_k)\in \RR^k$ satisfying
\begin{enumerate} \isep
\item $0 \le \xi_1 \le \cdots \le \xi_k < 1${\rm ;}
\item For $1\le i\le \sqrt{k-M}${\rm ,} $\xi_{M+i^2} > i/v$ and
$\xi_{k+1-(M+i^2)} < 1-i/v${\rm ;}
\item $\sum_{j=1}^k 2^{j-v\xi_j} \le 2^s$.
\end{enumerate}
\end{lem}

Condition (ii) in the definition of $Y_k(s,v)$ is very mild and does
not affect the volume very much.  It arises from taking numbers $a\in \AA$
which  have neither too many small prime factors nor
too many large prime factors.
Lemmas \ref{sumW} and \ref{lowvol} will be proved in Section~\ref{sec:lowtovol}.
The volume of $Y_k(s,v)$ can be estimated using bounds on uniform
order statistics (\S \ref{sec:uos}).  Once the volume has been bounded
from below, we are in position to complete the lower bounds in
Theorems \ref{thm1}--\ref{Hrbounds}.

%
%

\begin{lem}\label{volYnsv}
Suppose $v\ge 1${\rm ,} $10M \le k\le 100(v-1)${\rm ,} $s\ge M/2+1$ and 
$0\le k-v \le s-M/3-1$. 
Then
$$
\Vol(Y_k(s,v)) \gg \frac{k-v+1}{(k+1)!}.
$$
\end{lem}

Lemma \ref{volYnsv} will be proved in Section~\ref{sec:lowvol}.

%
%
%
%
\section{Proof of Theorems \ref{thm1}, \ref{thm:shorts},
  \ref{thm:squarefree}, \ref{H1bounds} and \ref{Hrbounds}}
  \label{sec:theorems}
%
%
%

Suppose throughout that $z\ge y+1$ and $y\ge y_0$ (Theorem \ref{thm1}
is trivial if $y<y_0$).   

\demo{Upper bounds in  Theorems {\rm \ref{thm1}, \ref{thm:shorts}} and
  {\rm \ref{thm:squarefree}} when $y\le \sqrt{x}$} 
If $0<\eta\le 1$ and $\Delta\ge \sqrt{x}$, \pagebreak then
\begin{align*}
H(x,y,z)-H(x-\Delta,y,z) &\le \sum_{y<d\le z} \(\frac{\Delta}{d}+1\) \\*
&\ll \eta \Delta + (z-y) \\*
&\ll \eta \Delta.
\end{align*}
This proves the upper bounds in the three theorems when $z\le z_0(y)$.
For $z_0(y) \le z \le y^{1.001}$, the desired bounds follow from Lemmas
\ref{Hupper}, \ref{Hupper2} and \ref{Sstarfinal}. 
When $\b \gg 1$, our upper bound coincides with that of
Theorem T1 (ii).
When $z\ge y^{1.001}$, the trivial bound
$H(x,y,z)-H(x-\Delta,y,z)\le \Delta+1$ suffices.

\demo{Lower bounds in  Theorems {\rm \ref{thm1}, \ref{thm:shorts}} and
  {\rm \ref{thm:squarefree}} when $y\le \sqrt{x}$} Assume $\frac{x}{\log^{10} z} \le \Delta \le x$ for the
estimation of
$H(x,y,z)-H(x-\Delta,y,z)$ and $\frac{x}{\log y}
\le \Delta \le x$ in the estimation of $H^*(x,y,z)-H^*(x-\Delta,y,z)$.

If $0<\eta \le \frac{1}{\log^{20} y}$, then
\begin{align*}
H(x,y,z)-H(x&-\Delta,y,z) \ge \sum_{y<d\le z}
  \fl{\frac{x}{d}}-\fl{\frac{x-\Delta}{d}} \\[3pt]
&\qquad\qquad\quad - \sum_{y<d_1<d_2\le z}
  \fl{\frac{x}{\lcm[d_1,d_2]}}-\fl{\frac{x-\Delta}{\lcm[d_1,d_2]}} \\[3pt]
&\ge \Delta \( \sum_{y<d\le z} \frac{1}{d} - \sum_{y<d_1<d_2\le z}
  \frac{1}{\lcm[d_1,d_2]} \) - 2(z-y+1)^2.
\end{align*}
Let $m=(d_1,d_2)$, so that $m\le z-y$.  Write $d_1=t_1 m$, $d_2=t_2 m$.  Then
\begin{align*}
H(x,y,z)-H(x-\Delta,y,z) 
&\ge \Delta \( \sum_{y<d\le z} \frac{1}{d} - \sum_{m\le z-y}
  \frac{1}{m} \sum_{\frac{y}{m} < t_1 < t_2 \le \frac{z}{m}}
  \frac{1}{t_1 t_2} \) \\[3pt] & \quad -  O(\eta^2 y^2) \\[3pt]
&= \Delta \sum_{y<d\le z} \frac{1}{d} - O(\Delta \eta^2 \log y +
  \eta^2 x) \\[3pt]
&= \Delta \sum_{y<d\le z} \frac{1}{d} - O( \eta \Delta (\log y)^{-10}).
\end{align*}
The sum on $d$ is $\gg \eta$, and we conclude the lower bound in
Theorem 1 (v) and Theorem 2 for this range of $y,z$.

Now suppose $g>0$, $y\ge y_0(g)$, $0< \eta \le \frac{1}{\log^{20} y}$ and
there are $\ge g(z-y)$ square-free
integers in $(y,z]$.  By a theorem of Filaseta and Trifonov
 \cite{FiTr}, this last condition holds unconditionally 
with $g=\frac12$ provided
$z\ge y + K y^{1/5}\log y$ for a large constant $K$.
We obtain
\begin{align*}
H^*(x,y,z)-H^*(x-\Delta,y,z) &\ge \sum_{y<d\le z} \;\;
 \sum_{\substack{x-\Delta < e
  \le x \\ d|e}} \mu^2(e) - \sum_{y<d_1<d_2\le z}
  \frac{x}{\lcm[d_1,d_2]} \\*
&\ge \sum_{y<d\le z} \mu^2(d) \;\;
  \sum_{\substack{\frac{x-\Delta}{d} < f \le \frac{x}{d} \\ (f,d)=1}}
  \mu^2(f) - O(\eta^2 x \log y).
\end{align*}
A simple elementary argument yields
$$
\sum_{\substack{f\le w \\ (d,f)=1}} \mu^2(f) = C_d w + O\(w^{1/2}\tau(d)\),
$$
where
$$
C_d = \frac{\phi(d)}{d} \prod_{p\nmid d} (1-1/p^2).
$$
Thus,
$$
H^*(x,y,z)-H^*(x-\Delta,y,z) \gg \frac{\Delta}{y^2} \sum_{y<d\le z}
\mu^2(d) \phi(d) - O(\eta \Delta (\log y)^{-18}).
$$
Now apply the estimate
\begin{equation}\label{SR}
\sum_{n\le x} \frac{1}{\phi(n)} = C_1 \log x + C_2 + O\pfrac{(\log x)^{2/3}}
{x}
\end{equation} 
due to Sitaramachandra Rao \cite{Rao}, where $C_1,C_2$ are certain
constants (Landau had in 1900 proved a weaker version with error term
$O(\frac{\log x}{x})$).  By the Cauchy-Schwarz inequality and our
assumption, 
$$
\sum_{y<d\le z} \mu^2(d) \phi(d) \ge \biggl( \sum_{y<d\le z} \mu^2(d)
\biggr)^2 \biggl( \sum_{y<d\le z} \frac{1}{\phi(d)} \biggr)^{-1} \gg
\eta y^2.
$$
We conclude that
$$
H^*(x,y,z)-H^*(x-\Delta,y,z) \gg \eta \Delta,
$$
which completes the proof of the lower bound in Theorem
\ref{thm:squarefree} in this case.

Next, suppose $\frac{1}{\log^{20} y} \le \eta \le \frac{1}{100}$ and
define $\b,\xi$ by \eqref{etabetaxi}.  Let $\sg=\eta$,
 $0<\a\le 1$, $g\ge 1$ and $y\ge y_0(\a,g)$.
In Lemma \ref{lowvol}, we have
\begin{align*}
v &= \fl{\frac{\log\log y}{\log 2} -2M-1/\a+1}, \\
s &= M+c_3-21-\frac{\log \eta}{\log 2} \ge \frac{M}{2}+1-\frac{\log \eta}
{\log 2}.
\end{align*}
We will apply Lemmas \ref{lowvol} and \ref{volYnsv} with all $k$ satisfying
$$
(1+\b/100)v \le k \le \min(1+\b,\log 4) v.
$$
This includes at least one value of $k$ since $\frac{\log 100}{\log\log y} \le \beta \le 20$.  Also, by
\eqref{etabetaxi},
$$
k-v \le \b v = \frac{-\log \eta}{\log\log y} v \le s-M/3-1,
$$
so that all hypotheses of Lemma \ref{volYnsv} are satisfied.  
For each such $k$ we obtain
$$
\sum_{a\in \AA} 
\frac{3\tau(a)-2W(a;\eta)}{a} \gg_{g,\a} \b \frac{(2v\log 2)^k}{k!},
$$
for some subset $\AA$ of the squarefree integers $a\le y^{\a}$ with
$\omega(a)=k$ and $P^-(a) > e^{\eta}$.  By Lemma \ref{Ilow},
\begin{equation}\label{Isum2}
\sum_{\substack{a\le y^{\a} \\ \omega(a)=k \\ \mu^2(a)=1}}
\frac{I(a;\eta)^g}{a} \gg_{g,\a} \b \frac{(2^g v\log 2)^k}{k!}.
\end{equation} 
When $g=1$, Lemma \ref{Nortonlem} gives
\begin{equation}\label{lowerbd5}
\frac{\eta}{\log^2 y} \sum_{\substack{a\le y^{\a} \\ \mu^2(a)=1}}
\frac{I(a;\eta)}{a} \gg_\a \frac{\b}{\max(1,-\xi) (\log y)^{G(\b)}}.
\end{equation} 
By Lemma \ref{Hlower2} and \eqref{lowerbd5} with $\a=\frac{1}{10}$, we obtain
\begin{equation}\label{lowerbd6}
\begin{split}
H(x,y,z)-H(x-\Delta,y,z) &\ge H^*(x,y,z)-H^*(x-\Delta,y,z) \\
&\gg  \frac{\b \Delta}{\max(1,-\xi) (\log y)^{G(\b)}}.
\end{split}
\end{equation} 

Next, let $0<\a\le 1$,  $\g=2^{-20M-1/\a}$,  
and suppose that $\frac{1}{100} \le \eta \le \g \log y$.
Let $g\ge 1$, $y\ge y_0(\a,g)$ and assume $\sg=\eta$ or $\sg=2\eta$.
In Lemmas \ref{lowvol} and
\ref{volYnsv}, $v \ge 18M$ and $s\ge M/2+1$ if $M$ is large enough.
Using the single term $k=v$, we have (note $e^\eta=y^u$)
$$
\sum_{a\in \AA}
 \frac{3\tau(a)-2W(a;\sg)}{a} \gg_{g,\a} \frac{(2v\log 2)^v}{v\cdot v!}\gg 
\frac{u^{\del-2}}{(-\log u)^{3/2}},
$$
for a subset $\AA$ of the squarefree integers $a\le y^{\a}$ satisfying
$\omega(a)=k$ and $P^-(a) > e^{\sg}$.
Lemma \ref{Ilow} then gives
\begin{equation}\label{Isum}
\sum_{\substack{a\le y^{\a} \\ P^-(a) > e^{\sg} \\ \mu^2(a)=1}}
\frac{I(a;\sg)^g}{a} \gg_{g,\a} \frac{2^{v(g-1)}
u^{\del-2}}{(-\log u)^{3/2}}
\gg_g \frac{u^{\del-1-g}}{(-\log u)^{3/2}}.
\end{equation} 
By Lemma \ref{Hlower2} and \eqref{Isum} with $g=1$ and $\a=\frac{1}{10}$,
\begin{equation}\label{lowerbd2}
H^*(x,y,z)-H^*(x-\Delta,y,z) \gg \Delta \frac{u^\del}{(-\log u)^{3/2}}.
\end{equation} 

Finally, if $z\ge y^{1+\g}$, then by \eqref{lowerbd2}
$$
H(x,y,z)-H(x-\Delta,y,z)\ge H(x,y,y^{1+\g})-H(x-\Delta,y,y^{1+\g}) \gg \Delta
$$
and
$$
H^*(x,y,z)-H^*(x-\Delta,y,z)\ge
H^*(x,y,y^{1+\g})-H^*(x-\Delta,y,y^{1+\g}) \gg \Delta.
$$

Corollary \ref{cor1} also follows in the case $y_i\le \sqrt{x_i}$ ($i=1,2$).

\demo{Proof of Theorem {\rm \ref{thm1} (vi)}} 
Assume throughout that $y>\sqrt{x}$ and $y+1 \le z \le x$. 
If $\frac{x}{y} < \frac{x}{z} + 1$ and $y<d_1<d_2\le z$, then
$$
\lcm[d_1,d_2] = \frac{d_1 d_2}{(d_1,d_2)} >
\frac{y(y+(d_1,d_2))}{(d_1,d_2)} \ge y + \frac{y^2}{z-y} \ge x.
$$
Hence
$$
H(x,y,z) = \sum_{y<d\le z} \fl{\frac{x}{d}} \asymp \eta x.
$$

Now assume $\frac{x}{y} \ge  \frac{x}{z} + 1$.  
If $y>x/y_0$, then $H(x,y,z)\gg \fl{z}-\fl{y}\gg  x$\break\vskip-11pt\noindent since 
$z\ge \frac{xy}{x-y} \ge \frac{y_0}{y_0-1}y$.
Thus 
$$
H(x,y,z) \asymp x \asymp H(x,\tfrac{x}{z},\tfrac{x}{y}).
$$

Next suppose $\sqrt{x} < y \le x/y_0$ and $\eta \le \frac{1}{\log^2 (x/y)}$.
Here, 
$$
H(x,y,z) = \sum_{y<d\le z} \fl{\frac{x}{d}} + O\(
\sum_{\substack{y<d_1<d_2\le z \\ \lcm[d_1,d_2]\le x}}
\frac{x}{\lcm[d_1,d_2]} \).
$$
The sum on $d$ is $\asymp \eta x$.  Writing
$m=(d_1,d_2)$, the big-$O$ term is
$$
\ll \sum_{y^2/x < m \le z-y} \frac{1}{m} \sum_{\frac{y}{m}<t_1<t_2 \le
    \frac{z}{m}} \frac{1}{t_1 t_2} \ll \eta^2 \log(x/y)\ll
  \frac{\eta}{\log(x/y)}. 
$$
Since $\log(x/y) \asymp \log(x/z)$,
$$
H(x,y,z) \asymp \eta x \asymp  H(x,\tfrac{x}{z},\tfrac{x}{y})
$$
in this case.

Lastly, suppose $\eta \ge \frac{1}{\log^2 (x/y)}$.
Partition $(\frac{x}{\log^2 (x/y)},x]$ into intervals $(x_1,x_2]$, where
$$
x_2-x_1 \in \left[ \frac{x_2}{\log^3 (x/y)}, \frac{2x_2}{\log^3
    (x/y)}\right]. 
$$
We have, for $n$ lying in such an interval $(x_1,x_2]$,
$$
\tau(n,\tfrac{x_2}{z},\tfrac{x_1}{y})\ge 1 \; \implies \; 
\tau(n,y,z) \ge 1 \; \implies \tau(n,\tfrac{x_1}{z},\tfrac{x_2}{y})\ge 1.
$$
We obtain the upper bound
$$
H(x,y,z) \le \frac{x}{\log^2 (x/y)} + 
\sum_{x_1,x_2} H(x_2,\tfrac{x_1}{z},\tfrac{x_2}{y}) - 
H(x_1,\tfrac{x_1}{z},\tfrac{x_2}{y}).
$$
For large enough $y_0$,
$$
\log \pfrac{x_2}{y} \ge \log \pfrac{x}{y\log^2(x/y)} \ge \frac12
\log\pfrac{x}{y},
$$
and thus 
$$
x_2-x_1 \ge \frac{x_2}{\log^{4} (x_2/y)}.
$$
Also, 
$$
\log\pfrac{x_2/y}{x_1/z} \asymp \eta, \qquad \frac{x_1}{z} \le
\frac{x_1}{\sqrt{x}} \le \sqrt{x_2},
$$
so that  by Theorem \ref{thm:shorts} and the part of Corollary \ref{cor1}
already proved,
$$
 H(x_2,\tfrac{x_1}{z},\tfrac{x_2}{y}) - H(x_1,\tfrac{x_1}{z},\tfrac{x_2}{y})
\ll \frac{x_2-x_1}{x_2} H(x_2,\tfrac{x_1}{z},\tfrac{x_2}{y})
\ll \frac{x_2-x_1}{x} H(x,\tfrac{x}{z},\tfrac{x}{y}).
$$
Summing over intervals $(x_1,x_2]$ gives the desired upper bound.
The lower estimate is obtained in the same way starting from
$$
H(x,y,z) \ge 
\sum_{x_1,x_2} H(x_2,\tfrac{x_2}{z},\tfrac{x_1}{y}) - 
H(x_1,\tfrac{x_2}{z},\tfrac{x_1}{y}).
$$
This completes the proof of Theorems \ref{thm1}, \ref{thm:shorts}, and
  \ref{thm:squarefree}.

\demo{Proof of Theorems {\rm \ref{H1bounds}} and {\rm \ref{Hrbounds}} when $z\le y^{C}$}
Suppose $y_0(r)+1 \le y+1 \le z \le x^{5/8}$.
When $z\le e^{100rC} y$, the trivial bound $H_r(x,y,z) \le H(x,y,z)$
suffices for an upper bound.
For $e^{100rC} y \le z \le y^{C}$, the desired upper bound follows from
Lemmas \ref{Hrupper} and \ref{Sstarfinal}, plus the lower bound for $H(x,y,z)$
given in Theorem \ref{thm1}.

If $r=1$ and $\eta \le \frac{1}{\log^2 y}$, arguing as in the lower
bounds for $H(x,y,z)$, we have
$$
H_1(x,y,z)\ge \sum_{y<d\le z} \fl{\frac{x}{d}} - O(\eta^2 x\log y) 
\gg \eta x \gg H(x,y,z).
$$
If $r\ge 1$ and $\frac{1}{\log^2 y} \le \eta \le \frac{1}{100}$,
combining \eqref{Isum2} (taking $\a=\frac{c}{300r}$ and $g=r$) 
and Lemma \ref{Hrlower1} with $c'=c/3$ gives
\begin{equation}\label{Hrlower5}
H_r(x,y,z) \gg_{r,c} \frac{\b \eta^r x}{(\log y)^{r+1}} 
\sum_{(1+\frac{\b}{100})v \le k \le \min(1+\b,\log 4)v} 
\frac{(2^r v\log 2)^k}{k!}.
\end{equation} 
When $r=1$, Lemma \ref{Nortonlem} gives
$$
H_1(x,y,z) \gg_{c} \frac{\b x}{\max(1,-\xi) (\log y)^{G(\b)}}.
$$
When $r\ge 2$, $z\ge z_0(y)$ and $\eta \le \frac{1}{100}$,
the sum on $k$ in \eqref{Hrlower5}
is dominated by the term $k=\fl{(1+\b)v}$,
whence by Stirling's formula, \eqref{Gdef} and Theorem 1,
\begin{align*}
H_r(x,y,z) &\gg_{r,c} \frac{\b x}{(\log\log y)^{1/2} (\log
  y)^{G(\b)}} \\
&\gg_{r,c} H(x,y,z) \frac{\max(1,-\xi)}{\sqrt{\log\log y}}.
\end{align*}

Next, let $\a=\frac{c}{300r}$, $\g=2^{-2M-1/\a}$ and suppose
$\frac{1}{100} \le \eta \le \g \log y$ and $c'=c/3$.
When $\eta \ge 1000r \cdot 3^{2r}$, apply Lemma \ref{Hrlower2} and the
$g=1$ case of \eqref{Isum}.  Otherwise apply Lemma  \ref{Hrlower1} and the
$g=r$ case of \eqref{Isum}.  In either case, we obtain
\begin{equation}\label{lowerbd3}
\begin{split}
H_r(x,y,z) &\gg_{r,c} x  
\frac{u^\del (\log(2+\eta))^{\nu(r)+1}}{\eta (-\log u)^{3/2}} \\
&\gg_{r,c} H(x,y,z) \frac{(\log(2+\eta))^{\nu(r)+1}}{\eta}.
\end{split}
\end{equation} 
Note that by \eqref{Gdef}, $G(\b)=\del+O(1/\log\log y)$
for $1.01y \le z \le 2y$, and so $(\log y)^{G(\b)} \asymp
(\frac{\log y}{\eta})^\del$ in this range.
The desired lower bound for $H_r(x,y,z)$ now follows from
\eqref{lowerbd3} and the upper bound for $H(x,y,z)$ in Theorem
\ref{thm1}.
When $y^{1+\g} \le z \le \min(y^{C},x^{1/2-c/3})$ (for $r=1$, take $C=10$), 
we take the $h=1$ term in the sum
in Lemma \ref{Hrlower2} (i), obtaining
\begin{equation}\label{lowerbd1}
H_r(x,y,z) \gg_{r,c,C} \frac{x (\log \log y)^{\nu(r)+1}}{\log y}.
\end{equation} 

Finally, when $x^{1/2-c/3} < z \le x^{5/8}$, $z\le y^C$ 
and $yz\le x^{1-c}$, the desired lower bound is
given in Lemma \ref{Hrlower2} (ii) since $y\le x^{1/2-2c/3}$ and thus
$\min(y,z/y) \ge y^{c/3}$.

\demo{Proof of Theorem {\rm \ref{Hrbounds} \eqref{Prhugez}}} 
Apply Lemma \ref{Hrlower2} (ii) with $c=\frac{1}{16}$.

\demo{Proof of Theorem {\rm \ref{H1bounds}} when $y^{10} \le z \le
  x^{5/8}$} 
This proof is quite simple and does not depend on the results of
Sections \ref{sec:upper} and \ref{sec:lower}.
In this range, $H(x,y,z)\break \gg x$.
Write each $n$ with $\tau(n,y,z)=1$ in the form 
$$
n=klm, \qquad P^+(k) \le y, \qquad l\in \PP(y,z), \qquad  P^-(m)>z.
$$
If $p^2|l$ for some prime $p$, then $p>\sqrt{z}$ and thus the number
of such $n$ is $\ll x/\sqrt{z}$.  
Otherwise, $l=1$ or $l$ is prime.  Also
$k\le y^2$, for otherwise $k$ has at least 2 divisors in $(y,z]$.
Thus, $kl\le y^2 z \le z^{3/2} \le x^{15/16}$.

The number of $n$ with $m=1$ is $\le x^{15/16}$.
Now suppose $m>1$.  For each $k,l$, $x/kl \ge x^{1/16} \ge z^{1/16}$.
By Lemma \ref{Phi}, the number of $m$ is $\ll x/(kl \log z)$.  
Clearly $k$ and $l$ can't both
be 1.  If $k=1$, then $l$ is prime and by Lemma \ref{sum1p},
the number of such $n$ is
$$
\ll\frac{x}{\log z} \sum_{y<p\le z} \frac{1}{p} \ll
\frac{x (\log\log z-\log\log y)}{\log z}.
$$
If $l>1$ and $k>1$, then $k\le y$ and also $z/P^-(k) < l \le z$.
For a given $k$, the sum over $l$ of
$1/l$ is $\ll \frac{\log P^-(k)}{\log z}$.  Also,
$$
\sum_{2\le k\le y} \frac{\log P^-(k)}{k} \le
\sum_{p\le y} \frac{\log p}{p} \sum_{k'\in \PP(p-1/2,y)} \frac{1}{k'}
\ll (\log y) \log\log y.
$$
Thus, the number of such $n$ is
$$
\ll \frac{x (\log y)\log\log y}{\log^2 z} \ll \frac{x \log\log y}{\log z}.
$$
The last case to consider is $l=1$ and $k>1$.  Here $k/P^-(k) \le y < k$.
and we write $k=pk'$, where $p=P^-(k)$,
$P^-(k') \ge p$ and $y/p < k' \le y$.  By Lemma~\ref{Phi} and partial
summation, the number of such $n$ is
$$
\ll \frac{x}{\log z} \sum_{p\le y} \frac{1}{p} \sum_{\substack{P^-(k')\ge p \\
y/p < k' \le y}} \frac{1}{k'}
\ll \frac{x}{\log z} \sum_{p\le y} \frac{1}{p} \ll \frac{x\log\log y}{\log z}.
$$
Putting these estimates together proves the upper bound.

For the lower bound, first note that if $kl \le \frac{x}{2z}$, then by 
Lemma \ref{Phi}, the number of $m$ is $\gg \frac{x}{kl \log z}$.
The number of $n$ with $k=1$, $l$ a prime in
$(y,z^{1/3})$ and $P^-(m) > z$ is
$$
\gg \frac{x}{\log z} \sum_{y<p\le z^{1/3}} \frac{1}{p} \gg
 \frac{x(\log\log z-\log\log y)}{\log z}.
$$
Next, let $l=1$ and put $k=ph$, where $10< p\le y^{1/4}$, $P^-(h) > p$
and $y/p < h\le y$.  The number of such $n$ is, by Lemma \ref{Phi}
and partial summation,
\begin{align*}
&\gg \frac{x}{\log z} \sum_{10 < p\le y^{1/4}} \frac{1}{p} \sum_{\substack{
j\in \ZZ \\ y/p < 2^j \le y/2}} \frac{\Phi(2^{j+1},p)-\Phi(2^j,p)}{2^j} \\
&\gg \frac{x}{\log z} \sum_{10 < p\le y^{1/4}} \frac{1}{p} 
\gg \frac{x\log\log y}{\log z}.
\end{align*}
This completes the lower bound.

\demo{Remarks} By extending the methods used to prove Theorem
\ref{H1bounds} when $y^{10}\le z \le x^{5/8}$, it should be possible to
determine the order of 
$H_r(x,y,z)$ for $y^{10}\le z\le x^{1/2}$ for any fixed $r\ge 2$.
We conjecture that for each $r\ge 1$ and $y^{10} \le z \le x^{1/2}$, 
$$
\frac{H_r(x,y,z)}{x} \asymp \frac{Q_r(\log\log y,\log\log z)}{\log z}
$$
for some polynomial $Q_r$.  
 
\section{Initial sums over $L(a;\sg)$ and $L_s(a;\sg)$}\label{sec:initial}
 
The object of this section is to prove Lemmas \ref{Hupper} and
\ref{basiclower}.
The upper bounds are more complicated,
due to having to count integers with $d|n$, $y<d\le z$ and
$P^+(d)=z^{o(1)}$.
It is convenient to work with divisors $d|n$ with $P^+(d) < P^+(n)$.
If this is not the case, then the complementary divisor $n/d$ does
satisfy $P^+(n/d) < P^+(n)$ (here we assume that $P^+(n)^2 \nmid n$,
the number of exceptions being very small).
We put $n$ into a very short interval, so that when $y<d\le z$, $n/d$
lies in an interval $(y',z']$ with $\log(z'/y') \approx \log(z/y)$.

%
%

\begin{lem}\label{Hshort}
Suppose $y,z,x_1,x_2$ are positive real numbers satisfying
$$
100 \le y< z = e^{\sg}y \le x_1^{3/4}, \quad z\le y^{\log\log y},
$$
and
$$
\sg \ge \frac{1}{10\log^2 z}, \quad 
\frac{3^{10} x_1}{\log^{10} z} \le x_2 -x_1 \le \frac{x_1}{\log^4 z}.
$$
Then
$$
H(x_2,y,z)-H(x_1,y,z) \ll (x_2-x_1) \bigl[ S(y;\sg) + S(x_2/z;\sg) \bigr].
$$
If in addition $r\ge 1$ and $y\ge y_0(r)${\rm ,} then
$$
H_r(x_2,y,z)-H_r(x_1,y,z) \ll_r (x_2-x_1) \sum_{\substack{1\le s\le r \\ 
\nu(s)\le \nu(r)}} \bigl[ S_s(y;\sg) + S_s(x_2/z;\sg) \bigr].
$$
\end{lem}

\Proof  Let $\AA$ be the set of integers $n\in (x_1,x_2]$ satisfying
\begin{enumerate} \isep
\item $\tau(n,y,z)\ge 1$;
\item $\tau(n,x_1/z,x_2/z)=0$;
\item if $p$ is prime, $p|n$ and $p>\log^{10} z$, then $p^2 \nmid n$;
\item if $d|n$ with $y<d\le z$, then $P^+(d) > \log^{20} z$ and
$P^+(n/d)>\log^{20} (x_2/y)$.
\end{enumerate}
Let $\AA_r$ be the set of $n\in \AA$ with $\tau(n,y,z)=r$.
Since $x_2-x_1 \ge x_1/z$, the number of integers in
$(x_1,x_2]$ not satisfying (ii) is at most
$$
\sum_{x_1/z \le d \le x_2/z} \( \frac{x_2-x_1}{d}+1 \) \ll \frac{(x_2-x_1)^2}
{x_1} \ll \frac{x_2-x_1}{\log^4 z}.
$$
The number of integers in $(x_1,x_2]$ failing (iii) is
$$
\le \sum_{\log^{10} z < p \le \sqrt{x_2}} \( \frac{x_2-x_1}{p^2} + 1 \)
\ll \frac{x_2-x_1}{\log^{10} z}.
$$
Put $y_1=y$, $z_1=z$, $y_2=x_2/z$ and $z_2=x_2/y$.
Every $n\in (x_1,x_2]$ satisfying (i) and (ii) can be written in the form
\begin{equation}\label{nm1m2}
n=m_1 m_2, \qquad y_j < m_j \le z_j \; (j=1,2).
\end{equation} 
Note that $z_2 \ge x_1^{1/4} \ge z_1^{1/3}$ and so
 $y_j \ge z_j^{1/(3\log\log z_j)}$
for $j=1,2$.  By Lemma \ref{Psixyz}, the number of integers $\le x_2$
not satisfying (iv) is
\begin{align*}
\le \sum_{j=1}^2 \sum_{\substack{m_j>y_j \\ P^+(m_j)\le \log^{20} z_j}}
\frac{x_2}{m_j}& \ll x_2 \sum_{j=1}^2 (\log\log z_j) e^{-\frac{\log y_j}
{80\log\log z_j}} \\
&\ll \frac{x_2}{\log^{100} z} \ll \frac{x_2-x_1}{\log^{10} z}.
\end{align*}
We conclude that
\begin{equation}\label{HHr1}
\begin{split}
H(x_2,y,z)-H(x_1,y,z) &\le |\AA| + O\pfrac{x_2-x_1}{\log^4 z}, \\
H_r(x_2,y,z)-H_r(x_1,y,z) &\le |\AA_r| + O\pfrac{x_2-x_1}{\log^4 z}.
\end{split}
\end{equation} 

We follow similar procedures for bounding $|\AA|$ and $|\AA_r|$.  First,
writing $n\in \AA$ in the form \eqref{nm1m2}, we see by (iv) that
$P^+(m_j) > \log^{20} z_j$ for $j=1,2$.  By (iii), $P^+(m_1) \ne
P^+(m_2)$.  Suppose that
\begin{equation}\label{pdef}
p = P^+(m_j) < P^+(m_{3-j}) \qquad (j=1 \text{ or } j=2),
\end{equation} 
and write
\begin{equation}\label{nrep1}
n=abp, \quad P^+(a) < p < P^-(b), \quad b>p.
\end{equation} 
We have $\tau(a,y_j/p,z_j/p) \ge 1$, which implies
$$
\log(y_j/p) \in \LL(a;\sg).
$$

Each $n\in \AA_r$ may be written uniquely as
$$
n=m_{11} m_{12} = m_{21} m_{22} = \cdots = m_{r1} m_{r2},
$$
where $y_1 < m_{11} < \cdots < m_{r1} \le z_1$ and
$y_2 < m_{r2} < \cdots < m_{12} \le z_2$.  There may be divisors of $n$
lying in $(y_1,z_1] \cap (y_2,z_2]$, in which case $m_{i1}=m_{i'2}$ for
some pairs $i,i'$.  Let $p_{ij} = P^+(m_{ij})$ for each $i,j$.
Since $p_{i1} \ne p_{i2}$ for each $i$, $p_{ij}\ne P^+(n)$ for exactly $r$
pairs of indices $i,j$.
Therefore, there is an integer $s$ with
$1\le s\le r$ and $\nu(s)\le \nu(r)$, a $j\in \{1,2\}$ and
a prime $p<P^+(n)$, so that
exactly $s$ of the primes $p_{1j},\dots, p_{rj}$ are equal to $p$.
Writing $n$ in the form \eqref{nrep1}, we see that
$\tau(a,y_j/p,z_j/p)=s$, and thus 
$$
\log(y_j/p) \in \LL_s(a;\sg).
$$

By (iv) and \eqref{nrep1},
\begin{equation}\label{bp}
b>p>\log^{20} z_j \ge \log^{20} (z^{1/3})
\end{equation} 
and so 
$$
\frac{x_2-x_1}{ap} \ge \frac{x_2-x_1}{x_2} p \ge 
\frac{p}{\log^{10} (z^{1/3})} > p^{1/2}.
$$
By Lemma \ref{Phi}, given $a$ and $p$, the number of possible $b$ is at most
$$
\Phi\(\frac{x_2}{ap},p\) - \Phi\(\frac{x_1}{ap},p\)
\ll \frac{x_2-x_1}{ap\log p}. 
$$
Let $Q_j(a) := \max(\log^{20} z_j, P^+(a))$, so that $p>Q_j(a)$.
By \eqref{nrep1} and \eqref{bp},
\begin{equation}\label{HHr2}
\begin{split}
|\AA| &\ll (x_2-x_1) \sum_{j=1}^2 \sum_{P^+(a)\le z_j} \frac{1}{a}
\sum_{\substack{\log(y_j/p)\in \LL(a;\sg) \\ p > Q_j(a)}} \frac{1}{p\log p}, \\
|\AA_r| &\ll (x_2-x_1) \sum_{j=1}^2 \sum_{\substack{1\le s\le r \\
\nu(s) \le \nu(r)}} \sum_{P^+(a)\le z_j} \frac{1}{a}
\sum_{\substack{\log(y_j/p)\in \LL_s(a;\sg) \\ p > Q_j(a)}} \frac{1}{p\log p}.
\end{split}
\end{equation} 
Let $\LL=\LL(a;\sg)$ or $\LL=\LL_s(a;\sg)$, as appropriate.  Then
$$
\log(y_j/p) \in \SS := \LL \cap \left[-\sg, \log\pfrac{y_j}{Q_j(a)} \right).
$$
The set $\LL$ is the disjoint union of
intervals, each with a left endpoint  $-\sg + \log d$ for
some $d|a$ or a right endpoint  $\log d$ for
some $d|a$.  Therefore,
$\SS$ has the same property with the possible exception of
one interval whose right endpoint is $\log(y_j/Q_j(a))$.
Breaking long intervals into many short ones,
we may partition $\SS$ into
intervals $I_i$, each of length $1/\log^{10} z_j$, and
intervals $I'_i$, each with left endpoint  $-\sg+\log d$ or right 
endpoint  $\log d$ for some $d|a$ (with one possible exception)
and of length $< 1/\log^{10} z_j$.
If $I_i=[A,B)$, then $y_j e^{-B} \ge \log^{20} z_j$.
By Lemma \ref{BrunTitch},
$$
\sum_{\log(y_j/p)\in I_i} \frac{1}{p\log p} \le \frac{\pi(y_j e^{-A})-
\pi(y_j e^{-B})}{y_j e^{-B} \log(y_j e^{-B})}
\ll \frac{B-A}{\log^2(y_j e^{-B})}.
$$
Since $B \le \log a$ and $B\le \log(y_j/P^+(a))$ when $a>1$, we have
$$
\log (y_j e^{-B}) \ge \log \max(y_j/a,P^+(a)).
$$
Adding the contributions of all intervals $I_i$ gives
$$
\sum_{\log(y_j/p)\in \cup I_i} \frac{1}{p\log p} \ll \frac{\text{meas}(\LL)}
{\log^2 (y_j/a + P^+(a))}.
$$
Trivially, for each $i$
$$
\sum_{\log(y_j/p)\in I_i'} \frac{1}{p\log p} \le 
\sum_{\log(y_j/m)\in I_i'} \frac{1}{m} \ll \frac{1}{\log^{10} z_j}.
$$
The number of intervals $I_i'$ is $\le 2\tau(a)+1 \le 3\tau(a)$, and  thus
$$
\sum_{\log(y_j/p)\in \SS} \frac{1}{p\log p} \ll
\frac{\text{meas}(\LL)}{\log^2(y_j/a+P^+(a))} + \frac{\tau(a)}{\log^{10} z_j}.
$$
Next, we sum on $a,j$ and $s$ in \eqref{HHr2} and use
$$
\sum_{P^+(a)\le z_j} \frac{\tau(a)}{a} = \prod_{p\le z_j} \( 1 + 
\frac{2}{p} + \frac{3}{p^2} + \cdots \) \ll \log^2 z_j.
$$
By \eqref{HHr1}, this gives
\begin{align*}
H(x_2,y,z)-H(x_1,y,z) &\ll (x_2-x_1)
\biggl[ \frac{1}{\log^4 z}+S(y_1;\sg)+S(y_2;\sg) \biggr], \\
H_r(x_2,y,z)-H_r(x_1,y,z) &\ll (x_2-x_1)
\biggl[ \frac{1}{\log^4 z} + \!\!\!
\sum_{\substack{1\le s\le r \\ \nu(s)\le \nu(r)}}
\!\!\! \bigl(S_s(y_1;\sg)+S_s(y_2;\sg)\bigr) \biggr].
\end{align*}
Since $L_1(1;\sg) = \sg$, we have 
$$
S(y_1;\sg) \ge S_1(y_1;\sg) \ge \frac{\sg}{\log^2 y_1}
\ge \frac{1}{10\log^4 z},
$$
and this completes the proof.
\Endproof\vskip4pt

\emph{Proof of Lemma} \ref{Hupper}.
Let $x_0=\max(x-\Delta,x/\log^{100} z)$,
partition $(x_0,x]$ into intervals  $(x_1,x_2]$ with
$x_2-x_1 \asymp x_1/\log^{9} z$ and apply Lemma \ref{Hshort} (with
$\sg=\eta$) to each.  For each pair $(x_1,x_2)$,
we have $z \le x^{5/8} \le x_1^{3/4}$, $\sg \ge \frac{1}{\log y} \ge
\frac{1}{\log^2 z}$ and $x_2/z \ge x^{1/3} \ge y^{1/2}$.
Thus
\begin{align*}
H(x,y,z)-H(x-\Delta,y,z) &\ll \frac{x}{\log^{100} z} + \Delta 
  \max_{y^{1/2} \le t\le x} S(t;\eta), \\
H_r(x,y,z) &\ll \frac{x}{\log^{100} z} + x \max_{\substack{1\le s\le r \\
\nu(s) \le \nu(r)}} \max_{y^{1/2} \le t\le x} S_s(t;\eta).
\end{align*}
Since $L_1(1;\eta)=\eta \ge \frac{1}{\log y}$, we have 
$S(y;\eta) \ge S_1(y;\eta) \ge 1/\log^3 y$.
The lemma follows.
\Endproof\vskip4pt

Turning now to the lower bounds, our task is simpler because we may
restrict the integers $n$ in any manner we choose.

\demo{Proof of Lemma {\rm \ref{basiclower}}}
First we derive the lower bound for $H^*(x,y,z)-H^*(x-\Delta,y,z)$.
Consider all square-free integers $n=abp \in (x-\Delta,x]$, where 
\begin{enumerate} \isep
\item $a \le y^{1/8}$;
\item $q|b$ and $q$ prime implies $y^{1/4} < q < y^{7/8}$ or $q>z$;
\item $p$ is prime and $\log (y/p) \in \LL(a;\eta)$.
\end{enumerate}

Condition (iii) ensures that $ap$ has a divisor in $(y,z]$,
and so does $n$  as well.
By (i) and (iii),
$y/z \le y/p \le a \le y^{1/8}$; thus $y^{7/8} \le p \le z$.
Hence each $n$ has a unique representation in the form $abp$.
For fixed $a$ and $p$,
$$
\frac{x}{ap} \ge \frac{x}{y^{1/8}z} \ge y^{1/2}.
$$
Also, $\Delta\ge x/\log^{10} z$, so that
$$
\frac{\Delta}{ap} \ge \frac{x}{ap \log^{10} (x/ap)}.
$$
If $y^{1/2} \le x/(ap) \le 2z$, we count $b$ with $P^-(b)>y^{1/4}$
(so automatically $P^+(b)< y^{7/8}$).
By Lemma \ref{Phi}, the number of such $b$ is $\gg \Delta/(ap\log y)$.
When $x/(ap) >2z$, take $b$ so that $P^-(b) > z$.  
By Lemma \ref{Phi}, the number of such $b$ is $\gg \frac{\Delta}{ap\log z}
\gg \frac{\Delta}{ap\log y}$ as well.
The number of $b$ which are not square-free is at most
$$
\sum_{q > y^{1/4}} \frac{x}{apq^2} \ll \frac{x}{ap y^{1/4}}
\ll \frac{\Delta}{ap y^{1/5}}.
$$
Thus, we obtain
$$
H^*(x,y,z)-H^*(x-\Delta,y,z) \gg \frac{\Delta}{\log y} \sum_{\substack{
a\le y^{1/8} \\ \mu^2(a)=1}} \frac{1}{a}
\sum_{\log(y/p)\in\LL(a;\eta)} \frac{1}{p}.
$$

For each $a$ satisfying (i), $\LL(a;\eta)$ consists of a disjoint
union of
intervals, each with length $\ge \eta$.
Breaking long intervals into shorter ones, we may partition $\LL(a;\eta)$
into intervals each with length between $\frac{1}{2\log^2 y}$ and $\frac{1}
{\log^2 y}$.
If $[v,v+w]$ is one such interval, then $v+w \le \log a \le 
\frac{\log y}{8}$ and thus by Lemma~\ref{sum1p},
$$
\sum_{\substack{v \le \log(y/p) < v+w}} \frac{1}{p} = 
\sum_{\substack{ye^{-v-w} < p \le ye^{-v}}} \frac{1}{p}
\gg \frac{w}{\log y}.
$$
Combining all such intervals gives 
$$ 
\sum_{\log(y/p)\in\LL(a;\eta)} \frac{1}{p} \gg \frac{L(a;\eta)}{\log y},
$$ 
which completes the proof  of the first assertion.

Next we derive the lower bound for $H_r(x,y,z)$.
Consider all integers $n=abp\le x$, where 
\begin{enumerate}
\item $a \le y^{2c}$;
\item $\log(y/p)\in \LL_r(a;\eta)$;
\item $P^-(b) > z$.
\end{enumerate}
We have $y^{1-2c} \le p \le z$, and so each $n$ has a unique
representation as $abp$.  If $d|n$ with $y<d\le z$, then $d=pd'$
for some $d'|a$, thus (ii) ensures that $\tau(n,y,z)=r$.
Since $ap \le z^{1+2c} \le x/4z$, for a given pair $a,p$,
Lemma \ref{Phi} implies that the number of $b$ is $\gg \frac{x}{ap\log z}
\gg \frac{x}{Cap\log y}$.
In contrast with $\LL(a;\eta)$, $\LL_r(a;\eta)$ may not consist only of
intervals of length $\ge \eta$. 
With $a$ fixed, partition $\LL_r(a;\eta)$ as
$$ 
\LL_r(a;\eta) = \(\bigcup_{i=1}^N I_i\) \cup \(\bigcup_{j=1}^M I'_j\),
$$ 
where each $I_i$ is an interval of length $w=1/\log^{3r+7} y$, and 
each $I'_j$ is an interval of length $< w$ and with a left
endpoint $-\eta + \log d$ or a right endpoint
$\log d$ for some $d|a$.  Clearly
$$
\sum_{j=1}^M \text{meas}(I'_j) \le \frac{2\tau(a)}{\log^{3r+7} y}.
$$
Consider one interval $I_i=[v,v+w]$.
Since $v\ge -\eta$, $e^{v+w} \le a \le y^{2c}$ and $y\ge y_0(r,c)$, 
Lemma \ref{sum1p} gives
$$ 
\sum_{\log(y/p)\in I_i} \frac{1}{p} \gg_C \frac{w}{\log y}.
$$  
Adding the contributions of all intervals $I_i$ gives
$$
\sum_{\log(y/p)\in \LL_r(a;\eta)} \frac{1}{p} \ge K
\( \frac{L_r(a;\eta)}{\log y} - \frac{2\tau(a)}{\log^{3r+8} y} \),
$$
where $K$ is a positive constant which depends on $C$.  
Summing over $a$ gives
$$
H_r(x,y,z) \gg_{r,c,C} 
\frac{x}{\log^2 y} \( \sum_{a\le y^{2c}} \frac{L_r(a;\eta)}a
- \frac{2}{\log^{3r+7} y} \sum_{a\le y^{2c}} \frac{\tau(a)}{a} \).
$$
The second sum on $a$ is $O(\log^2 y)$, thus it suffices to show that
the first sum on $a$ is $\gg (\log y)^{-3r-2}$.
By standard prime number estimates,
there is an interval of length $100r \log\log y$ contained in $[\log^3 y,
2\log^3 y]$ which contains  $r$ primes $p_1<\cdots< p_r$.
If $a_0=p_1\cdots p_r$ and $y$ is large enough,
then $\log p_r -\log p_1 \le \eta/10$ and
$\log(p_1p_2)-\log(p_r) \ge \log\log y$.  Hence
$$
L_r(a_0;\eta) \ge \min(\eta,\log\log y) - \log(p_r/p_1) \ge 
\frac{1}{2\log^2 y}.
$$
Therefore
$$
\frac{L_r(a_0;\eta)}{a_0} \gg_r \frac{1}{(\log y)^{3r+2}},
$$
and this completes the proof.
\hfill\qed
 \section{Upper bounds in terms of $S^*(t;\sg)$}\label{sec:Sstar}
 
In this section we prove Lemmas \ref{Hupper2} and \ref{Hrupper}.

\demo{Proof of Lemma {\rm \ref{Hupper2}}}  
We will show that for $0<\sg\le \log t$,
\begin{equation}\label{Shat}
S(t;\sg) \ll (1+\sg) \hat{S}(t;\sg),
\end{equation} 
where
\begin{equation}\label{Shatdef}
\hat{S}(t;\sg)= \!\!\!\sum_{a\in
\PP(e^\sg,te^\sg)} \frac{L(a;\sg)}{a\log^2(t^{7/8}/a+P^+(a))},
\end{equation} 
and then prove
\begin{equation}\label{hatstar}
\hat{S}(t;\sg) \ll S^*(t;\sg) \qquad (t\text{ large}, \sg>0).
\end{equation} 
First, assume $\sg \ge \log 2$, else \eqref{Shat} is trivial.
In \eqref{Sdef}, write each $a=a_1a_2$, where
$P^+(a_1) \le e^{\sg} < P^-(a_2)$.  
By Lemma \ref{Lineq} (iii) and (iv),
\begin{equation}\label{L34}
L(a;\sg) \le L(a_2;\sg+\log a_1) \le \(1+\frac{\log a_1}{\sg}\) L(a_2;\sg).
\end{equation} 
The contribution to $S(t;\sg)$ from those $a$ with $a_1\ge t^{1/8}$ is thus
\begin{equation}\label{a1a2}
\ll \frac{1}{\sg} \sum_{\substack{P^+(a_1) \le e^{\sg} \\ a_1 > t^{1/8}}}
\frac{\log a_1}{a_1\log^2 P^+(a_1)}
\sum_{a_2\in \PP(e^\sg,te^{\sg})}\frac{L(a_2;\sg)}{a_2}.
\end{equation} 
By Lemma \ref{Psixyz}, the sum on $a_1$ in \eqref{a1a2} is
\begin{align*}
&\le 4 \sum_{j\in\ZZ:t^{2^{-j}} \le e^{2\sg}} \frac{2^{2j}}
  {\log^2 t} \sum_{\substack{a_1>t^{1/8} \\ t^{2^{-j-1}} < P^+(a_1) \le
   t^{2^{-j}}}} \frac{\log a_1}{a_1} \\
&\ll \sum_{j\in\ZZ:2^j\ge \frac{\log t}{2\sg}}  2^j  
  e^{-2^{j-5}} \\
&\ll \frac{\sg^2}{\log^2 t} \left[ \pfrac{\log t}{\sg}^3 \exp \left\{ - 
  \frac{1}{64} \frac{\log t}{\sg} \right\} \right] \ll \frac{\sg^2}{\log^2 t}.
\end{align*}
It follows that the expression in \eqref{a1a2} is $\ll \sg \hat{S}(t;\sg)$.
By \eqref{L34}, $P^+(a) \ge P^+(a_2)$ and Lemma \ref{Psixyz},
the contribution to $S(t;\sg)$ from those $a$
with $a_1 \le t^{1/8}$ is at most
$$
\sum_{P^+(a_1)\le e^\sg} \frac{1+\tfrac{1}{\sg}\log a_1}{a_1} \!\!
\sum_{a_2\in\PP(e^\sg,te^\sg)} \frac{L(a_2;\sg)}{a_2\log^2(t^{7/8}/a_2+
P^+(a_2))} \ll \sg \hat{S}(t;\sg).
$$
This proves \eqref{Shat}.
Next, in \eqref{Shatdef} write $a=a_1a_2$ with
$$
a_1=\prod_{p^\b \| a, \b \ge 2} p^\b,
$$
so that $(a_1,a_2)=1$ and  $a_2$ is square-free.  By Lemma \ref{Lineq},
\begin{equation}\label{Sstar2}
L(a;\sg) \le \tau(a_1) L(a_2;\sg) \le \sg \tau(a_1) \tau(a_2).
\end{equation} 
Also
\begin{equation}\label{suma1}
\sum_{a_1 > w} \frac{\tau(a_1)}{a_1} \ll \sum_{a_1>w} a_1^{-3/4} \ll
w^{-1/4} \quad (w\ge 1).
\end{equation} 
The contribution to $\hat{S}(t;\sg)$
coming from those $a$ with $a_1 \ge t^{1/8}$ is thus
\begin{equation}\label{Sstar4}
\ll \sg \!\! \sum_{a_2\in\PP^*(e^\sg,te^{\sg})} \frac{\tau(a_2)}{a_2} 
\sum_{a_1 \ge t^{1/8}}
\frac{\tau(a_1)}{a_1} \ll \frac{\sg \log^2 t}{t^{1/32}} \ll S^*(t;\sg),
\end{equation} 
where we used \eqref{Sstarlower} in the last step.  Using
$P^+(a) \ge P^+(a_2)$ and \eqref{Sstar2}, we see that 
the contribution to $\hat{S}(t;\sg)$ from those $a$ with $a_1 < t^{1/8}$ is
$$
\le \sum_{a_1} \frac{\tau(a_1)}{a_1} \sum_{a_2\in \PP^*(e^{\sg},te^{\sg})}
\frac{L(a_2;\sg)}{a_2 \log^2 (t^{3/4}/a_2+P^+(a_2))} \ll S^*(t;\sg).
$$
This proves \eqref{hatstar}.
\Endproof\vskip4pt

%
%

Lemma \ref{Hrupper} depends on the distribution of the first
$r+1$ divisors of typical integers.  Throughout the remainder of this section,
let $d_j(n)$ denote the $j$-th smallest divisor of $n$.

\begin{lem}\label{djn}
Let $n$ have prime factorization $n=p_1^{e_1} \cdots p_f^{e_f}${\rm ,} where
$p_1 < \cdots < p_f$.  Let $N=\tau(n)${\rm ,} $1\le \ell \le N-1$ and define
$v$ uniquely by
\begin{equation}\label{e1}
(e_1+1) \cdots (e_{v-1}+1) | \ell, \qquad (e_1+1) \cdots (e_{v}+1) \nmid \ell.
\end{equation} 
Then $d_{\ell+1}(n)/d_{\ell}(n) \le p_v$.
\end{lem}

{\it Remarks}.  Lemma \ref{djn} is nearly best possible, e.g.\ if one takes
$e_1=\cdots=e_f=1$ and chooses the prime divisors of $n$ so that
$p_i > p_{i-1}^4$ for $2\le i\le f$, then $v=\nu(\ell)+1$ and
$$
\frac{d_{\ell+1}(n)}{d_{\ell}(n)} = 
\frac{p_v}{p_1 \cdots p_{v-1}} \ge p_v^{2/3}
$$
for every $\ell$.  
Moreover, this represents the typical case and is the origin of the exponent
$\nu(r)+1$ appearing in Lemma \ref{Hrupper} and Theorem \ref{Hrbounds}.

\Proof
We apply induction on $f$, the case $f=1$ being trivial.  Assume the statement
is true for $f=m$ and take an integer $n$ with prime factorization
$n=p_1^{e_1} \cdots p_{m+1}^{e_{m+1}}$, where
$p_1 < \cdots < p_{m+1}$.  Put $N=\tau(n)$, suppose $1\le \ell \le N-1$
and define $v$ by \eqref{e1}.  The conclusion is trivial if $v=m+1$; so
suppose $v\le m$.  Let $h=p_1^{e_1} \cdots p_m^{e_m}$ and
$E=e_{m+1}+1$.  For $0\le j\le E-1$, let 
$$
a_j = | \{ i\le \ell : p_{m+1}^j \| d_i(n) \}| = | \{ d|h: dp_{m+1}^j
\le d_{\ell}(n) \}|.
$$
Since $a_0+ \cdots + a_{E-1}=\ell$,   by
\eqref{e1}  there is a $j$ so that $(e_1+1) \cdots (e_v+1) \nmid a_j$
(in particular $1\le a_j\le \tau(h)-1$).  Define $u$ by
$$
(e_1+1) \cdots (e_{u-1}+1) | a_j, \qquad (e_1+1) \cdots (e_{u}+1) 
\nmid a_j,
$$
so that $1\le u \le v$. Finally, $d_{a_j+1}(h) p_{m+1}^j$ is a divisor
of $n$ that is larger than $d_{\ell}(n)$.  Therefore,
by the induction hypothesis,
\vskip12pt
\hfill $
\displaystyle{\frac{d_{\ell+1}(n)}{d_\ell(n)} \le \frac{d_{a_j+1}(h) p_{m+1}^j}{d_{a_j}(h)
p_{m+1}^j} \le p_{u} \le p_v.}
$ 
\Endproof

Using Lemma \ref{djn}, we can obtain upper bounds of averages of
ratios of consecutive divisors of numbers.

\begin{lem}\label{sumlogdl}
Fix $\ell \ge 1$.  Uniformly in $x\ge 1${\rm ,} $y\ge 4${\rm ,} 
$$
\sum_{\substack{P^+(n) \le y \\ n \ge x \\ \tau(n) \ge \ell+1}}
\frac{\log(d_{\ell+1}(n)/d_\ell(n))}{n} \ll_\ell (\log y) (\log\log y)
^{\nu(\ell)+1} \exp \biggl\{ - \frac{\log x}{4\log y} \biggr\}.
$$
\end{lem}

\Proof Denote by $S$ the sum in the lemma.  Let $V$ be the largest
number so that $\ell$ is the product of $V$ integers, each at least $2$ (if $\ell=1$, set $V=0$).
Suppose $n>x$, $\tau(n)\ge \ell+1$ and $P^+(n) \le y$.  
For each such $n$ there is
a unique $v$ between 1 and $V+1$, so that
$n=p_1^{e_1} \cdots p_v^{e_v} m$, where $p_1 < \cdots < p_v$, $m\in\PP(p_v,y)$
and
\begin{equation}\label{sumlogstar}
(e_1+1) \cdots (e_{v-1}+1) | \ell, \qquad (e_1+1) \cdots (e_{v}+1) \nmid \ell.
\end{equation} 
By Lemma \ref{djn}, $d_{\ell+1}(n)/d_{\ell}(n) \le p_v$.
Since $e_1+\cdots+e_{v-1} \le \ell$,
$m>x/y^{\ell+e_v}$.  Thus $S\le AB$, where (ignoring the second
 condition in \eqref{sumlogstar} and the condition $p_v>p_{v-1}$)
\begin{align*}
A &= \sum_{v=1}^{V+1} \sum_{e_1,\dots,e_{v-1}} \; 
\sum_{p_1<\cdots<p_{v-1}\le y} \frac{1}{p_1^{e_1} \cdots p_{v-1}^{e_{v-1}}}, \\
B &= \sum_{\substack{p\le y \\ f\ge 1}} \frac{\log p}{p^f} \sum_{\substack{
m\in \PP(p,y) \\ m>x/y^{\ell+f}}} \frac{1}{m}.
\end{align*}
For each tuple $v,e_1,\dots,e_{v-1}$ satisfying the first condition of 
\eqref{sumlogstar}, at most $\nu(\ell)$ of the numbers $e_i$ can equal 1.
Thus, by Lemma \ref{sum1p}, 
$$
A \ll_v (\log\log y)^{\nu(\ell)}.
$$
By Lemmas \ref{sum1p} and \ref{Psixyz},
\begin{align*}
B &\ll (\log y) \sum_{\substack{p\le y \\ f\ge 1}} \frac{1}{p^f} \exp
  \left\{ - \frac{\log(x/y^{\ell+f})}{4\log y} \right\} \\
&\le e^{\frac{\ell}{4} - \frac{\log x}{4\log y}} (\log y)
  \sum_{f=1}^\infty e^{f/4} \sum_{p\le y} p^{-f} \\
&\ll_{\ell} (\log y) (\log \log y)  e^{- \frac{\log x}{4\log y}}.
\end{align*}
\vglue-18pt
\Endproof

%
%

\demo{Proof of Lemma {\rm \ref{Hrupper}}} 
In light of Lemma \ref{Hupper}, it suffices to show that
\begin{equation}\label{Hrupmain}
S_s(t;\eta) \ll_r (\log \eta)^{\nu(r)+1} S^*(t;\eta)
\end{equation} 
for each $s$ satisfying $1\le s\le r$, $\nu(s) \le \nu(r)$ and each 
$t\ge y^{1/2}$.  Let $\th=\frac{\eta}{100Cr}$.  By hypothesis,
$\th \ge 1$.
In \eqref{Srdef}, write each $a=gh$, where
$P^+(g) \le e^{\th} < P^-(h)$.  Put $m=\tau(g)$.  We have
$S_s(t,\eta) = T_1+T_2+T_3$, where $T_1$ is the sum over numbers $a$
with $m\le s$ and $\nu(m)\le \nu(s)$; $T_2$ is the sum over those $a$
with $g\le t^{1/16}$ and with $m>s$ or $\nu(m) > \nu(s)$; and $T_3$ is the
sum over the remaining~$a$.

For $T_1$, we use 
$L_s(gh;\eta) \le L(gh;\eta) \le m L(h;\eta)$,
which is a consequence of Lemma \ref{Lineq}.  Also,
$$
g \le e^{m\th} \le e^{r\th} \le (z/y)^{1/100C} \le y^{1/100} \le t^{1/50};
$$
thus $t/(gh) + P^+(gh) \ge t^{15/16}/h+P^+(h)$.  This gives
\begin{equation}\label{HrT1}
T_1 \le r \sum_{\substack{P^+(g) \le e^{\th} \\ \tau(g)=m\le s, 
  \nu(m)\le \nu(s)}}
  \frac{1}{g} \sum_{h\in \PP(e^{\th},te^{\eta})} \frac{L(h;\eta)}
  {h \log^2 (t^{15/16}/h+P^+(h))}.
\end{equation} 
Factor each $g$ as $p_1^{e_1} \cdots p_j^{e_j}$.  Since
$m=(e_1+1)\cdots (e_j+1)$ and $j\le \frac{\log s}{\log 2}$, 
there are at most $\nu(m)$ indices $j$
with $e_j=1$.  Therefore, by Lemma \ref{sum1p}, the sum on $g$ in
\eqref{HrT1} is
$$
\le \sum_j \sum_{\substack{e_1,\cdots,e_j \\ p_1,\dots,p_j \le e^{\th}}}
\frac{1}{p_1^{e_1} \cdots p_j^{e_j}} \ll_m (\log 2\th)^{\nu(m)} \ll_r
(\log \eta)^{\nu(r)}.
$$
Consequently, 
\begin{equation}\label{T1bound}
T_1 \ll_r (\log \eta)^{\nu(r)} \sum_{h\in \PP(e^{\th},te^{\eta})} 
  \frac{L(h;\eta)}{h \log^2 (t^{15/16}/h+P^+(h))}.
\end{equation} 

For $g$ satisfying $m>s$ or $\nu(m)>\nu(s)$, we first prove that
\begin{equation}\label{HrL}
L_s(gh;\eta) \le 4r \frac{L(h;\eta)}{\th} \max_{\substack{1\le \ell \le
 \min(s,m-1) \\ \nu(\ell)\le \nu(s)}} \log\pfrac{d_{\ell+1}(g)}{d_\ell(g)}.
\end{equation} 
Write $D_j=d_j(gh)$ for $1\le j\le \tau(gh)$.
If $\tau(gh,e^u,e^{u+\eta})=s$, then $[u,u+\eta]$ contains
an interval $[\a,\b]$ so that $\tau(gh,e^\a,e^\b)=0$ and $\b-\a \ge 
\eta/(s+1) \ge 2\th$.  Let 
\begin{align*}
J&= \{ 1\le j\le \tau(gh)-1: \log(D_{j+1}/D_j) \ge 2\th \}, \\
M_j &= \{ u\in\RR: e^u < D_j \le e^{u+\eta}, \tau(gh,e^u,e^{u+\eta})=s \}\quad (1\le j\le
\tau(gh)).
\end{align*}
If $\tau(gh,e^u,e^{u+\eta})=s$, then either $u\in M_1 \cup M_{\tau(gh)}$ or
for some $j\in J$, $u\in M_j \cup M_{j+1}$.  Hence
\begin{equation}\label{LJM}
L_s(a;\eta) \le (2 |J|+2) \max_j \text{meas}(M_j).
\end{equation} 

Since $P^+(g) \le e^{\th}$, if $j\in J$,
then $g|D_j$ and $(g,D_{j+1})=1$, and thus $D_j/g$ and $D_{j+1}$ are 
consecutive divisors of $h$.  Also $D_{j+1}/(D_j/g) \ge e^{2\th}$,
hence the intervals $(\log D_{j+1}-\th,\log D_{j+1})$ are disjoint.
All such intervals lie in
$\LL(h;\th)$, thus $L(h;\th) \ge \th |J|$, and, since $L(h;\th)\ge \th$,
\begin{equation}\label{Jsize}
2|J|+2 \le 2\frac{L(h;\th)}{\th}+2 \le \frac{4L(h;\th)}{\th}.
\end{equation} 
Fix $j$ and let 
$$
\BB = \{ b\in \ZZ: \max(1,j+1-s) \le b \le \min(j,\tau(gh)+1-s) \},
$$
i.e., if $u\in M_j$, then the $s$ divisors
of $gh$ lying in $(e^u,e^{u+\eta}]$ are $D_b,\dots,D_{b+s-1}$ for some
$b\in \BB$.  Fix one such $b\in \BB$ and let
$$
M_{j,b} = \{ u\in M_j : e^u < D_b < \cdots < D_{b+s-1} \le e^{u+\eta} \}.
$$
For each $f|h$, let
$$
\KK_f = \{ j : d_j(g) f \in \{D_b,\dots,D_{b+s-1} \} \}, 
$$
which is a set of consecutive integers.  Also, since
$d_{j+1}(g)/d_j(g) \le e^{\th} < e^{\frac{\eta}{s+1}}$, each
set $\KK_f$ contains 1 or $m$.
Let $K_f = |\KK_f|$.  Since $\sum_f K_f = s$, there is some
$f$ so that $K_f\ge 1$ and $\nu(K_f) \le \nu(s)$.  
We have $K_f<m$, otherwise $s\ge m$ and $\nu(K_f)=\nu(m) > \nu(s)$.
If $\KK_f = \{1, 2,\dots, K_f\}$,
then $D_{b+s} \le d_{K_f+1}(g) f$, so that
$$
\text{meas}(M_{j,b}) \le \log \pfrac{d_{K_f+1}(g)}{d_{K_f}(g)}.
$$
Likewise, if $\KK_f = \{ m-K_f + 1,\dots, m\}$, then
$$
\text{meas}(M_{j,b}) \le  \log \pfrac{d_{m-K_f+1}(g)}{d_{m-K_f}(g)} =
\log \pfrac{d_{K_f+1}(g)}{d_{K_f}(g)},
$$
by the symmetry of the divisors of $g$.  Since $|\BB| \le s$,
$$
\text{meas}(M_j) \le s 
\max_{\substack{1\le \ell\le \min(s,m-1) \\ \nu(\ell)\le \nu(s)}}
 \log\pfrac{d_{\ell+1}(g)}{d_\ell(g)}.
$$
Combined with \eqref{LJM} and \eqref{Jsize}, this proves \eqref{HrL}.

If $g\le t^{1/16}$, then $t/gh+P^+(gh) \ge t^{15/16}/h+P^+(h)$.  Thus,
by \eqref{HrL} and Lemma \ref{sumlogdl},
\begin{equation}\label{T2bound}
\begin{split}
T_2 &\le \frac{4r}{\th} \max_{\substack{1\le \ell\le s \\ \nu(\ell)\le \nu(s)}}
  \sum_{\substack{P^+(g)\le e^{\th} \\ \tau(g)\ge \ell+1}}
  \frac{\log(\tfrac{d_{\ell+1}(g)}{d_\ell(g)})}{g}
  \sum_{h\in \PP(e^{\th},te^{\eta})} \frac{L(h;\eta)}
  {h\log^2(t^{15/16}/h+P^+(h))} \\
&\ll_r (\log \eta)^{\nu(s)+1}
  \sum_{h\in \PP(e^{\th},te^{\eta})} \frac{L(h;\eta)}{h\log^2(t^{15/16}/h+
  P^+(h))}.
\end{split}
\end{equation} 

If $g>t^{1/16}$, we use the bound $t/(gh) + P^+(gh) \ge P^+(g)$.  By Lemma
\ref{Lineq}, $L(h;\eta) \le \eta \tau(h)$.  Since
$$
\sum_{h\in \PP(e^{\th},te^{\eta})} \frac{\tau(h)}{h} \ll 
\frac{\log^2(te^\eta)}{\th^2} \ll_{r,C} \pfrac{\log t}{\th}^2,
$$
we have, by \eqref{HrL},
\begin{equation}\label{T3part1}
T_3 \ll_{r,C} \pfrac{\log t}{\th}^2 \max_{\substack{1\le \ell \le s \\ 
  \nu(\ell)\le \nu(s)}} \sum_{\substack{P^+(g) \le e^{\th} \\ \tau(g) \ge
 \ell+1 
  \\ g>t^{1/16}}} \frac{\log(d_{\ell+1}(g)/d_\ell(g))}{g\log^2 P^+(g)}.
\end{equation} 
By Lemma \ref{sumlogdl} and partial summation (with respect to $P^+(g)$),
the sum on $g$ in \eqref{T3part1} is
\begin{align*}
&\ll (\log \th)^{\nu(\ell)+1} \biggl[ \frac{1}{\th} e^{-\frac{\log t}{64\th}}
+\int_2^{e^\th} \frac{e^{-\frac{\log t}{64\log u}}}{u\log^2 u}\, du \biggr]\\
&\le (\log \th)^{\nu(\ell)+1} \( \frac{1}{\th} + \frac{64}{\log t} \) 
 e^{-\frac{\log t}{64\th}}.
\end{align*}
Since $\th \ll_{r,C} \log t$, \eqref{T3part1} and \eqref{Sstarlower} give
\begin{align*}
T_3 &\ll_{r,C} (\log \eta)^{\nu(s)+1} \frac{\th}{\log^2 t} \biggl[ 
\pfrac{\log t}{\th}^4  e^{-\frac{\log t}{64\th}} \biggr] \\
&\ll_{r,C} (\log \eta)^{\nu(r)+1} \frac{\th}{\log^2 t} 
\ll_{r,C}  (\log \eta)^{\nu(r)+1}
S^*(t;\eta).
\end{align*}
Combining this with \eqref{T1bound} and \eqref{T2bound}, we obtain
$$
S_s(t,\eta) \ll_{r,C} (\log \eta)^{\nu(r)+1} \( S^*(t,\eta) + 
\sum_{h\in \PP(e^{\th},te^{\eta})} \frac{L(h;\eta)}{h\log^2 (t^{15/16}/h+
P^+(h))} \).
$$

As in the proof of \eqref{Shat}, write each $h$ in the above sum
as $h=a_1a_2$, where $a_1\in \PP(e^{\th},e^{\eta})$ and
$a_2\in\PP(e^{\eta},te^{\eta})$.   Let $T_4$ be the contribution to the sum
from those $h$ with $a_1 \le t^{1/16}$, and let $T_5$
be the remaining portion of the sum.  By Lemma \ref{Lineq}, 
$L(h;\eta) \le \tau(a_1) L(a_2;\eta)$, and we also have
$$
\sum_{a_1\in \PP(e^\th,e^{\eta})} \frac{\tau(a_1)}{a_1} \ll \pfrac{\eta}
{\th}^2 \ll_{r,C} 1.
$$
Thus $T_4 \ll_{r,C} \hat{S}(t;\eta)$.  When $a_1>t^{1/16}$, Lemma \ref{Lineq}
gives
$$
L(h;\eta) \le \(1+ \frac{\log a_1}{\eta} \) L(a_2;\eta) \ll_{r,C} 
\frac{\log a_1}{\eta} L(a_2;\eta).
$$
Also, $t^{15/16}/h+P^+(h) \ge P^+(a_1) \ge e^\th$, so
$$
T_5 \ll_{r,C} \frac1{\eta^3} \sum_{\substack{a_1\in \PP(e^{\th},e^{\eta}) \\
a_1 > t^{1/16}}} \frac{\log a_1}
{a_1} \sum_{a_2\in\PP(e^\eta,te^{\eta})} \frac{L(a_2;\eta)}{a_2}.
$$
By Lemma \ref{Psixyz}, the sum on $a_1$ is
$$
\ll_{r,C} \frac{\eta^3}{\log^2 t} \biggl[ \pfrac{\log t}{\eta}^3 
e^{-\frac{\log t}{64
\eta}} \biggr] \ll_{r,C} \frac{\eta^3}{\log^2 t}.
$$
Thus, $T_5 \ll_{r,C} \hat{S}(t;\eta)$ and hence
$$
S_s(t;\eta) \ll_{r,C} (\log \eta)^{\nu(r)+1} \bigl( S^*(t;\eta)+\hat{S}(t;\eta)
\bigr).
$$
The lemma now follows from \eqref{hatstar}.
\hfill\qed

%
%
\section{Upper bounds: reduction to an integral}\label{sec:uppertoint}
%
%

In this section we prove Lemma \ref{Tnint}.
We have $T_0(\sigma,P,Q)=0$ if $Q>1$ and
$T_0(\sigma,P,1)=\sigma$.  Now let $k\ge 1$ and
put $\g = \frac{1}{\log P}$.  We have
$$
T_k(\sg,P,Q) \le Q^{-\g} \sum_{\substack{a \in \PP^*(e^{\sg},P) \\
\omega(a)=k}} \frac{L(a;\sg)}{a^{1-\g}}.
$$
The parameter $\g$ has been chosen so that the sum is only 
a constant multiplicative factor larger than the
corresponding sum with $\g=0$.  Since $p^\g \le 1 + 2\g \log p$
for $p\le P$, we have
$$
\sum_{e^{\sg} <p \le P} \frac{1}{p^{1-\g}} \le \sum_{e^{\sg} <p \le P}
\frac{1}{p} + 2\g \sum_{p\le P} \frac{\log p}{p} \le \log\log P -
\log (1+\sg) + O(1).
$$
By an argument similar to that used to construct the sets $D_j$ in Section 
\ref{sec:lower}, we find that there is an absolute constant $K$ so that the
following holds for all $\sg,P$:  the 
interval $(e^{\sg},P]$ may be partitioned into subintervals
$E_0,\dots,E_{v+K-1}$ with $v$ as given in Lemma \ref{Tnint} and
for each $j$,
$$
\sum_{p\in E_{j}} \frac{1}{p^{1-\g}} \le \log 2
$$
and
\begin{equation}\label{81}
p\in E_j \, \implies \, \frac{\log\log p - \log (1+\sg)}{\log 2} \le j+K.
\end{equation} 

Consider $a=p_1 \cdots p_k$, $e^{\sg} < p_1 < \cdots < p_k \le P$ 
and define $j_i$ by 
$p_i\in E_{j_i}$  ($1\le i\le k$).  Put $l_i = \frac{\log\log p_i}{\log 2}$.
By Lemma \ref{Lineq} (v) and \eqref{81},
\begin{align*}
L(a;\sg) &\le 2^k \min_{0\le g\le k} 2^{-g} (2^{l_1}+\cdots +2^{l_g}+\sg) \\*
&\le (\sg+1) 2^{k+K} F(\bj),
\end{align*}
where
$$
F(\bj) = \min_{0\le g\le k} 2^{-g} (2^{j_1}+\cdots+ 2^{j_g}+\min(1,\sg)).
$$
Let $J$ denote the set of vectors $\bj$ satisfying
$0\le j_1 \le \cdots \le j_k \le v+K-1.$  Then
$$
T_k(\sg,P,Q) \le Q^{-\g} 2^{k+K} \sum_{\bj \in J} F(\bj) 
\sum_{\substack{p_1< \cdots < p_k \\ p_i\in E_{j_i} \; (1\le i\le k)}} \frac{1}
{(p_1\cdots p_k)^{1-\g}}.
$$
If $b_j$ is the number of primes $p_i$ in $E_{j}$ for $0\le j\le v+K-1$,
the sum over $p_1,\cdots,p_k$ above is at most
\begin{align*}
\prod_{j=0}^{v+K-1} \frac{1}{b_j!} \biggl( \sum_{p\in E_j}
\frac{1}{p^{1-\g}} \biggr)
^{b_j} &\le \frac{(\log 2)^k}{b_0! \cdots b_{v+K-1}!} \\
&= ((v+K)\log 2)^k
\int_{R(\bj)} 1\, d\bx \\
&\le e^{10K} (v\log 2)^k \int_{R(\bj)} 1\, d\bx,
\end{align*}
where
$$
R(\bj) = \{ 0\le \xi_1\le \cdots \le \xi_k\le 1: 
j_i \le (v+K) \xi_i \le j_i+1 \,\; \forall i\} \subseteq R_k.
$$
Finally, since $2^{j_i} \le 2^{(v+K)\xi_i} \le 2^K 2^{v\xi_i}$ for each $i$,
$$
\sum_{\bj \in J} F(\bj) \int_{R(\bj)} 1 d\bx \le 2^K U_k(v;\a). 
$$

%
%
\section{Lower bounds: isolated divisors}\label{sec:isolated}
%
%

In this section we prove Lemmas \ref{Hrlower1} and \ref{Hrlower2}.

\demo{Proof of Lemma {\rm \ref{Hrlower1}}}
Write $a=a' p_1 \cdots p_r$, where
$$
a' \le y^{\frac{c'}{100r}}, \qquad y^{\frac{c'}{2r}} < p_1 < \cdots <
p_r \le y^{\frac{c'}{r}}.
$$
We have
\begin{equation}\label{Hrlow1}
\sum_{a\le y^{2c'}} \frac{L_r(a;\eta)}{a} \ge \sum_{a'} 
\frac{1}{a'} \int \sum_{\substack{p_1,\dots,p_r \\ \tau(a'p_1\cdots p_r,
e^u,e^{u+\eta})=r}} \frac{1}{p_1\cdots p_r} \; du.
\end{equation} 
Fix $a'$ and fix $u$ so that $y^{\frac{0.6 c'}{r}} \le e^u \le 
y^{\frac{0.7c'}{r}}$.  The measure of such $u$ is $\gg_{r,c'} \log y$.
Since $e^{u+\eta} \le y^{\frac{0.8c'}{r}}$, a divisor of $a'p_1\cdots p_r$ 
lying in
$(e^u,e^{u+\eta})$ must have the form $dp_i$, where $d|a'$ and $1\le i\le r$.
If 
\begin{equation}\label{Hrlow2}
\log\pfrac{e^u}{p_i} \in \LL_1(a';\eta) \qquad (1\le i\le r),
\end{equation} 
there are exactly $r$ such divisors.  Assume that $I(a';\eta) \ge 1$.
If we take, for each $i$,
$e^u/d_i < p_i \le e^{u+\eta}/d_i$ for some $\eta$-isolated divisor $d_i$ of
$a'$, then \eqref{Hrlow2} will be satisfied.  By Lemma \ref{sum1p},
the sum over $p_1,\dots,p_r$ in \eqref{Hrlow1} is
$$
\ge \frac{1}{r!} \prod_{i=1}^r \biggl( 
\sum_{\log(e^u/p_i)\in \LL_1(a';\eta)} \frac{1}{p_i} \biggr) -
O_{r,c'}\pfrac{1}{y^{c'/2r}} \gg_{r,c'} 
\pfrac{\eta I(a';\eta)}{\log y}^r.
$$
By inserting this into \eqref{Hrlow1} and using Lemma \ref{basiclower}, we
complete  the proof.
\Endproof\vskip4pt

\emph{Proof of Lemma} \ref{Hrlower2}.
(i) Write $a=gh$, where $g \le K := \min(e^{\eta},y^{c'})$, 
$e^{2\eta} < P^-(h)$, $h\le y^{c'}$ and
$\tau(g) \ge r+1$.  In particular, if $\eta \ge \frac{c'}{2}\log y$,
then $h=1$.
Let $f$ be a $2\eta$-isolated divisor of $h$.  
As before, let $d_j(g)$ be the $j$th smallest divisor of $g$.
If $fd_r(g) \le e^{u+\eta} \le fd_{r+1}(g)$, then 
$\tau(gh,e^u,e^{u+\eta})=r$ and hence
$$
L_r(gh;\eta) \ge \log\pfrac{d_{r+1}(g)}{d_r(g)} I(h;2\eta).
$$
By Lemma \ref{basiclower},
\begin{equation}\label{Hrlow3}
H_r(x,y,z) \gg_{r,c',C} 
\frac{x}{\log^2 y} \, A_r(K) \, \sum_{\substack{h\le y^{c'}
\\ P^-(h) > e^{2\eta}}} \frac{I(h;2\eta)}{h},
\end{equation} 
where
$$
A_r(K) = \sum_{\substack{g\le K \\ \tau(g) \ge r+1}} \frac{\log(d_{r+1}(g)/
d_r(g))}{g}.
$$
To complete the proof of part (i) of the lemma, it suffices to show
\begin{equation}\label{Aineq}
A_r(K) \gg_{r} (\log K) (\log \log K)^{\nu(r)+1} \qquad (K \ge
e^{1000\cdot 3^{2r}})
\end{equation} 
and note that $K \ge e^{\eta c'/C}$.  
We restrict attention to those $g=p_1 \cdots p_r m$, with
\begin{equation}\label{Hrlowp}
p_1 \le K^{1/3^{2r}}, \;\;\;
p_{j-1}^3 < p_j < K^{1/3^{2r-2j+2}} \quad (2\le j\le r),
\end{equation} 
$P^-(m) > p_r$ and $m \le K^{1/2}$.
Inequality \eqref{Hrlowp} implies that 
$$
\frac{d_{r+1}(g)}{d_r(g)} = \frac{p_v}{p_1 \cdots p_{v-1}} \ge p_v^{1/2}, 
\quad v=\nu(r)+1.
$$
Hence
$$
A \ge \frac12 \sum_{\substack{p_1,\dots,p_r \\ (9.5) }} \frac{\log p_v}{p_1\cdots p_r} 
\sum_{\substack{P^-(m) > p_r \\ m \le K^{1/2}}} \frac{1}{m}.
$$
Since $p_r \le K^{1/9}$,
by Lemma \ref{Phi} and partial summation, the sum on $m$ is
$\gg \frac{\log K}{\log p_r}.$
If $v=r$, then $r\in \{1,2\}$ and we obtain
$$
A \gg \log K \sum_{p_1,\cdots,p_r} \frac{1}{p_1\cdots p_r} \gg
\log K (\log\log K)^{r},
$$
proving \eqref{Aineq} in this case.  Otherwise, $v<r$.
By Lemma \ref{sum1p} and Bertrand's Postulate, if $u>w^2>4$, then 
$$
\sum_{w < p \le u} \frac{1}{p\log p} \ge
\frac{1}{2\log w} \sum_{w<p\le w^2} \frac{1}{p} \gg \frac{1}{\log w}.
$$
Applying this iteratively and using \eqref{Hrlowp}, we obtain
$$
\sum_{p_{v+1},\dots, p_r} \frac{1}{p_{v+1} \cdots p_r \log p_r} \gg_r
\frac{1}{\log p_v}.
$$
Finally, we have
$$
\sum_{p_1,\cdots,p_v} \frac{1}{p_1\cdots p_v} \gg_r
(\log \log K)^{v} = (\log \log K)^{\nu(r)+1}.
$$
This is trivially true when $K \le K_0(r)$, for a large constant
$K_0(r)$, and for larger $K$ the sum includes $p_1,\dots,p_v$ satisfying
$$
p_i \in \( \exp[(\log K)^{i/2r}], \exp[(\log K)^{(i+1/2)/2r}] \) (1\le
i\le v)
$$
by \eqref{Hrlowp}.  This proves \eqref{Aineq} in the second case.

(ii) Suppose $y_0(r,c) \le y$ and $y^2 \le z \le x^{1-c}/y$.
Consider $n=gpq$, where $g\le y^{c/2}$, $\tau(g)
\ge r+1$, $p$ is prime, $yd_r(g)/g < p \le yd_{r+1}(g)/g$,
$P^-(q)>z$ and $\frac{x}{2pg} < q \le \frac{x}{pg}$.
We have $y\ge p > y/g > g$, so each $n$ has at most one factorization
of this type.   If $d|n$ and $y<d\le z$, then $p|d$.  Thus
$\tau(n,y,z)=r$, because 
$$
\frac{pg}{d_{r+1}(g)} \le y < \frac{pg}{d_r(g)} <
\frac{pg}{d_{r-1}(g)} < \cdots < pg \le z.
$$
Since $pg \le yx^{c/2} \le \frac{x}{4z}$, Lemma \ref{Phi}
implies that for each pair $g,p$, the number of $q$ is $\gg
\frac{x}{pg\log z}$.  By Lemma \ref{sum1p}, for each $g$,
$$
\sum \frac{1}{p} \gg \frac{\log(d_{r+1}(g)/d_r(g))}{\log y}.
$$
Using \eqref{Aineq}, the number of such $n$ is
\vskip12pt \hfill
$ 
\displaystyle{\gg \frac{x}{(\log y)(\log z)} A_r(y^{c/2}) \gg_{r,c} \frac{x
(\log\log y)^{\nu(r)+1}}{\log z}.}
$ 
\Endproof

\ifdraft
\vfil\eject
\else
\fi

%
%
\section{Lower bounds: reduction to a volume}\label{sec:lowtovol}
%
%

\emph{Proof of Lemma} \ref{sumW}.
Recall the definitions of the sets $D_j$ and numbers $\lam_j$
from Section \ref{sec:lower}.
For $j\ge 0$ let $b_j'=\sum_{i\le j} b_j$.  Let $a=p_1\cdots p_k$, where
 $k=b_1+\cdots+b_h$,
\begin{equation}\label{sumW_a}
p_{b_{j-1}'+1},\dots, p_{b_j'} \in D_j \qquad (m\le j\le h)
\end{equation} 
and the primes in each interval $D_j$ are unordered. 
Observe that $W(p_1\cdots p_k;\sg)$
is the number of pairs $Y,Z\subseteq \{1,\dots,k\}$ with
\begin{equation}\label{sumW_b}
\left| \sum_{i\in Y} \log p_i - \sum_{i\in Z} \log p_i \right| \le \sg.
\end{equation} 
We thus have
\begin{equation}\label{sumW_c}
\sum_{a\in \AA(\bb)} \frac{W(a;\sg)}{a} \le \frac{1}{b_m! \cdots b_h!}
\sum_{Y,Z \subseteq \{1,\dots,k\}} \;\; \sum_{\substack{p_1,\dots,p_k \\
\eqref{sumW_a}, \eqref{sumW_b}}} \frac{1}{p_1 \cdots p_k}.
\end{equation} 
The contribution from terms with $Y=Z$ is
\begin{equation}\label{sumW_trivial}
\frac{2^k}{b_m!\cdots b_h!} \sum_{\substack{p_1,\dots,p_k \\ \eqref{sumW_a}}}
\frac{1}{p_1 \cdots p_k} 
\le \frac{(2\log 2)^k}{b_m!\cdots b_h!}.
\end{equation} 

If $Y\ne Z$, let $I=\max [(Y\cup Z)-(Y\cap Z)]$ and define $E(I)$ by 
$p_I\in D_{E(I)}$, i.e., $b'_{E(I)-1} <  I \le b'_{E(I)}$.
Let
$$
\ell = \min \{ j : \lam_j \ge \sg^{-2} \}.
$$
We distinguish two cases: (i) $E(I) > \ell$; (ii) $m \le E(I) \le \ell$.

In case (i), fix all of the $p_i$ except for $p_I$. 
Inequality \eqref{sumW_b} implies that
$U \le p_I \le e^{2\sg} U$ for some number $U\ge \lam_{E(I)}$.  
If $\sg\ge 1$, then by \eqref{sum1p},
$$
\sum_{U\le p_I \le e^{2\sg} U} \frac{1}{p_I} \ll \log\(1 + \frac{2\sg}{\log
  U}\) \ll \frac{\sg}{\log U} \ll \sg 2^{-E(I)}.
$$
If $\sg<1$, Lemma \ref{BrunTitch} implies
$$
\sum_{U\le p_I \le e^{2\sg} U} \frac{1}{p_I} \ll \frac{\sg}{\log (\sg U)}
\ll \frac{\sg}{\log U} \ll \sg 2^{-E(I)},
$$
where the second inequality follows from $U\ge \lam_\ell \ge \sg^{-2}$.
Thus, by \eqref{Dj} 
the inner sum on the right side of
\eqref{sumW_c} is $\ll \sg 2^{-E(I)} (\log 2)^k$.  With $I$
fixed, there correspond $2^{k-I+1} 4^{I-1}=2^{k+I-1}$ pairs $Y,Z$.
We find that the contribution to the right side of \eqref{sumW_c}
from those $Y,Z$ counted in case (i) is
\begin{equation}\label{sumW_casei}
\begin{split}
&\ll \frac{\sg (2\log 2)^k}{b_m! \cdots b_h!} \sum_{I=1}^k 2^{I-E(I)} \\
&\ll \frac{\sg (2\log 2)^k}{b_m! \cdots b_h!} \sum_{j=m}^h 2^{-j} 
\sum_{b'_{j-1}<I\le b'_j} 2^I
\ll \frac{\sg (2\log 2)^k}{b_m! \cdots b_h!} \sum_{j=m}^h 
2^{-j+b_m+\cdots+b_j}.
\end{split}
\end{equation} 

For case (ii), note that $\sg<1$ and write
$$
a = a' p_{b'_{\ell}+1} \cdots p_k, \qquad a'=p_1\cdots p_{b'_{\ell}}.
$$
By hypothesis, $Y \cap \{b'_\ell+1,\dots,k\} = 
Z \cap \{b'_\ell+1,\dots,k\}$.  By \eqref{Dj}, the contribution to
the right side of \eqref{sumW_c} from those $Y,Z$ counted in case (ii)
is
$$
\le \frac{(2\log 2)^{b_{\ell+1}+\cdots+b_h}}{b_{\ell+1}! \ldots b_h!}
\sum_{a'} \frac{W(a';\sg)-\tau(a')}{a'}.
$$
Suppose $d_1|a'$, $d_2|a'$ and $1 < d_2/d_1 \le e^{\sg}$.  Let
$d=(d_1,d_2)$, $d_1=f_1 d$, $d_2=f_2 d$ and $a'=df_1f_2 a''$.
Then
\begin{align*}
\sum_{a'} \frac{W(a';\sg)-\tau(a')}{a'} &\le
2 \sum_{P^+(a''d f_1) \le \lam_{\ell}} \frac{1}{a'' d f_1} \sum_{f_1 < f_2 \le
    e^{\sg} f_1} \frac{1}{f_2} \\
&\le 4\sg  \sum_{P^+(a''d f_1) \le \lam_{\ell}} \frac{1}{a'' d f_1} \\
&= 4\sg \prod_{p\le \lam_{\ell}} \( 1 + \frac{1}{p} \)^3 \\
&\le 4\sg \exp \( 3 \sum_{p\le \lam_{\ell}} \frac{1}{p} \) 
  \le 2^{3\ell+2} \sg.
\end{align*}
Since $\sg > \lam_{\ell-1}^{-1/2}$, we have 
$2^{3\ell+2} \sg < \exp \{ - 2^{2\ell/3} \}$
if $m$ is large.  
On the other hand, we have assumed that $b_j \le 2^{j/2}$ for every $j$,
and thus
$$
\frac{1}{b_m! \cdots b_\ell!} \ge \( 2^{(\ell/2) \cdot 2^{\ell/2} } \)^{-\ell}
\ge 100 \exp \{ -2^{2\ell/3} \} > 100 \cdot 2^{3\ell+2}\sg
$$
for large $m$.  Therefore, the contribution to the right side of
\eqref{sumW_c} from those $Y,Z$ counted in case (ii) is
\begin{equation}\label{sumW_caseii}
\le 0.01 \frac{(2\log 2)^{k}}{b_{m}! \cdots b_h!}.
\end{equation} 

Together, \eqref{sumW_c},
\eqref{sumW_trivial}, \eqref{sumW_casei}, and \eqref{sumW_caseii} imply that
\vskip12pt \hfill
$ 
\displaystyle{\sum_{a\in \AA(\bb)} \frac{W(a;\sg)}{a} \le \frac{(2\log 2)^{k}}{b_{m}! 
\ldots b_h!} \left[ 1.01 + O\( \sg \sum_{j=m}^h
2^{-j+b_m+\cdots+b_j}\)\right]. }
$ 
\Endproof\vskip12pt

%
%

\emph{Proof of Lemma} \ref{lowvol}.
Let $m = M + \max\( 0, \fl{\frac{\log \sg}{\log 2}} \),$
put $h=v+m-1$ and suppose $\bb$ satisfies $b_j=0$ ($1\le j\le m-1$),
$b_1+\cdots+b_h=k$ and
\begin{equation}\label{bcond}
\text{for all } j\ge 1, \; b_j\le 2^{j/10} \text{ and } b_{h+1-j} \le
2^{(M+j)/10}.
\end{equation} 
Since \pagebreak $h\le \frac{\log\log y}{\log 2}-M-1/\a$, 
\begin{align*}
\sum_j b_j 2^j &\le 2^{h+\frac{M+1}{10}} (1+2^{-0.9}+2^{-1.8}+\cdots) \\
&\le 2^{h+M/10+2} \\
&\le\frac{\log y}{2^{\frac{9}{10}M-2+1/\a}}
\le \frac{\a \log y}{2^{c_3+c_4}}.
\end{align*}
Using Lemma \ref{muj}, it follows that $a\le y^{\a}$ for $a\in\AA(\bb)$.
Also $P^-(a) > \lam_{m-1}\break > e^{\sg}$.  By the definition of the sets $D_j$,
\begin{equation}\label{Rlower}
\begin{split}
\sum_{a\in\AA(\bb)} \frac{1}{a} &= \prod_{j=m}^h \frac{1}{b_j!} \biggl(
  \sum_{p_1\in D_j} \frac{1}{p_1} \sum_{\substack{p_2\in D_j \\ p_2 \ne p_1}}
  \frac{1}{p_2} \cdots \sum_{\substack{p_{b_j}\in D_j \\ p_{b_j} \not\in
  \{ p_1,\dots, p_{b_j-1} \} }} \frac{1}{p_{b_j}} \biggr) \\
&\ge\prod_{j=m}^h \frac{1}{b_j!} 
   \biggl( \sum_{p\in D_j} \frac{1}{p} - \frac{b_j-1}{\lam_{j-1}} 
   \biggr)^{b_j} \\
&\ge \prod_{j=m}^h \frac{1}{b_j!}  \biggl( \log 2 - \frac{b_j}{\lam_{j-1}}
   \biggr)^{b_j} \\
&\ge \frac{(\log 2)^k}{b_m! \cdots b_h!} \prod_{j=m}^h \( 1 -
  \frac{2^{j/10}}{\exp\{2^{j-1+c_3-c_4} \}} \)^{2^{j/10}} \\
&\ge 0.999 \frac{(\log 2)^k}{b_m! \cdots b_h!}.
\end{split}
\end{equation} 
By Lemma \ref{sumW} and \eqref{Rlower},
$$
\sum_{a\in\AA(\bb)} \frac{3\tau(a)-2W(a;\sg)}{a} \ge 
  \frac{(2\log 2)^k}{b_m! \cdots b_h!}
\biggl[ 0.9 - 2^{c_5+1} \sg \sum_{j=m}^h 2^{-j+b_m+\cdots+b_j} 
  \biggr].
$$
For $1\le i\le v$, set $g_i=b_{m-1+i}$.  We have
$2^{(M+i-1)/10} \ge M+i^2$ for $M\ge 200$ and $i\ge 1$, 
so \eqref{bcond} is implied by
\begin{equation}\label{gcond}
\forall i\ge 1, g_i \le M+i^2 \text{ and } g_{v+1-i} \le M+i^2.
\end{equation} 
Assume in addition that
\begin{equation}\label{gcond2}
2^{m-1}\sum_{j=m}^h 2^{-j+b_m+\cdots+b_j} = 
\sum_{i=1}^v 2^{-i+g_1+\cdots+g_i} \le 2^{s+1}.
\end{equation} 
Since $2^{s-m+2}=\frac{1}{\sg} 2^{-c_5-8}$, we obtain
\begin{equation}\label{L1a}
\sum_{a\in\AA(\bb)} \frac{3\tau(a)-2W(a;\sg)}{a} \ge  \frac{(2\log 2)^k}
{3 g_1! \cdots g_v!}.
\end{equation} 
Let $\GG$ be the set of $\bg=(g_1,\cdots,g_v)$ with $g_1+\cdots+g_v=k$
and also satisfying \eqref{gcond} and \eqref{gcond2}.
We have
\begin{equation}\label{L1b}
\frac{1}{g_1! \cdots g_v!} = \Vol(R(\bg)),
\end{equation} 
where
$R(\bg)$ is the set of $\xx\in\RR^k$ with
$0\le x_1\le\cdots\le x_k<v$ and exactly $g_i$  of the 
variables $x_j$ lie in $[i-1,i)$ for each $i$.  We claim that
\begin{equation}\label{L1claim}
 \Vol(\cup_{\bg\in\GG} R(\bg))  \ge v^k \Vol(Y_k(s,v)).
\end{equation} 
Take $\bx \in Y_k(s,v)$ with $\xi_k<1$ 
and let $x_j=v\xi_j$ for each $j$.  Let $g_i$ be the
 number of $x_j$ lying in $[i-1,i)$.
By condition (ii) in the definition of $Y_k(s,v)$,
$$
x_{M+i^2} > i\; \text{ and }\; x_{k+1-(M+i^2)} < v-i \qquad (1\le i\le
\sqrt{k-M}),
$$
which implies \eqref{gcond}.
Condition (iii) in the definition of $Y_k(s,v)$ implies
$$
2^s \ge \sum_{j=1}^k 2^{j-x_j} \ge
\sum_{i=1}^v 2^{-i} \sum_{j:x_j\in [i-1,i)} 2^j
\ge \sum_{i=1}^v 2^{-i + g_1 + \cdots + g_i-1}.
$$
Because \eqref{gcond} and \eqref{gcond2} hold,
$\bg \in \GG$ and $\xx\in R(\bg)$.  This proves the
claim \eqref{L1claim}.  Finally, combining \eqref{L1a}, \eqref{L1b}
and \eqref{L1claim} proves the lemma.
\hfill\qed

\ifdraft
\vfil\eject
\else
\fi

%
%
%
\section{Uniform order statistics}\label{sec:uos}
%
%
%

Let $X_1,\dots, X_k$ be independent, uniformly distributed random
variables in $[0,1]$, and let $\xi_1,\dots,\xi_k$ be their order
statistics, so that
$0 \le \xi_1 \le \cdots \le \xi_k \le 1.$
Our main interest in this section is to estimate $Q_k(u,v)$,
the probability that $\xi_i \ge \frac{i-u}{v}$ for every $i$.
The special case $Q_k(\lam \sqrt{k},k)$ is the distribution function (as a
function of $\lambda$) of the quantity
$$
D_k^+ = - \sqrt{k} \inf_{1\le i\le k} (\xi_i-i/k) = \sqrt{k} 
\sup_{0\le y\le 1} \( k^{-1} \sum_{X_i\le y} 1 - y \),
$$
which is one of the Kolmogorov-Smirnov statistics.
In 1939, N. V. Smirnov proved \cite{Sm1} for each fixed
$\lam\ge 0$ the asymptotic formula
$$
Q_k(\lam \sqrt{k},k) \sim 1 - e^{-2\lam^2} \qquad (k\to\infty).
$$
We need bounds on $Q_k(u,v)$ which are uniform in $k,u$ and $v$,
particularly when $u$ is small or $u+v-n$ is small.

\begin{lem}\label{Q2}  For a given pair $u,v${\rm ,} let $w = w(u,v) = u+v-n$.
Uniformly in $u>0${\rm ,} $w>0$ and $k\ge 1${\rm ,} 
\begin{equation}\label{Qupper}
Q_k(u,v) \ll \frac{(u+1)(w+1)^2}{k}.
\end{equation} 
If $k\ge 1$ and $1\le u\le k${\rm ,} then
\begin{equation}\label{Qlower}
Q_k(u,k+1-u) \ge \frac{u-1/2}{k+1/2}.
\end{equation} 
\end{lem}

{\it Remarks}.  The more precise estimate
$$
Q_k(u,v) = 1 - e^{-2uw/k} + O\pfrac{u+w}{k} \qquad (k\ge 1, u\ge 0, w\ge 0),
$$
which was proved in an earlier version of this paper, has a much longer proof
which will be given in a separate paper \cite{uos}.

\Proof
If $\min_{1\le i\le k} (\xi_i - \frac{i-u}{v}) \le 0$, let $l$ be 
the smallest index with $\xi_l \le \frac{l-u}{v}$ and write 
$\xi_l=\frac{l-u-\lam}{v}$, so that $0\le \lam \le 1$.  Also let
$$
R_l(\lam) = \Vol\left\{ 0\le \xi_1 \le \cdots \le \xi_{l-1} \le
\frac{l-u-\lam}{v} : \xi_i \ge \frac{i-u}{v} \, (1\le i\le l-1) \right\}.
$$
Then
\begin{align*}
Q_k(u,v) &= 1 - \frac{k!}{v} \int_0^1 \sum_{u+\lam \le l \le k}
R_l(\lam) \Vol \left\{ \frac{l-u-\lam}{v} \le \xi_{l+1} \le \cdots \le
\xi_k \le 1 \right\} \, d\lam \\
&= 1 - \frac{k!}{v} \int_0^1 \sum_{u+\lam \le l \le k}
\frac{R_l(\lam)}{(k-l)!} \pfrac{k+w+\lam-l}{v}^{k-l}\, d\lam.
\end{align*}

Let $0\le a\le k-u$ and suppose that 
$\xi_k \le 1 - \frac{w+a}{v} = \frac{k-u-a}{v}$.
Then $\min_{1\le i\le k} \xi_i-\frac{i-u}{v} \le 0$.  Defining $l$ and $\lam$
as before, we have
\begin{align*}
\(1 - \frac{w+a}{v}\)^k &= k! \Vol \biggl\{ 
0\le \xi_1\le \cdots \le \xi_k \le 1 - \frac{w+a}{v} \biggr\} \\
&= \frac{k!}{v} \int_0^1 \sum_{u+\lam \le l \le k-a+\lam} 
\frac{R_l(\lam)}{(k-l)!} \pfrac{k-l-a+\lam}{v}^{k-l}\, d\lam.
\end{align*}
Thus, for any $A>0$, 
\begin{multline}\label{QkA}
Q_k(u,v) = 1 - A \( 1 - \frac{w+a}{v} \)^k \\ \phantom{Q_k(u)=}
- \frac{k!}{v} \int_0^1
\sum_{k-a+\lam <l \le k} \frac{R_{l}(\lam)}{(k-l)!} 
\pfrac{k+w+\lam-l}{v}^{k-l}\, d\lam \\
+ \frac{k!}{v} \int_0^1 \! \sum_{l=\fl{u+\lam+1}}^{\fl{k-a+\lam}}
  \frac{R_{l}(\lam)}{(k-l)! v^{k-l}} \left[
    A(k-l-a+\lam)^{k-l}\!-\!(k-l+w+\lam)^{k-l} \right] d\lam. 
\end{multline}

To prove \eqref{Qupper}, we may assume without loss of generality that $k\ge 10$,
$u\le k/10$ and $w\le \sqrt{k}$.  Let $a=2w+2$, and 
note that $2-\lam \ge \lam$. \pagebreak Now,
\begin{align*}
\biggl( \frac{k-l-w-2+\lam}{k-l+w+\lam} & \biggr)^{k-l} = 
\(1 - \frac{w+2-\lam}{k-l}\)^{k-l}
\(1 + \frac{w+\lam}{k-l}\)^{-(k-l)} \\
&\nts\nts\nts\nts= \exp \left\{ -(2w+2) +
\sum_{j=2}^\infty \frac{-(w+2-\lam)^j + (-1)^j(w+\lam)^j}{j(k-l)^{j-1}}
\right\} \\
&\nts\nts\nts\nts\le e^{-(2w+2)}.
\end{align*}
Thus, taking $A=e^{2w+2}$ in \eqref{QkA}, we conclude that
\begin{align*}
Q_k(u,v) &\le 1 - e^{2w+2} \(1 - \frac{2w+2}{v}\)^k \\
&= 1 - \exp \left\{ \frac{2w+2}{v} \( v - k + O(w)\) 
  \right\} \\
&= 1 - \exp \left\{ \frac{-2uw+O(u+w^2+1)}{v} \right\} \\
&\le \frac{2uw+O(u+w^2+1)}{v} \ll \frac{(u+1)(w+1)^2}{k}. 
\end{align*}

To prove \eqref{Qlower}, let $a=0$ in \eqref{QkA}.  The first sum on $l$ in
\eqref{QkA} is empty, and 
$$
\pfrac{k-l+w+\lam}{k-l+\lam}^{k-l} = \( 1 + \frac{w}{k-l+\lam} \)^{k-l} \le
e^w. 
$$
Thus, taking $A=e^w$ in \eqref{QkA}, we obtain
$$
Q_k(u,v) \ge 1 - e^w \( 1 - \frac{w}{v} \)^k.
$$
Now let $v=k+1-u$ (so that $w=1$) and let $f(x)=e^x (1- \frac{x}{k+1-u})^k$.
For $0\le x\le 1$, 
$$
f'(x) = - f(x) \frac{u-1+x}{k+1-u-x} \le 0;
$$
hence
\begin{align*}
f(1) &= 1 - \int_0^1 f(x) \ \frac{u-1+x}{k+1-u-x} \, dx \\
&\le 1 - \frac{f(1)}{k+1-u} \int_0^1 u-1+x\, dx =  1 - \frac{u-1/2}{k+1-u}
f(1). 
\end{align*}
This implies that $f(1) \le \frac{k+1-u}{k+1/2}$ and \eqref{Qlower} follows.
\hfill\qed

%
\section{The lower bound volume}\label{sec:lowvol}
%

In this section we prove Lemma \ref{volYnsv} using Lemma \ref{Q2}.  
We begin with a  crude upper bound on a combinatorial sum which will
also be needed for the upper bound integral in Section~\ref{sec:upperint}.

\begin{lem}\label{combsum} 
Let $0<\eps\le 1${\rm ,} $t\ge 10${\rm ,} and suppose $a,b$ are real numbers with 
$a+b > -(1-\eps) t$. 
Then
$$
\sum_{\substack{1\le j\le t-1 \\ -a < j < b+t}} \binom{t}{j} (a+j)^{j-1}
(b+t-j)^{t-j-1} \le e^{5+1/\eps} (t+a+b)^{t-1}.
$$
\end{lem}

\Proof 
Let $C_t(a,b)$ denote the sum in the lemma.  Since $C_t(a,b)=C_t(b,a)$,
we may suppose $a\le b.$  We also assume that $a>1-t$ and $b>1-t$,
otherwise $C_t(a,b)=0$.  The associated ``complete'' sum is evaluated
exactly using one of Abel's
identities (\cite[p.~20, Eq.~(20)]{Riordan})
\begin{equation}\label{Abel}
\sum_{j=0}^t \binom{t}{j} (a+j)^{j-1} (b+t-j)^{t-j-1} = \( \frac{1}{a} +
\frac{1}{b} \) (t+a+b)^{t-1} \qquad (ab\ne 0).
\end{equation} 
Note that $C_t(a,b)$ is increasing in $a$ and   in $b$.
If $-1 \le a\le b$, put $A=\max(1/2,a)$ and $B=\max(1/2,b)$.  By
\eqref{Abel},
\begin{equation}\label{Ctab}
\begin{split}
C_t(a,b) &\le C_t(A,B) \le \( \frac{1}{A} + \frac{1}{B} \) (t+A+B)^{t-1} \\
&\le 4 (t+a+b+3)^{t-1} \\
&\le 4e^{\frac{3(t-1)}{t+a+b}} (t+a+b)^{t-1} < e^5 (t+a+b)^{t-1}. 
\end{split}
\end{equation} 
Next assume $a< -1\le b$.  Since $(1-1/x)^x$ is an increasing function 
for $x>1$, when $j>-a$ we have
$$ 
(a+j)^{j-1}=(j-1)^{j-1} \( 1- \frac{-a-1}{j-1} \)^{j-1} \le (j-1)^{j-1}
\( 1 - \frac{-a-1}{t-1} \)^{t-1}.
$$ 
Also, $(a+1)b \le (1-\eps) t < \frac{1-\eps}{\eps}(t+a+b)$. 
Thus, by \eqref{Ctab},
\begin{align*}
C_t(a,b) &\le \pfrac{t+a}{t-1}^{t-1}  C_t(-1,b) \\
&\le e^5 \pfrac{(t+a)(t+b-1)}{t-1}^{t-1} = e^5 \(t+a+b + 
  \frac{(a+1)b}{t-1} \)^{t-1} \\
&< e^{4+1/\eps} (t+a+b)^{t-1}.
\end{align*}
Lastly, suppose $a \le b < -1$.  By \eqref{Ctab},
\begin{align*}
C_t(a,b) &\le \max_{-a<j<b+t} \pfrac{j+a}{j-1}^{j-1} \pfrac{b+t-j}{t-j-1}
^{t-j-1} C_t(-1,-1) \\
&< e^5 (t-2)^{t-1} \max_{-a<j<b+t} \pfrac{j+a}{j-1}^{j-1} 
\pfrac{b+t-j}{t-j-1}^{t-j-1}.
\end{align*}
As a function of real $j$, the maximum above occurs when
$\frac{j+a}{j-1}=\frac{b+t-j}{t-j-1}$; that is,
$j=\frac{(t-1)(a+1)+b+1}{a+b+2}$.
By our assumptions on $a$, $b$ and $t$, it is clear that $1<j<t-1$.
Hence $j+a$ and $b+t-j$ have the same sign, so that
$-a < j < b+t$.  Then
\begin{align*}
C_t(a,b) &\le e^5 \pfrac{t+a+b}{t-2}^{t-2} (t-2)^{t-1} \\
&< \eps^{-1} e^{5} (t+a+b)^{t-1} \\
&< e^{5+1/\eps} (t+a+b)^{t-1}.
\end{align*}
\vglue-22pt
\Endproof\vskip12pt


\emph{Proof of Lemma} \ref{volYnsv}.
For brevity, write
\begin{equation}\label{Snuv}
S_{k}(u,v) = \{ \bx: 0 \le \xi_1 \le \cdots \le \xi_k \le 1 : \xi_i \ge
\tfrac{i-u}{v} \, \forall i \},
\end{equation} 
so that $Q_k(u,v) = k! \Vol(S_k(u,v))$.  Let
\begin{equation}\label{Fxidef}
F_{k,v}(\bx) = \sum_{j=1}^k 2^{j-v \xi_j}.
\end{equation} 
Assume $v,k,s$ satisfy the hypotheses of Lemma \ref{volYnsv}, and
put  $u=k+1-v$.
For $1\le a\le k$, $0\le b\le k$, let 
\begin{align*}
V_1(a,b) &= \Vol\{ \bx\in S_k(u,v): \xi_a \le b/v \}, \\
V_2(a,b) &= \Vol \{ \bx\in S_k(u,v): \xi_{k+1-a} \ge 1-b/v \}.
\end{align*}
Make the change of variables $\theta_i=\xi_{a+i}$, so that $\theta_i \ge
\frac{i-(u-a)}{v}$ 
($1\le i\le k-a$).  By Lemma \ref{Q2}, we have
\begin{equation}\label{V1ab}
\begin{split}
V_1(a,b)&\le \Vol \{ 0 \le \xi_1\le \cdots \le \xi_a \le \tfrac{b}{v}\} 
  \Vol(S_{k-a}(u-a,v)) \\
&= \frac{(b/v)^a}{a!(k-a)!} Q_{k-a} (u-a,v) \\
&\le \frac{b^a}{a! k!} \pfrac{k}{v}^a  Q_{k-a}(u-a,v) \\
&\ll \frac{(100b)^a u}{a! k! (k-a+1)}.
\end{split}
\end{equation} 
Similarly,
\begin{equation}\label{V2ab}
\begin{split}
V_2(a,b)&\le \Vol(S_{k-a}(u,v)) \Vol \{ 1-\tfrac{b}{v} \le \xi_{k+1-a}
\le \cdots \le \xi_k \le 1\} \\
&\le \frac{(b/v)^a}{a!(k-a)!} Q_{k-a}(u,v) \\
&\ll \frac{a^2 (100b)^a u}{a! k! (k-a+1)}.
\end{split}
\end{equation} 

Next, we show that $F_{k,v}(\bx)$ is $\ll 2^u$ on average for 
$\bx\in S_k(u,v)$.
We integrate each term in the sum \eqref{Fxidef} separately, introducing 
$y=v\xi_j -j+u$, so that 
$$
\max(0,u-j)\le y\le u+v-j
$$
and
$$
0\le \xi_1 \le \cdots \le \xi_{j-1} \le \tfrac{j-u+y}{v} \le \xi_{j+1} \le
\cdots \le \xi_{k}\le 1.
$$
Let 
\begin{align*}
\theta_i &= \tfrac{v\xi_i}{j-u+y}\quad (1\le i\le j-1), \\
\zeta_i &= \tfrac{v}{u+v-j-y} \( \xi_{j+i} - \tfrac{j-u+y}{v}\) \quad 
(1\le i\le k-j).
\end{align*}
Then $\boldsymbol{\theta} \in S_{j-1}(u,j-u+y)$ and $\boldsymbol{\zeta} 
\in S_{k-j}(y,u+v-j-y)$.
By Lemma \ref{Q2}, we obtain
\begin{align*}
\int\limits_{S_k(u,v)} F_{k,v}(\bx)\, d\bx &= \frac{2^u}{v} \sum_{j=1}^k 
  \int_{\max(0,u-j)}^{u+v-j} 2^{-y} \pfrac{j-u+y}{v}^{j-1} 
  \frac{Q_{j-1}(u,j-u+y)}{(j-1)!} \\
&\qquad\qquad \times \pfrac{u+v-j-y}{v}^{k-j} \frac{Q_{k-j}(y,u+v-j-y)}{(k-j)!}
  \,dy \\
&\ll \frac{2^u u}{v^k (k+1)!} \int_0^{u+v-1} (y+1)^3 2^{-y} \sum_{\substack{
  1\le j\le k \\ u-y < j < u+v-y}} \binom{k+1}{j} \\
&\qquad\qquad \times (j+y-u)^{j-1} (u+v-y-j)^{k-j}\, dy.
\end{align*}
Since $k\le 100v-100$, $-u=v-1-k \ge -0.99 k$. 
By Lemma \ref{combsum} (with $t=k+1$, $a=y-u$, $b=-y$ and $\eps=0.01$),
for each $y$ the sum on $j$ on the right side 
is $\ll v^k$.  Since $\int_0^\infty
(y+1)^3 2^{-y}\, dy = O(1)$, we conclude that
\begin{equation}\label{intFineq}
\int\limits_{S_k(u,v)} F_{k,v}(\bx)\, d\bx \ll \frac{2^u u}{(k+1)!}.
\end{equation} 
Next,
\begin{eqnarray}\label{Y1}
\Vol(Y_k(s,v)) &\ge& \Vol \{ \bx\in S_k(u,v) : F_{k,v}(\bx) \le 2^s \} \\
&&- \sum_{1\le i\le \sqrt{k-M}} \bigl( V_1(M+i^2,i) + V_2(M+i^2,i) \bigr).
\nonumber
\end{eqnarray}
By \eqref{V1ab}, \eqref{V2ab} and the simple inequality $h! > (h/e)^h$,
\begin{equation}\label{sumV}
\begin{split}
\sum_{j=1}^2 \sum_{1\le i\le \sqrt{k-M}} V_j(M+i^2,i) &\ll \frac{u}{k!} 
\sum_{1\le i\le \sqrt{k-M}} \frac{(M+i^2)^2(100i)^{M+i^2}}
  {(M+i^2)! (k+1-M-i^2)} \\ 
&\ll \frac{u}{2^M \cdot (k+1)!}. 
\end{split}
\end{equation} 
Also,
$$
\Vol \{ \bx\in S_k(u,v) : F_{k,v}(\bx)\le 2^s \} \ge 
\Vol (S_k(u,v)) - 2^{-s} \int\limits_{S_k(u,v)} F_{k,v}(\bx)\, d\bx
$$
and, by Lemma \ref{Q2},
$$
\Vol(S_k(u,v)) \ge \frac{u}{3(k+1)!}.
$$
Combining this  with \eqref{intFineq}, \eqref{Y1} and \eqref{sumV}, we conclude that
$$
\Vol(Y_k(s,v)) \ge \frac{u}{(k+1)!} \( \frac14 - \frac{K_1}{2^M} - K_2 
2^{u-s} \),
$$
where $K_1$ and $K_2$ are positive absolute constants.  By hypothesis,
$u= k+1-v \le s-M/3$,
and the lemma follows provided $M$ is large enough.
\hfill\qed

%
%
%
\section{The upper bound integral}\label{sec:upperint}
%
%
%

The purpose of this section is to use the bounds on uniform order 
statistics proved in Section~\ref{sec:uos} to prove Lemma \ref{Unlem}.
The primary tools are bounds for the volumes of subsets of $S_k(u,v)$
(defined in \eqref{Snuv}) in which one or more coordinates is abnormally
small or large.

\begin{lem}\label{vol1}
Suppose $g,k,t,u,v \in \ZZ$ satisfy
$$
2\le g\le k/2, \; t\ge -1, \; v\ge k/10, \; u\ge 0, \; u+v\ge k+1.
$$
Let $R$ be the subset of $\bx \in S_k(u,v)$ where{\rm ,} for some $l\ge g+1${\rm ,}
\begin{equation}\label{lcond1}
\frac{l-u}{v} \le \xi_l \le \frac{l-u+1}{v}, \qquad \xi_{l-g} \ge 
\frac{l-u-t}{v}.
\end{equation} 
Then
$$
\Vol(R) \ll \frac{(10(t+1))^g}{(g-2)!} \, \frac{(u+1)(u+v-k)^2}{(k+1)!}.
$$
\end{lem}

\Proof  Fix $l$ satisfying $\max(u,g+1) \le l \le k$.  Let $R_l$
be the subset of  $\bx \in S_k(u,v)$ satisfying \eqref{lcond1} for
this particular $l$. 
We have $\Vol(R_l) \le V_1 V_2 V_3 V_4$, where by Lemma \ref{Q2},
\begin{align*}
V_1 &= \Vol\{ 0\le \xi_1 \le \cdots \le \xi_{l-g-1} \le \tfrac{l-u+1}{v}:
  \xi_i \ge \tfrac{i-u}{v}\, \forall i \} \\
&= \pfrac{l-u+1}{v}^{l-g-1} \frac{Q_{l-g-1}(u,l-u+1)}{(l-g-1)!} \\
&\ll \pfrac{l-u+1}{v}^{l-g-1} \frac{(u+1)g^2}{(l-g)!}, \\
V_2 &= \Vol \{ \tfrac{l-u-t}{v} \le \xi_{l-g} \le \cdots \le \xi_{l-1} \le
  \tfrac{l-u+1}{v} \} = \frac{1}{g!} \pfrac{t+1}{v}^g, \\
V_3 &= \Vol \{ \tfrac{l-u}{v} \le \xi_l \le \tfrac{l-u+1}{v} \} = 
  \frac{1}{v}, \\
V_4 &= \Vol\{ \xi_{l+1} \le \cdots \le \xi_k \le 1 : \xi_i \ge \tfrac{i-u}{v}
  \, \forall i \} \\
&= \frac{1}{(k-l)!} \pfrac{u+v-l}{v}^{k-l} Q_{k-l}(0,u+v-l) \\
&\ll \frac{(u+v-l)^{k-l-1} (u+v-k)^2}{v^{k-l} (k-l)!}.
\end{align*}
Thus
$$
\Vol(R) \! \ll \! \frac{(t+1)^g (u+1)(u+v-k)^2}{(g-2)! v^k (k-g)!}
\sum_l \! \binom{k-g}{l-g} (l-u+1)^{l-g-1} (u+v-l)^{k-l-1}.
$$
By Lemma \ref{combsum} (with $t=k-g$, $a=g+1-u$, $b=u+v-k$ and
$\eps=\frac{1}{10}$), the sum on $l$ is 
$$
\ll (v+1)^{k-g-1} \ll \frac{v^{k-g}}{k} \le 
\frac{v^k (k-g)!}{k\cdot k!} \, \pfrac{k}{v}^g \ll \frac{v^k 10^g (k-g)!}
{(k+1)!},
$$
and the lemma follows.
\Endproof\vskip4pt

To bound $U_k(v;\a)$, we will bound the volume of the set
\begin{equation}\label{TTUU}
\begin{split}
\TT(k,v,\g) &= \{ \bx\in \RR^k : 0 \le \xi_1 \le \cdots \le \xi_k \le 1, \\
&\qquad\qquad 
  2^{v\xi_1} + \cdots + 2^{v\xi_j} \ge 2^{j-\g}\; (1\le j\le k) \}.
\end{split}
\end{equation} 

\begin{lem}\label{UUlem}
Suppose $k,v,\g$ are integers with $1\le k\le 10v$ and $\g\ge 0$.  
Set $b=k-v$.  Then
$$
\Vol(\TT(k,v,\g)) \ll \frac{Y}{2^{2^{b-\g}} (k+1)!},
$$
where
$$
Y = \begin{cases} b & \text{ if } b\ge \g+5 \\ (\g+5-b)^2(\g+1) & 
\text{ if } b\le \g+4. \end{cases}
$$
\end{lem}

\Proof  
Let $r=\max(5,b-\g)$ and $\bx \in \TT(k,v,\g)$.  Then either
\begin{equation}\label{U-A1}
\xi_j > \tfrac{j-\g-r}{v} \quad (1\le j\le k)
\end{equation} 
or
\begin{multline}\label{U-A2}
\min_{1\le j\le k} ( \xi_j - \tfrac{j-\g}{v} ) = \xi_l - \tfrac{l-\g}{v} \in 
[\tfrac{-h}{v}, \tfrac{1-h}{v} ] \\ 
\text{ for some integers } h\ge r+1, 1\le l \le k.
\end{multline}
Let $V_1$ be the volume of
$\bx \in \TT(k,v,\g)$ satisfying \eqref{U-A1}.
If $b\ge \g+5$, \eqref{U-A1} is not possible, so $b\le \g+4$ and $r=5$.
By Lemma \ref{Q2},
$$
V_1 \le \frac{Q_k(\g+5,v)}{k!} \ll \frac{(\g+6)(\g+6-b)^2}{(k+1)!}
\ll \frac{Y}{2^{2^{b-\g}} (k+1)!}.
$$

If \eqref{U-A2} holds, then there is an integer $m$ satisfying
\begin{equation}\label{U-A3}
m\ge h-3, \; 2^m < \tfrac{l}2, \; \xi_{l-2^m} \ge \tfrac{l-\g-2m}{v}.
\end{equation} 
To see \eqref{U-A3}, suppose such an $m$ does not exist.  Then
\begin{align*}
2^{v\xi_1}+\cdots+2^{v\xi_l} &\le 2 \sum_{l/2<j\le l} 2^{v \xi_j} \\
&< 2 \biggl( 2^{h-3} 2^{l-\g-h+1} + \sum_{m\ge h-3} 2^m 2^{l-\g-2m} \biggr)\\
&\le 2^{l-\g},
\end{align*}
a contradiction.  
Let $V_2$  be the volume of
$\bx \in \TT(k,v,\g)$ satisfying \eqref{U-A2}.
Fix $h$ and $m$ satisfying \eqref{U-A3}
and apply Lemma \ref{vol1} with $u=\g+h$,
$g=2^m$, $t=2m$.  The volume of such $\bx$ is
\begin{align*}
&\ll \frac{(\g+h+1)(\g+h-b)^2}{(k+1)!}\, \frac{(20m+10)^{2^m}}{(2^m-2)!} \\
&\ll \frac{(\g+h+1)(\g+h-b)^2}{2^{2^{m+3}} (k+1)!}.
\end{align*}
The sum of $2^{-2^{m+3}}$ over $m\ge h-3$ is $\ll 2^{-2^{h}}$.
Summing over $h\ge r+1$ gives
\vskip12pt\hfill $
\displaystyle{V_2 \ll \frac{(\g+r+2)(\g-b+r+2)^2}{2^{2^{r+1}} (k+1)!}
\ll \frac{Y}{2^{2^{b-\g}} (k+1)!}.}
$ 
\Endproof

\vskip12pt 
\emph{Proof of Lemma} \ref{Unlem}.
Assume $k\ge 1$, since the lemma is trivial when $k=0$.
Put $b=k-v$ and define
$$
F(\bx) = \min_{0\le j\le k} 2^{-j}\( 2^{v\xi_1}+\cdots+2^{v\xi_j}+\a \).
$$
Let $t=\fl{\frac{\log \a}{\log 2}} \le 0$.  For integers $m\ge 0$, 
consider $\bx \in R_k$ satisfying
$2^{-m}\a \le F(\bx) < 2^{1-m}\a.$  For $1\le j\le k$, 
$$
2^{-j} \( 2^{v\xi_1}+\cdots+2^{v\xi_j} \) \ge
\max(2^{-j}, (2^{-m}-2^{-j})\a) \ge 2^{t-m-1},
$$
so that $\bx \in \TT(k,v,m+1-t)$.  Hence, by Lemma \ref{UUlem},
\begin{align*}
U_k(v;\a)  &\le \sum_{m\ge 0} 2^{1-m} \a \Vol(\TT(k,v,m+1-t)) \\
&\ll \frac{\a}{(k+1)!} \sum_{m\ge 0} \frac{2^{-m} Y_m}
  {2^{2^{b+t-m-1}}}, \\
Y_m &= \begin{cases} b & \text{ if } m\le b+t-6 \\ (m+6-t-b)^2(m+2-t) & 
\text{ if } m\ge b+t-5 .
  \end{cases}
\end{align*}
Next,
\begin{eqnarray*}
\sum_{m\ge 0} \frac{2^{-m} Y_m}{2^{2^{b+t-m-1}}}& = &
\sum_{0\le m\le b+t-6} \frac{b}{2^m 2^{2^{b+t-m-1}}} \\
&&+  \sum_{m\ge  \max(0,b+t-5)} \frac{(m+6-t-b)^2(m+2-t)}{2^m}.
\end{eqnarray*}
If $b\ge 6-t$, each sum on the right side is $\ll b 2^{-b-t}$.
If $b\le 5-t$, the first sum is empty and the second is $\ll (6-t-b)^2(2-t)$.
In both cases 
$$
\sum_{m\ge 0} \frac{2^{-m} Y_m}{2^{2^{b+t-m-1}}} \ll \frac{(1+|b+t|^2) (1-t)}
{2^{b+t}+1},
$$
whence
$$
U_k(v;\a) \ll \frac{\a (1+|v-k-\frac{\log \a}{\log 2}|^2) \log(2/\a)}
{(k+1)! (\a 2^{k-v}+1)}.
$$
Sometimes this is worse than the simpler bound
\begin{equation}\label{Un2}
U_k(v;\a) \ll \frac{\a}{k! (\a 2^{k-v}+1)},
\end{equation} 
which we now prove.  
When $\a 2^{k-v} \le 1$,  \eqref{Un2} follows from the trivial bound
$U_k(v;\a) \le \a / k!$.
Otherwise, 
\begin{align*}
U_k(v;\a) &\le 2^{1-k} \int\limits_{R_k}
  2^{v\xi_1}+\cdots+2^{v\xi_k}\, d\bx \\
&=2^{1-k} \sum_{j=1}^k \int_0^1 2^{vy} \frac{y^{j-1}
  (1-y)^{k-j}}{(j-1)! (k-j)!}\, dy \\
&= \frac{2^{1-k} (2^v-1)}{v (k-1)!} \le \frac{20}{2^{k-v} k!}.
\end{align*}
\vskip-22pt
\Endproof
 \vskip8pt

%
%
\section{Divisors of shifted primes}\label{sec:shifted}
%

\emph{Proof of  Theorem} \ref{thm:shifted}.
We first take care of some easy cases.   When $z\le z_0(y)$,
$H(x,y,z) \asymp \eta x$.  By the Brun-Titchmarsh inequality for primes in
arithmetic progressions (e.g.\ \cite{MV}),
$$
H(x,y,z;P_\lam) \ll \frac{x}{\log x} \sum_{y<d\le z} \frac{1}{\phi(d)}.
$$
The theorem in this case follows from the asymptotic formula \eqref{SR}.

If $z\ge y^{1.001}$, then $H(x,y,z) \asymp x$ and
$H(x,y,z;P_\lam) \le \pi(x) \ll \frac{x}{\log x}$.
The most difficult case is when $z_0(y) \le z \le y^{1.001}$.
We follow the outline from Section~\ref{sec:upper},
inserting a sieve estimate at the appropriate point.  

\begin{lem}\label{lem:shifted}
Suppose $10 \le y<z=e^{\eta} y${\rm ,} $y\le \sqrt{x}${\rm ,} 
$(\log y)^{-1/2} \le \eta \le \frac{\log y}{1000}$.  Then
$$
H(x,y,z;P_\lam) \ll x \max_{y^{1/2} \le t \le x} \; \sum_{P^+(a) \le te^{\eta}}
\frac{L(a;\eta)/\phi(a)}{\log^2 \( \frac{t}{a}+P^+(a) \) 
\log\( \frac{x^{0.49}}{a}+P^+(a) \)}.
$$
\end{lem}

The proof of Lemma \ref{lem:shifted} is similar to Lemma \ref{Hupper},
and will be given at the end of this section.
By the proof of Lemma \ref{Hupper2} and the Cauchy-Schwarz
inequality, we obtain
\begin{equation}\label{HP}
H(x,y,z;P_\lam) \ll x (1+\eta) \max_{y^{1/2} \le t\le x} S^*(t;\eta)^{1/2} 
\widetilde{S}(t;\eta)^{1/2},
\end{equation} 
where
$$
\widetilde{S}(t;\eta) =  \sum_{a \in \PP^*(e^\eta,te^{\eta})} 
\frac{a}{\phi^2(a)}\, \frac{L(a;\eta)}{\log^2 (x^{0.4}/a+P^+(a))}.
$$
In the same way that \eqref{ST} was proved, we obtain
$$
\widetilde{S}(t;\eta) \ll \frac{\widetilde{T}(\eta,te^{\eta},1)}{\log^2 x} + 
\sum_{\substack{
  k\in\ZZ, k\ge 1 \\ e^{\eta} \le e^{e^{k-1}} \le te^{\eta}}}
  e^{-2k} \widetilde{T}(\sg,e^{e^k},x^{1/10}),
$$
where
$$
\widetilde{T}(\sg,P,Q) = \sum_{\substack{a\in\PP^*(e^{\sg},P) \\ a\ge Q}}
\frac{a}{\phi^2(a)} \, L(a;\sg).
$$
Similarly, we define $\widetilde{T}_k(\sg,P,Q)$.
The bound given in Lemma \ref{Tnint} holds with
$T_k(\sg,P,Q)$ replaced by $\widetilde{T}_k(\sg,P,Q)$.   
The only change in the proof is to define $E_j$ so that
$$
\sum_{p \in E_j} \( \frac{p}{(p-1)^2} \)^{1-\g} \le \log 2.
$$
Therefore, Lemma \ref{Sstarfinal} holds with $S^*(t;\eta)$ replaced by
$\widetilde{S}(t;\eta)$, and the theorem follows from \eqref{HP}.\hfill
\qed

\demo{Proof of Lemma {\rm \ref{lem:shifted}}} 
Suppose $q\le x$, $q$ is prime and $\tau(q+\lam,y,z)\ge 1$.
We distinguish two cases.  First, assume $z \le \exp\{(\log x)^{9/10}\}$.
Let $y<d\le z$, $d|(q+\lam)$ and set $p=P^+(d)$.  
Write $q+\lam=apb$, where $P^+(a) \le p \le  P^-(b)$.  The number of $q$
with $ap>\sqrt{x}$ is, by Lemma \ref{Psixyz}, at most
$$
\sum_{\substack{m>\sqrt{x} \\ P^+(m)\le z}} \frac{x}{m} \ll 
\frac{x\log x}{\log z} e^{-\frac{\log x}{8\log z}} \ll \frac{x}{\log^{10} z}. 
$$
Now suppose $ap \le \sqrt{x}$.  Given $a,p$, we wish to count 
the number of $b\le x/(ap)$ with $P^-(b) \ge p$ and
$abp-\lam$ prime.  By the arithmetic form of the large sieve inequality
\cite{MV}, the number of $b$ is
$$
\ll \frac{x}{ \phi(ap) \log x \log p} \ll \frac{x}{\log x} \, \frac{1}
{\phi(a) p\log p}.
$$
As in the proof of Lemma \ref{Hshort}, for each $a$ we have
$$
\sum_{\substack{\tau(a,y/p,z/p) \ge 1 \\ p \ge P^+(a)}} \frac{1}{p\log p}
\ll \frac{L(a;\eta)}{\log^2(y/a+P^+(a))}
$$
and Lemma \ref{lem:shifted} follows in this case.

Now suppose $z>\exp \{ (\log x)^{9/10} \}$.  Put $\sg=\eta$ and suppose 
$x_1,x_2$ satisfy the hypotheses of Lemma \ref{Hshort}.  Define $\AA$
as in that lemma, and let $\AA^*$ be the set of $n\in\AA$ for which
$n-\lam$ is prime.   By the argument leading
to \eqref{HHr1}, and since $\log z \ge (\log x)^{9/10}$, we obtain
$$
H(x_2,y,z;P_\lam) - H(x_1,y,z;P_\lam) \le |\AA^*| + O\pfrac{x_2-x_1}
{\log^{3.6} x_1}.
$$
Assume \eqref{pdef} and write $q+\lam$ in the form \eqref{nrep1}.
Given $a$ and $p$, we wish to count the number of $b\in(\frac{x_1}{ap},
\frac{x_2}{ap}]$ with $P^-(b) > p$, $b>p$, and $apb-\lam$ prime.
Here
$$
\frac{x_2-x_1}{ap} \ge \frac{x_1}{ap\log^{10} z} \ge \max\(p^{1/2},
\frac{x_1^{0.999}}{az_j}\) := Q.
$$
We apply the arithmetic form of the large sieve \cite{MV}, eliminating those
$b$ with $b\equiv 0\pmod{\varpi}$ for some prime $\varpi < p^{1/10}$
and those $b$ with $apb-\lam\equiv 0\pmod{\varpi}$ for some prime
$\varpi \le Q^{1/3}$.  The number of remaining $b$ is
$$
\ll \frac{x_2-x_1}{\phi(ap) \log p \log Q} \ll \frac{x_2-x_1}{\phi(a)
p\log p \log (x_1^{0.999}/(az_j)+y_j/a+P^+(a))}.
$$
As in the proof of Lemma \ref{Hshort}, the sum over $p$ of $\frac{1}{p\log p}$
is
$$
\ll \frac{L(a;\eta)}{\log^2 (y_j/a+P^+(a))}.
$$
Thus
\begin{multline*}
H(x_2,y,z;P_\lam) - H(x_1,y,z;P_\lam) \le (x_2-x_1) \biggl[
   O\pfrac{1}{\log^{3.6} x}   \\
+ \sum_{j=1}^2 \sum_{P^+(a) \le z_j}
  \frac{L(a;\eta)}{\phi(a) \log^2 (y_j/a+P^+(a))
   \log(x_1^{0.999}/az_j + y_j/a + P^+(a))} \biggr].
$$
\end{multline*}
Since $z_j \le y_j^{1.001} \le x^{0.001} y_j$, we have
$$
\frac{x_1^{0.999}}{z_j} + y_j \ge x_1^{0.499}. 
$$ 
By considering the term $a=1$ alone, we see that
the double sum exceeds 
$$
\eta/\log^3 y_1 \ge (\log x_1)^{-3.6}.
$$
Lastly, we sum over short intervals $(x_1,x_2] \subseteq [x/\log^5 x, x]$.
\Endproof
\vskip4pt 

\emph{Proof of Theorem} \ref{thm:shiftedlong}.
First, suppose $0<\a<\b<\frac12$.
Let $N_p$ denote the number of primes $q\in (x/2,x]$ with $q\equiv
 -\lam \pmod{p}$.  For large $x$, a given number $q+\lam$ can be
 divisible by at most $\fl{1/\a}+1$ primes $p>x^{\a}$.  Thus, using
 the Bombieri-Vinogradov theorem (see e.g., Chapter 28 of \cite{Dav}), we have
\begin{equation}\label{Hlab}
\begin{split}
H(x,x^\a,x^\b;P_\lam) - H(x/2, x^\a,x^\b;P_\lam) &\ge
\frac{1}{\fl{1/\a}+1} \sum_{x^\a<p\le x^\b} N_p \\
&\gg_{\a,\lam} \frac{x}{\log x} \sum_{x^\a<p\le x^\b} \frac{1}{p-1} \\
&\gg_{\a,\b,\lam} \frac{x}{\log x}.
\end{split}
\end{equation} 
The theorem follows when $0\le a <b \le \frac12$ by applying
\eqref{Hlab} with $\a,\b$ chosen so that $a<\a<\b<b$.
When $0 \le a < \frac12 < b$, the theorem follows by applying
\eqref{Hlab} with $a<\a<\frac12$ and $\b=\frac12$.
Finally, suppose $\frac12 \le a < b \le 1$ and let $a',b'$ be numbers
with $a<a'<b'<b$.  For large $x$, $q\in (x/2,x]$ implies 
$q+\lam \in (x/4,2x]$.  Thus, if $d|(q+\lam)$ and $x^{1-b'}<d\le
  x^{1-a'}$, then the complementary divisor $\frac{q+\lam}{d}$ lies in
$(\frac14 x^{a'}, 2x^{b'}]$.  For large $x$, $\frac{q+\lam}{d}\in
  (x^a,x^b]$ and thus the theorem follows by taking $\a=1-b'$ and
$\b=1-a'$ in \eqref{Hlab}.
\Endproof

\end{document}